\documentclass[11pt,letterpaper]{article}

\usepackage{MathStyle}
\usepackage[long]{PaperStyle} 

\usepackage{bm}
\usepackage{needspace}

\newcommand{\Dvp}{\mathsf{D}}

\theoremstyle{plain}
\newcommand{\restatedtheoremheading}{}
\newtheorem*{restatedtheoreminner}{\restatedtheoremheading}
\newenvironment{restatedtheorem}[2][]{%
    \renewcommand{\restatedtheoremheading}{Theorem~\ref{#2}}%
    \begin{restatedtheoreminner}[#1]%
}{%
    \end{restatedtheoreminner}%
}

\title{Davenport Constants of High-Dimensional Integer Sets}
\author{Stefan Kuhlmann\textsuperscript{1} \and
        Timm Oertel\textsuperscript{2} \and
        Nicola Wengiel\textsuperscript{2}}
\affiliation{\textsuperscript{1}Department of Mathematics, ETH Zürich\\
             \textsuperscript{2}Department of Mathematics, FAU Erlangen-Nürnberg}
\emails{\emailaddress{stefan.kuhlmann@math.ethz.ch}\quad
        \emailaddress{timm.oertel@fau.de}\quad
        \emailaddress{nicola.wengiel@fau.de}}
\hypersetup{pdftitle={Davenport Constants of High-Dimensional Integer Sets},
            pdfauthor={Stefan Kuhlmann, Timm Oertel, Nicola Wengiel}}

\begin{document}

\maketitle

\begin{abstract}
    The Davenport constant of a finite set of integer vectors is the maximum length of an irreducible zero-sum sequence over that set.
    For nonempty sets closed under negation, it equals the largest $\ell_1$-norm of a Graver basis element of the matrix formed by those vectors.
    We study extremal and average-case versions in growing dimension, with applications to integer programming.
    For integer cubes $[-R,R]^d \cap \Z^d$ and the sign cube $\{-1,1\}^d$, totally unimodular encodings improve on the upper bounds from direct circuit decomposition, while random perturbations of Hadamard matrices give new lower bounds.
    A key ingredient is a uniform affine Bernoulli variant of Nguyen and Wood's surjectivity theorem.
    In every sufficiently large dimension the sign-cube lower bound is $d^{d/2}2^{-d-o(d)}$, and the integer-cube lower bound is $R^d d^{d/2}2^{-o(d)}$ uniformly for polynomially bounded integer radii.
    For random sign-cube sequences of even length conditioned on summing to zero, we prove a sharp threshold from irreducibility to reducibility at $d\log d/\log 4$, far below the extremal scale.
    For fixed $p$ and $R$, we determine exact exponential rates for lattice $\ell_p$-balls when $p=1$, $R\geq3$, or $2\leq p<\infty$, $\lfloor R^p\rfloor\geq3$.
    For each fixed integer $k\geq3$, we determine the exponential rate of the worst-case $\ell_1$-distance from an optimal LP vertex to the nearest integer optimum in standard-form programs with ternary matrices having at most $k$ nonzero entries per column.

\end{abstract}

\makepapercontents

\section{Introduction}
\label{sec:Introduction}

The Davenport constant has roots in zero-sum theory over finite abelian groups and non-unique factorization: prime-ideal factorizations of irreducible algebraic integers give irreducible zero-sum sequences in their ideal class groups; see the survey \cite{Gao2006}. We study the corresponding invariant for finite sets of integer vectors.
For a finite $\calX \subseteq \Z^d$, a nonempty zero-sum sequence over $\calX$ is \emph{irreducible} if no nonempty proper subsequence sums to zero. The Davenport constant $\Dvp(\calX)$ is the maximum length of an irreducible sequence over $\calX$, counting repetitions, and is zero if no such sequence exists.

We focus on integer cubes $\calQ_d(R) \coloneqq [-R,R]^d \cap \Z^d$, lattice $\ell_p$-balls, and the sign cube $\calS_d \coloneqq \{-1,1\}^d$. This choice is motivated by indecomposable integer relations in integer programming and by the geometric questions of Baeth and Geroldinger \cite[p.~315, B.2(a)--(b)]{Baeth2014}. Our main regime is growing dimension at fixed radius, with cube bounds uniform for polynomially growing integer radii.
We combine integral rounding for cube upper bounds with random perturbations of structured matrices for lower bounds. We also determine a sharp threshold from irreducibility to reducibility for random sequences in the sign cube conditioned on having sum zero.

\paragraph{Integer programming and coin weighing.}
Irreducible zero-sum sequences encode indecomposable nonnegative integer dependencies among matrix columns. They correspond to Hilbert basis elements and to nonnegative Graver basis elements, which underpin augmentation methods in integer programming \cite{DeLoera2008,Onn2010,Eisenbrand2019,Brianski2024}.
For the ternary cube $\{-1,0,1\}^d$, Alon, Kozlov, and V\~u \cite[Section~5]{Alon1996} relate $\Dvp(\{-1,0,1\}^d)$ to testing whether coins of two possible weights all have equal weight.

\paragraph{Integer cubes and probabilistic constructions.}
For $d \geq 3$, Plagne and Tringali \cite[Theorems~4--5]{Plagne2015} give
\begin{equation*}
    (2R - 1)^d \leq \Dvp(\calQ_d(R)) \leq (2(d - 1 + 1/d)R + 1)^d
\end{equation*}
for integer $R \geq 2$. For $R = 1$, the lower bound is $\Dvp(\calQ_d(R)) \geq 2^d$, and the same upper bound remains valid.
The upper bound uses the Steinitz lemma \cite{Steinitz1913,Grinberg1980}. The same method gives the related Graver bound $(2dR+1)^d$ \cite[Lemma~27]{Eisenbrand2019}, while circuit decomposition gives $\Dvp(\calQ_d(R)) \leq (2R+1)^d R^d(d+1)^{(d+1)/2}$ (see \cref{sec:Preliminaries,sec:IntegerCubes}). For comparison, \cite[Lemma~3]{Hunkenschroeder2026} gives a Graver element of $\ell_1$-norm $(R^{d+1}-1)/(R-1)$ for a $d$-row matrix with entries in $[-R,R]$, for every integer $R\geq2$.

For the sign cube, Alon and V\~u's anti-Hadamard construction \cite[Proposition~3.4.3]{Alon1997} already yields the superexponential lower bound $\Dvp(\calS_d) \geq d^{d/2}2^{-2d-o(d)}$ in all sufficiently large dimensions; see also \cite{Buchin2012}. For $d = 2^m$, their stronger construction from Theorem~2.1.1 gives $d^{d/2}2^{-d-o(d)}$, matching our lower bound for the sign cube up to subexponential factors.
For $d=2^m$, an adaptation of their construction to integer cubes gives the deterministic lower bound $\Dvp(\calQ_d(R)) \geq R^d d^{d/2}2^{-o(d)}$ uniformly over all integers $R \geq 1$. Our probabilistic construction gives the lower bounds for both cubes in every sufficiently large dimension. The resulting circuits give large Graver coefficients.

The construction uses \emph{circuits}, primitive integer kernel vectors of minimal support. After restricting to that support and selecting a maximal set of linearly independent rows, their coefficients are signed maximal minors divided by their common gcd. A large determinant therefore does not suffice unless this divisor is controlled. Alon and V\~u get around this issue by deterministically changing entries of a special Hadamard matrix that exists when $d = 2^m$.
We use a uniform affine Bernoulli variant of Nguyen and Wood's surjectivity theorem \cite{Nguyen2022}. Each entry of a polynomially bounded deterministic $n\times(n+1)$ integer matrix is independently shifted by a prescribed $+1$ or $-1$ with probability $q_n=n^{-1/6+\eps}$ for fixed $0<\eps<1/6$. The limiting probability of surjectivity over $\Z$ is $\prod_{j=2}^\infty\zeta(j)^{-1}>0$, uniformly in the deterministic data. See \cref{thm:AffineBernoulliSurjectivity,app:AffineBernoulliSurjectivity}.

For integer cubes, we perturb a scaled Hadamard matrix and adjoin a Bernoulli column. Surjectivity implies that the maximal minors have greatest common divisor one. For the sign cube, the perturbation keeps all entries in $\{-1,1\}$ and surjectivity of an auxiliary matrix gives maximal-minor gcd $2^{h-1}$ in Hadamard dimension $h$. Matrix concentration preserves a large determinant with probability $1-o(1)$ \cite{Tropp2012}, so the determinant and gcd conditions hold together with probability bounded away from zero along Hadamard orders. Nearby Hadamard orders extend the bounds to all sufficiently large dimensions \cite{Brent2013}.

For the upper bounds, we split coordinates into two blocks and encode vectors as edges of a bipartite graph. Total unimodularity and a lemma of Paat, Weismantel, and Weltge \cite{Paat2020} allow rounding circuit combinations while preserving sums and lengths. This encoding also bounds integrality numbers and affine TU-dimensions \cite{Paat2022,Bader2018}.

\begin{samepage}
    \begin{theorem}[Davenport constants of integer cubes]
        \label{thm:IntegerCubes}
        For every fixed $T > 0$, as $d \to \infty$, uniformly over integers $1 \leq R \leq d^T$,
        \begin{equation*}
            2^{-o(d)}R^d d^{d/2} \leq \Dvp(\calQ_d(R)) \leq \paren*{(2R+1)^{\floor*{d/2}} + (2R+1)^{\ceil*{d/2}}}R^d(d+1)^{(d+1)/2}.
        \end{equation*}
        For the sign cube,
        \begin{equation*}
            2^{-o(d)}\frac{d^{d/2}}{2^{d-1}} \leq \Dvp(\calS_d) \leq \paren*{2^{\floor*{d/2}} + 2^{\ceil*{d/2}}}\frac{(d+1)^{(d+1)/2}}{2^{d-1}}.
        \end{equation*}
    \end{theorem}
\end{samepage}

For the sign cube, this improves the circuit bound $2(d+1)^{(d+1)/2}$ exponentially, giving $\log_2{\Dvp(\calS_d)} \leq \tfrac{d}{2}\log_2{d} - \tfrac{d}{2} + O(\log{d})$.
A dimension-independent LP-to-IP $\ell_\infty$-proximity bound linear in the largest absolute square subdeterminant \cite[after Conjecture~1]{Paat2020} would make both lower bounds sharp up to subexponential factors.
See \cref{rem:SteinitzCubeUpper} for a separate conditional Steinitz bound.

\paragraph{A sharp average-case threshold.}
We study exact integer sums conditioned on a zero total. Related random zero-sum problems in finite abelian groups appear in \cite{Balachandran2021,Bushaw2024}. Let $\nu_d$ be the uniform distribution on $\calS_d$. For even $n \geq 2$, sample $n$ independent vectors from $\nu_d$ and condition on their sum being zero. Let $p_d(n)$ be the conditional probability of a nonempty proper zero-sum subsequence.
For $\eps \in (0,1)$, the average-case Davenport constant $D_\eps(\nu_d)$ is the first even length $n \geq 2$ with $p_d(n) \geq 1-\eps$.

\begin{theorem}[Average-case Davenport constant of the sign cube]
    \label{thm:AverageCase}
    For every fixed $\eps \in (0,1)$,
    \begin{equation*}
        \lim_{d \to \infty} \frac{D_\eps(\nu_d)}{d\log{d}} = \frac{1}{\log{4}}.
    \end{equation*}
    More precisely, for every fixed $\eta \in (0,1)$ and every sequence of even integers $n = n(d) \geq 2$,
    \begin{equation*}
        \lim_{d \to \infty} p_d(n(d)) =
        \begin{cases}
            0, & \displaystyle n(d) \leq (1 - \eta)\frac{d\log{d}}{\log{4}}, \\[6pt]
            1, & \displaystyle n(d) \geq (1 + \eta)\frac{d\log{d}}{\log{4}},
        \end{cases}
    \end{equation*}
    whenever the corresponding inequality holds for all sufficiently large $d$.
\end{theorem}

This threshold lies far below the extremal scale $\exp\paren*{\tfrac{d}{2}\log{d} + O(d)}$.
The proof uses first and second moments, controlling regular overlaps by local approximations and concentration and boundary overlaps by entropy estimates.

\paragraph{Exact exponential rates for lattice balls.}
Girard and Plagne study low-dimensional balls \cite{Girard2025} and boxes with a prescribed number of distinct vector types \cite{Girard2026}. We fix $p$ and $R$ and let the dimension grow, without restricting the number of vector types. For $1 \leq p < \infty$ and $R \geq 0$, let $\calB_{p,d}(R) \coloneqq \{v \in \Z^d : \norm*{v}_p \leq R\}$. For integers $k \geq 3$, put
\begin{equation*}
    \rho_k \coloneqq \sqrt{\frac{k + \sqrt{k^2 - 4}}{2}}.
\end{equation*}

\begin{theorem}[Davenport constants of lattice balls]
    \label{thm:LatticeBalls}
    For every fixed $2 \leq p < \infty$ and $R \geq 0$ with $k \coloneqq \floor*{R^p} \geq 3$,
    \begin{equation*}
        \lim_{d \to \infty} \Dvp(\calB_{p,d}(R))^{1/d} = \rho_k.
    \end{equation*}
    For every fixed $R \geq 3$,
    \begin{equation*}
        \lim_{d \to \infty} \Dvp(\calB_{1,d}(R))^{1/d} = \floor*{R} - 1.
    \end{equation*}
\end{theorem}

For integers $k \geq 0$, define the sparse ternary cubes $\calT_d(k) \coloneqq \{v \in \calQ_d(1) : \abs*{\supp(v)} \leq k\}$. For $k \geq 3$, integrality and $2^k > k$ give $\calT_d(k) = \calB_{k,d}(k^{1/k})$. Applying \cref{thm:LatticeBalls} yields the following corollary.

\begin{corollary}[Davenport constants of sparse ternary cubes]
    \label{cor:SparseTernaryCubes}
    For every fixed integer $k \geq 3$,
    \begin{equation*}
        \lim_{d \to \infty} \Dvp(\calT_d(k))^{1/d} = \rho_k.
    \end{equation*}
\end{corollary}

Here deterministic constructions match a Gram-determinant bound for $p \geq 2$ and a graph-based bound for $p = 1$. For fixed parameters the number of allowed integer vectors is polynomial in $d$, so the transfer from circuits preserves the exponential rate.
The corollary also applies to fixed all-equal weighing plans with two possible weights and at most $k$ weighings per coin \cite{Alon1996}.

\paragraph{Proximity in integer programming.}
The sparse ternary rate also controls the distance between LP and integer optima. For a finite column set $\calX \subseteq \Z^d$, let $\Pi(\calX)$ be the supremum of the $\ell_1$-distance from an optimal LP vertex to a nearest optimal integer solution. The supremum ranges over all optimal LP vertices of standard-form programs $\min\{c^\top x : Ax=b,\ x\geq0\}$ with columns in $\calX$, integral $b,c$, a finite LP optimum, and an integer feasible point.
Eisenbrand and Weismantel use LP-to-IP proximity to narrow the search for an optimal integer solution in their dynamic programming algorithm for integer programs with upper bounds on the variables \cite[Section~4]{Eisenbrand2020}.

\begin{corollary}[Proximity for column-sparse ternary matrices]
    \label{cor:SparseTernaryProximity}
    For every fixed integer $k\geq3$,
    \begin{equation*}
        \lim_{d\to\infty}\Pi(\calT_d(k))^{1/d}=\rho_k.
    \end{equation*}
\end{corollary}

The proof is at the end of \cref{sec:LatticeBalls}.
After collecting the common circuit tools in \cref{sec:Preliminaries}, we prove the cube bounds in \cref{sec:IntegerCubes}, the average-case threshold in \cref{sec:AverageCase}, and the lattice-ball results in \cref{sec:LatticeBalls}. The longer supporting proofs are given in the appendices.

\section{Circuits and Zero-Sum Sequences}
\label{sec:Preliminaries}

We collect the circuit tools used for cubes and lattice balls. These are standard in discrete optimization; see \cite{Onn2010}. Logarithms are natural unless otherwise specified.

We say that a vector $y$ \emph{conforms to} $x$ if $\supp(y) \subseteq \supp(x)$ and $x_i y_i > 0$ for every $i \in \supp(y)$. In particular, $y_i = 0$ whenever $x_i = 0$. A decomposition is \emph{conformal} if each summand conforms to the sum.

For a finite set $\calX = \{a_1,\ldots,a_s\} \subseteq \Z^d$, let $A_{\calX} \coloneqq [a_1 \ \cdots \ a_s]$. A nonempty zero-sum sequence over $\calX$ is represented by a nonzero multiplicity vector $x \in \Z_{\geq 0}^s$ satisfying $A_{\calX}x = 0$, and its length is $\norm*{x}_1$. It is irreducible precisely when $x$ cannot be written as the sum of two nonzero nonnegative integer kernel vectors. The multiplicity vectors of irreducible zero-sum sequences over $\calX$ therefore form the Hilbert basis of $\{v \in \R_{\geq 0}^s : A_{\calX}v = 0\}$. More generally, a Graver basis element is a nonzero integer kernel vector admitting no decomposition into two nonzero integer kernel vectors that both conform to it. For nonempty $\calX = -\calX$, orienting columns according to coefficient signs identifies the largest Graver $\ell_1$-norm with $\Dvp(\calX)$.

\paragraph{Circuits and cofactors.}
For an integer matrix $M \in \Z^{d \times n}$, a \emph{circuit} is a nonzero vector $c \in \ker M \cap \Z^n$ whose support is inclusion-wise minimal and whose coordinates have greatest common divisor one. Its support indexes a minimally linearly dependent set of columns. We call such a column configuration a \emph{positive circuit} if its circuit relation can be chosen positive on that support. A zero column gives a singleton circuit of $\ell_1$-norm one. Circuits are Graver basis elements.

Put $r \coloneqq \rank M$ and retain $r$ linearly independent coordinate rows to obtain $\widehat{M}$, which has the same kernel as $M$. Extend the support of a circuit $c$, if necessary, to a column index set $I = \{i_1,\ldots,i_{r+1}\}$ such that $\widehat{M}_I$ has rank $r$. Then, up to a common sign,
\begin{equation}
    \label{eq:CircuitCofactors}
    c_{i_j} = \frac{(-1)^{j-1}\det \widehat{M}^{(j)}_I}{D} \quad (1 \leq j \leq r+1), \qquad
    D = \gcd\set*{\abs*{\det \widehat{M}^{(k)}_I} : 1 \leq k \leq r+1},
\end{equation}
with $c_i = 0$ outside $I$. Here $\widehat{M}^{(j)}_I$ is obtained by deleting column $j$ of $\widehat{M}_I$. The denominator is the gcd of the maximal minors of this selected matrix. The formula also covers $r = 0$ with the convention that the determinant of the empty matrix is one.

\paragraph{From circuits to irreducible sequences.}
The following standard transfer bound uses the actual rank of the occurring vectors; compare \cite[Lemma~3.20]{Onn2010}.
We include the proof for completeness.

\begin{lemma}[From circuits to irreducible sequences]
    \label{lem:CircuitTransfer}
    Let $a_1,\ldots,a_s \in \Z^d$ be the distinct vectors occurring in an irreducible zero-sum sequence of length $L$, and let $B \coloneqq [a_1 \ \cdots \ a_s]$ have rank $r$. If every primitive nonnegative integer circuit $g$ of $B$ satisfies $\norm*{g}_1 \leq F$, then
    \begin{equation*}
        L \leq (s-r)F.
    \end{equation*}
\end{lemma}

\begin{proof}
    Let $x \in \Z_{>0}^s$ be the multiplicity vector. The cone $\{v \in \R_{\geq 0}^s : Bv = 0\}$ is generated by primitive nonnegative integer circuits and lies in the $(s-r)$-dimensional space $\ker B$. Conic Carath\'eodory therefore gives
    \begin{equation*}
        x = \sum_{j=1}^t \alpha_j g_j, \qquad \alpha_j \geq 0, \qquad t \leq s-r.
    \end{equation*}
    If some $\alpha_j \geq 1$, then $0 \neq g_j \leq x$, so irreducibility forces $g_j = x$. Hence $L \leq F \leq (s-r)F$. Otherwise all coefficients are less than one, and
    \begin{equation*}
        L = \norm*{x}_1 = \sum_{j=1}^t \alpha_j\norm*{g_j}_1 \leq tF \leq (s-r)F.
        \qedhere
    \end{equation*}
\end{proof}

For a finite $\calX\subseteq \Z^d$, define its maximum circuit $\ell_1$-norm by
\begin{equation}
    \label{def:Gamma}
    \gamma(\calX) \coloneqq \max_{c \text{ circuit of } A_{\calX}}\norm*{c}_1,
\end{equation}
with value zero if there is no circuit.

\begin{corollary}[Circuit bounds for symmetric finite sets of integer vectors]
    \label{cor:CircuitDavenport}
    Let $\calX = -\calX \subseteq \Z^d$ be nonempty and finite. Put $r \coloneqq \rank A_{\calX}$ and $q \coloneqq \abs*{\calX}-r$. Then $q \geq 1$ and
    \begin{equation*}
        \gamma(\calX) \leq \Dvp(\calX) \leq q\gamma(\calX).
    \end{equation*}
\end{corollary}

\begin{proof}
    The set $\calX$ contains either zero or an opposite pair, so $q \geq 1$. Orienting the columns of any circuit $c$ according to its coefficient signs gives an irreducible zero-sum sequence over $\calX$ of length $\norm*{c}_1$. This proves the lower bound, including singleton circuits at zero. For the upper bound, apply \cref{lem:CircuitTransfer} to the distinct types of an irreducible sequence. If their column matrix has rank $r'$ and contains $s$ columns, then $r \leq r' + (\abs*{\calX}-s)$, so $s-r' \leq q$. Its circuits, padded with zeros, are circuits of $A_{\calX}$ and hence have $\ell_1$-norm at most $\gamma(\calX)$.
\end{proof}

\section{Davenport Constants of Integer Cubes}
\label{sec:IntegerCubes}
We prove \cref{thm:IntegerCubes} using the circuit parameter and transfer bounds from \cref{sec:Preliminaries}. Throughout this section, suppose that $\calX = -\calX \subseteq \Z^d$ is finite and spans $\R^d$. Then \cref{cor:CircuitDavenport} gives the factor $|\calX|-d$ in the circuit-to-Davenport comparison.

Applying Hadamard's inequality to a $(d+1)\times (d+1)$ matrix gives the following circuit bound. Here, $\gcd(\calX)$ is the greatest common divisor of all $d\times d$ subdeterminants of $A_{\calX}$ and therefore divides $D$ in \eqref{eq:CircuitCofactors}.

\begin{lemma}[{cf. \cite[Section~5]{Alon1996}, \cite[Corollary~3.5]{Kozlov1997}}]
	\label{lem:GeneralCircuitUpperBound}
	Let $\calX$ be as above and $R$ be the largest entry in absolute value of a vector in $\calX$. Then 
	\begin{equation*}
		\gamma(\calX)\leq \frac{R^d (d+1)^{(d+1)/2}}{\gcd(\calX)}.
	\end{equation*}
\end{lemma}

Combining \cref{cor:CircuitDavenport} and \cref{lem:GeneralCircuitUpperBound} gives the standard estimates
\begin{enumerate}
	\item $\Dvp(\mathcal{Q}_d(R))\leq (2R+1)^d R^d (d+1)^{(d+1)/2}$ and
	\item $\Dvp(\mathcal{S}_d)\leq 2^d \frac{(d+1)^{(d+1)/2}}{2^{d-1}} = 2(d+1)^{(d+1)/2}$,
\end{enumerate}
where we use $\gcd(\mathcal{Q}_d(R)) = 1$ and $\gcd(\mathcal{S}_d) = 2^{d-1}$, and estimate $|\calX|-d$ by $|\calX|$ because keeping it does not change the asymptotic dependence on $d$. 

Let $\gamma_{\infty}(\calX)$ denote the \emph{max-circuit imbalance measure} of $A_{\calX}$, that is, the maximum absolute coefficient of a primitive circuit \cite{Ekbatani2022}. It is obtained by replacing the $\ell_1$-norm with the $\ell_\infty$-norm in \eqref{def:Gamma}. Then $\gamma_{\infty}(\calX)\leq \gamma(\calX)\leq (d+1) \gamma_\infty(\calX)$, so the two parameters agree up to a subexponential factor in the regime considered here. The construction in \cref{sec:HadamardLowerBounds} provides a lower bound for $\gamma_{\infty}(\calX)$ itself.

\subsection{Lower Bounds}
\label{sec:HadamardLowerBounds}

We prove the lower bounds in Hadamard dimensions and then extend them to arbitrary sufficiently large dimensions. The resulting cube circuits attain the bound in \cref{lem:GeneralCircuitUpperBound} up to a subexponential factor. To control the gcd of their maximal minors, we use the following affine version of Nguyen and Wood's surjectivity theorem \cite{Nguyen2022}.

\begin{samepage}
\begin{theorem}[{Affine Bernoulli Surjectivity; cf. \cite[Theorem~1.1]{Nguyen2022}}]
	\label{thm:AffineBernoulliSurjectivity}
	Fix $0 < \eps < 1/6$ and $T > 0$, and set $q_n \coloneqq n^{-1/6+\eps}$. For every pair of deterministic matrices $A_n = (a_{ij}) \in \Z^{n \times (n+1)}$, $S_n = (s_{ij}) \in \{\pm1\}^{n\times(n+1)}$ satisfying $\abs*{a_{ij}} \leq n^T$, let $M_n = (m_{ij}) \in \mathbb Z^{n\times(n+1)}$ be given by $m_{ij} \coloneqq a_{ij} + s_{ij}\xi_{ij}$ with $\xi_{ij} \sim \Ber(q_n)$,
	where all the $\xi_{ij}$ are independent. Then, with uniform convergence in $A_n$ and $S_n$,
	\begin{equation*}
		\lim_{n \to \infty }\Prob\paren*{M_n \colon \Z^{n+1} \to \Z^n \text{ is surjective}} = \prod_{k=2}^\infty \zeta(k)^{-1}.
	\end{equation*}
\end{theorem}
\end{samepage}

The proof in \cref{app:AffineBernoulliSurjectivity} adapts Nguyen and Wood's arguments \cite{Nguyen2022}.

For matrices, $\Vert A \Vert = \sigma_{\max}(A)$ denotes the Euclidean operator norm. We use the following special case of the matrix Bernstein inequality \cite[Theorem~1.4]{Tropp2012}, obtained from its rectangular form \cite[Theorem~1.6]{Tropp2012} with $R=1$ and $d_1 = d = d_2$.

\begin{theorem}[{cf. \cite[Theorem~1.6]{Tropp2012}}]
	\label{thm:MatrixBernstein}
	Consider a finite sequence $\lbrace Z_k\rbrace$ of independent, random matrices with dimension $d\times d$. Assume that $\E Z_k = 0$ and $\norm*{Z_k}\leq 1$ almost surely. Define the matrix-variance parameter
	\begin{equation*}
		\sigma^2:= \max\left\lbrace\norm*{\sum_{k}\E(Z_kZ^\top_k)}, \norm*{\sum_{k}\E(Z^\top_kZ_k)} \right\rbrace.
	\end{equation*}
	Then, for all $t>0$, 
	\begin{equation*}
		\Prob\paren*{\norm*{\sum_{k}Z_k}\geq t} \leq 2d \exp\paren*{\frac{-t^2}{2(\sigma^2 + t/3)}}.
	\end{equation*}
\end{theorem}
\begin{proof}[Proof of \cref{thm:IntegerCubes} (Lower Bounds)]
	We first treat the integer cubes $\mathcal{Q}_d(R)$.

	Let $h$ be a Hadamard order and $H\in \lbrace -1,1\rbrace^{h\times h}$ be a Hadamard matrix.
	Fix $0 < \eps < 1/6$, set $q \coloneqq h^{-1/6+\eps}$, and let $\xi_{ij}\sim\Ber(q)$ be independent Bernoulli random variables with parameter $q$ for $1\leq i \leq h$ and $1\leq j\leq h+1$. Define the random matrix $M = [A|b]$ with
	\begin{equation*}
		A_{ij} = H_{ij}(R-\xi_{ij}) \text{ for }1\leq i,j \leq h, \quad b_i = \xi_{i,h+1} \text{ for }1\leq i\leq h.
	\end{equation*}
	The matrix $M$ satisfies the assumptions of \cref{thm:AffineBernoulliSurjectivity} with $A_h = [RH|0]$ and $S_h = [-H|\bm{1}]$. Applying this result gives 
	\begin{equation*}
		\Prob\paren*{M:\Z^{h+1}\to\Z^h \text{ is surjective}} = \prod_{k=2}^\infty \zeta(k)^{-1} + o(1),
	\end{equation*}
	where the error term tends to zero along Hadamard orders, uniformly in $H$ and $1\leq R\leq h^T$ for fixed $T>0$. We complement this by showing that $\left|\det A\right|$ is large with probability tending to $1$. 
	
	Since $\E A = (R-q)H$ and $H$ is a Hadamard matrix, $\left|\det (\E A)\right| = (R-q)^h h^{h/2}$ and each singular value of $\E A$ equals $(R-q)\sqrt{h}$. We control the perturbation using \cref{thm:MatrixBernstein}.
	Abbreviate $W := A - \E A$. Its entries are $W_{ij} = (q - \xi_{ij}) H_{ij}$ and are therefore independent and satisfy $\E W_{ij} = 0$ with $\abs*{W_{ij}} \leq 1$ and $\E W_{ij}^2 = q(1 - q)$. Write $W = \sum_{i,j} W_{ij}e_ie_j^\top$ and set $Z_{i,j} := W_{ij}e_ie_j^\top$. By construction, we have $\E Z_{i,j} = 0$ and $\norm*{Z_{i,j}} = |W_{ij}|\leq 1$. Using that $\E W_{ij}^2 = q(1 - q)$, the matrix-variance parameter $\sigma^2$ in \cref{thm:MatrixBernstein} equals $hq(1-q)$. So, for all $t>0$, \cref{thm:MatrixBernstein} applied to $\sigma^2 \leq hq$ implies 
	\begin{equation*}
		\Prob\paren*{\norm*{W} \geq t} \leq 2h\exp\paren*{-\frac{t^2}{2(hq + t/3)}}.
	\end{equation*}
	Taking $t = 4 \sqrt{hq \log h}$, one can estimate that $\norm*{W} \leq 4 \sqrt{hq \log h}$ with probability $1 - o(1)$.
	On this event, consider $E \coloneqq (\E A)^{-1}W$, which captures the deviation of $A$ from $\E A$ as $A = \E A(I+E)$, and 
	\begin{equation*}
		\delta \coloneqq \norm*{E} \leq \norm*{(\E A)^{-1}}\norm*{W} = O\paren*{\frac{\sqrt{hq\log h}}{(R - q)\sqrt{h}}} = O\paren*{\sqrt{q\log h}}= o(1)
	\end{equation*}
	as $q = h^{-1/6 + \eps}$ and $R - q\geq 1/2$ for $h$ sufficiently large.
	By the reverse triangle inequality, every singular value of $I + E$ is at least $1 - \delta$. Thus, for all sufficiently large $h$, $\abs*{\det(I + E)} \geq (1 - \delta)^h \geq \exp\paren*{-2h\delta}$. 
	Combining the preceding estimates yields
	\begin{equation*}
		\abs*{\det{A}} = \abs*{\det{\E A}}\abs*{\det(I + E)} \geq (R-q)^hh^{h/2} \exp(-2h\delta)\geq R^h h^{h/2} \exp(-2h(\delta + q)).
	\end{equation*}
	Since $q = h^{-1/6+\eps}$ with $\eps < 1/6$, the exponent in the last line is $o(h)$, uniformly for $R \geq 1$. We conclude that, with probability $1 - o(1)$,
	\begin{equation*}
		\abs*{\det{A}} \geq 2^{-o(h)}R^{h} h^{h/2}.
	\end{equation*}
	Denote the event that $\abs*{\det{A}}$ satisfies the inequality above by $\calE_{\det}$ and the event that $M$ is surjective from $\Z^{h + 1} \to \Z^h$ by $\calE_{\gcd}$. Since the event $\calE_{\det}$ has probability $1 - o(1)$, both events occur simultaneously with positive probability for all sufficiently large $h$. Let $M$ be such a realization. Consider the element $c\in\ker M$ with $c_{i} = (-1)^i\det M^{(i)}$ for $1\leq i\leq h +1$, where $M^{(i)}$ is the matrix $M$ without column $i$. Since the event $\calE_{\gcd}$ holds, we know that $\gcd\set*{|\det M^{(i)}| : 1\leq i\leq h + 1} = 1$. So $c$ is a circuit by \eqref{eq:CircuitCofactors} with $D = 1$. From the event $\calE_{\det}$, we obtain
	\begin{equation*}
		\gamma(\mathcal{Q}_h(R))\geq \Vert c\Vert_1 \geq \abs*{c_{h+1}} = \abs*{\det A} \geq 2^{-o(h)}R^{h} h^{h/2}.
	\end{equation*}
	For an arbitrary sufficiently large dimension $d$, let $h \leq d$ be the largest Hadamard order not exceeding $d$. The known bound on gaps between consecutive Hadamard orders \cite[Section~2]{Brent2013} gives $d - h = O(d^{1/6})$.
	To cover $1\leq R\leq d^T$, apply the preceding argument with the fixed exponent $T+1$: for all sufficiently large $d$, we have $h\geq d/2$ and $R\leq d^T\leq 2^T h^T\leq h^{T+1}$.
	Consider the injective linear embedding $\iota \colon \Z^h \to \Z^d$ defined by $\iota(x_1,\ldots,x_h) \coloneqq (x_1,\ldots,x_h,x_1,\ldots,x_1)$. 
	It maps $\calQ_h(R)$ into $\calQ_d(R)$ and $\calS_h$ into $\calS_d$. Let $\iota(M)$ be the matrix with $\iota$ applied column-wise to $M$. Observe that $c$ remains a circuit of $\iota(M)$, since $\ker\iota(M)=\ker M$. Therefore, we get
	\begin{equation*}
		\gamma(\mathcal{Q}_d(R)) \geq  2^{-o(h)}R^{h}h^{h/2} = 2^{-o(d)}R^{d}d^{d/2}.
	\end{equation*}	
	Indeed, the logarithm of $R^{d}d^{d/2}/(R^{h}h^{h/2})$ is $O((d - h)\log(dR)) = O(d^{1/6}\log d) = o(d)$ uniformly for $1 \leq R \leq d^T$.
	The lower bound for $\mathcal{Q}_d(R)$ now follows from \cref{cor:CircuitDavenport}.
	
	For the sign cube, the preceding construction with $R=1$ can produce zero entries. Instead, let $h$ be a Hadamard order and $H$ a corresponding Hadamard matrix. After multiplying suitable columns by $-1$, we may assume that the first row of $H$ consists of ones. We perturb only the $(h-1)\times h$ submatrix obtained by deleting that row.
	Set $q=(h-1)^{-1/6+\eps}$ and let $\xi_{ij}\sim\Ber(q)$ be independent for $2\leq i\leq h$ and $1\leq j\leq h$. Define $M=[A|\bm{1}]$ by
	\begin{equation*}
		A_{1j}= H_{1j} = 1 , \qquad A_{ij}=H_{ij}(1-2\xi_{ij}) \text{ for }2\leq i\leq h, 1\leq j \leq h.
	\end{equation*}
	All columns of $M$ lie in $\mathcal{S}_h$. The row operations $R_i\leftarrow R_1-R_i$ for $2\leq i\leq h$ give
	\begin{equation*}
		M'=\begin{pmatrix}
			\bm{1}^\top & 1 \\
			B'& 0
		\end{pmatrix},
	\end{equation*}
	where $B'\in\lbrace 0,2\rbrace^{(h-1)\times h}$. Let $B = 1/2 B'$. Its entries are of the form $B_{i-1,j}= (1-H_{ij})/2 + H_{i,j}\xi_{ij}$ for $2\leq i\leq h, 1\leq j\leq h$. Applying \cref{thm:AffineBernoulliSurjectivity} to $B$ shows that $B$ defines a surjection $\Z^h\to\Z^{h-1}$ with constant probability. On this event, $B'$ maps $\Z^h$ onto $2\Z^{h-1}$. The last column of $M'$ is $e_1$, so the first coordinate can be set independently and $M'\Z^{h+1}=\Z\oplus2\Z^{h-1}$. The unimodular row operations preserve the lattice index, so $M\Z^{h+1}$ has index $2^{h-1}$ in $\Z^h$. Equivalently, the gcd of the $h\times h$ minors of $M$ is $2^{h-1}$.
	For the determinant estimate, observe that $\E A=\operatorname{diag}(1,1-2q,\ldots,1-2q)H$, and thus $\abs*{\det(\E A)}=(1-2q)^{h-1}h^{h/2}$.
	The smallest singular value of $\E A$ is $(1-2q)\sqrt h$. Applying \cref{thm:MatrixBernstein} and the same argument to $(A-\E A)/2$ yields $\norm*{A-\E A}=O(\sqrt{hq\log h})$ with probability $1-o(1)$. One can then follow the preceding arguments to obtain $\abs*{\det A}\geq 2^{-o(h)}h^{h/2}$. Choosing a realization where both the gcd and determinant events hold, we get
	\begin{equation*}
		\gamma(\mathcal{S}_h)\geq\frac{\abs*{\det A}}{2^{h-1}}\geq 2^{-o(h)}\frac{h^{h/2}}{2^{h-1}}.
	\end{equation*}
	The same embedding into dimension $d$ and \cref{cor:CircuitDavenport} yield the claimed lower bound for $\mathcal{S}_d$.
	
\end{proof}

\subsection{Upper Bounds}
\label{sec:CubeUpperBounds}

We use the following lemma of Paat, Weismantel, and Weltge \cite{Paat2020}. Its proof uses Olson's theorem on the Davenport constant of $\Z^d/p\Z^d$ for $p$ prime \cite{Olson1969}; the lemma was introduced to obtain proximity bounds for mixed-integer programs.

\begin{lemma}[{\cite[Lemma~1]{Paat2020}}]
	\label{lem:OlsonType}
	Let $d,k\in\Z_{\geq 1}$, $z_{1},\ldots,z_{k} \in\Z^d$, and $\alpha_1,\ldots,\alpha_k\geq 0$. If $\sum_{i=1}^k \alpha_i \geq d$, then there exist $\beta_i\in\lbrack 0,\alpha_i\rbrack$ for $i=1,\ldots,k$, not all zero, such that $\sum_{i=1}^k\beta_iz_i\in\Z^d$.
\end{lemma}

\begin{proof}[Proof of \cref{thm:IntegerCubes} (Upper Bounds)]
	Let $\mathcal{Q}_d \in \lbrace\mathcal{Q}_d(R),\mathcal{S}_d\rbrace$ for some $R\geq 1$ and $\mathcal{Q}_{\lfloor d/2\rfloor}$, $\mathcal{Q}_{\lceil d/2\rceil}$ be the corresponding sets in dimensions $\lfloor d/2\rfloor$ and $\lceil d/2\rceil$ respectively. Abbreviate $K_d := |\mathcal{Q}_{\lfloor d/2\rfloor}| + |\mathcal{Q}_{\lceil d/2\rceil}|$. 
	Suppose, for a contradiction, that an irreducible zero-sum sequence over $\mathcal{Q}_d$ has length larger than $K_d \gamma(\mathcal{Q}_d)$. Let $a_1,\ldots,a_s\in \mathcal{Q}_d$ be its pairwise distinct vectors with multiplicities $g_1,\ldots,g_s\in \N_{\geq 1}$, so that $\sum_{i=1}^{s} g_ia_i=0$. Set $M:= [a_1,\ldots,a_s]$ and $g := (g_1,\ldots,g_s)^\top$; thus $\mathbf{1}^\top g > K_d \gamma(\mathcal{Q}_d)$. We will use \cref{lem:OlsonType} to decompose $g$ into two nonzero nonnegative integer vectors in $\ker M$, contradicting irreducibility.
	
	We begin by constructing a totally unimodular matrix $W \in \lbrace 0,1\rbrace^{k\times s}$ such that $\mathbf{1}^\top$ and all rows of $M$ are linear combinations of the rows of $W$.
	Fixing $Wx$ will therefore preserve both $Mx$ and the total length $\mathbf{1}^\top x$.
	Partition the $d$ coordinates of $M$ into blocks of sizes $\floor*{d/2}$ and $\ceil*{d/2}$, for instance, take index sets $I_1 = \lbrace 1,\ldots, \floor*{d/2}\rbrace$ and $I_2 = \lbrace \floor*{d/2} + 1,\ldots,d\rbrace$. For $i=1,2$, define
	\begin{equation*}
		V_i \coloneqq \set*{(i,a_j\vert_{I_i}) : 1\leq j\leq s}.
	\end{equation*}
	The first component makes $V_1$ and $V_2$ disjoint, even when the two coordinate blocks contain identical patterns. Associate each column $a_j$ with an edge joining $(1,a_j\vert_{I_1})$ and $(2,a_j\vert_{I_2})$, and let $G = (V_1\sqcup V_2,E)$ be the resulting bipartite graph. Distinct columns of $M$ give distinct edges. Notice that $k := |V_1| + |V_2|\leq |\mathcal{Q}_{\floor*{d/2}}| + |\mathcal{Q}_{\ceil*{d/2}}| = K_d$. Let $W \in \{0,1\}^{k \times s}$ be its vertex-edge incidence matrix, with rows ordered first by $V_1$ and then by $V_2$, and column $j$ corresponding to $a_j$. Since $G$ is bipartite, the matrix $W$ is totally unimodular; cf. \cite{Schrijver1986}.
	For $i=1,2$, let $P_i$ have as its columns the underlying coordinate patterns of the vertices in $V_i$, omitting their first component and using the same order as the corresponding rows of $W$. Then
	\begin{equation}
		\label{eq:UpperBoundTransformation}
		M =
		\begin{pmatrix}
			P_1 & 0 \\
			0 & P_2
		\end{pmatrix}W,
		\qquad
		\mathbf{1}^\top = \left(\mathbf{1}_{|V_1|}^\top, 0\right)W.
	\end{equation}
	Now, take $g$ and normalize it by $\theta := K_d \gamma(\mathcal{Q}_d)/\mathbf{1}^\top g$. By assumption, we have $\theta < 1$. Express $\theta g$ as a nonnegative combination of nonnegative circuits $c_j$ of $\ker M$, i.e.,
	\begin{equation*}
		\theta g = \sum_j \alpha_jc_j, \qquad \alpha_j \geq 0.
	\end{equation*}
	The triangle inequality implies $\theta \mathbf{1}^\top g \leq \gamma(\mathcal{Q}_d) \sum_{j}\alpha_j$. Rearranging and using the definition of $\theta$ gives $K_d\leq \sum_{j} \alpha_j$. Applying \cref{lem:OlsonType} to the integer vectors $Wc_j\in\Z^{k}$, we obtain $0\leq \beta_j\leq \alpha_j$, not all zero, such that $\sum_{j}\beta_jWc_j\in\Z^{k}$. Set $v:= \sum_{j}\beta_jc_j$. Then, nonnegativity and $\theta < 1$ yield
	\begin{equation*}
		0 \neq v \leq \theta g, \qquad Mv = 0, \qquad Wv \in \Z^k.
	\end{equation*}
	To round $v$ to an integer vector, we use the total unimodularity of $W$. Consider the polytope $P := \lbrace x\in \R^s: Wx = Wv, 0\leq x\leq g\rbrace$. The properties of $v$ above ensure that $P$ has an integer-valued right-hand side and is nonempty as $v\in P$. Since $W$ is totally unimodular, there exists $y\in P\cap \Z^s$. Moreover, the equation $Wy = Wv$ implies $My = Mv = 0$ and $\mathbf{1}^\top y = \mathbf{1}^\top v$ by \eqref{eq:UpperBoundTransformation}. Since $v \leq \theta g$ and $v \neq 0$, we have $0 < \mathbf{1}^\top y < \mathbf{1}^\top g$, which implies $y\neq g$ and $y\neq 0$. So $g = y + (g - y)$ is a proper decomposition of $g$, a contradiction.
	
	Thus each irreducible zero-sum sequence has length at most $K_d\gamma(\mathcal{Q}_d)$. The claimed bounds in \cref{thm:IntegerCubes} for $\mathcal{Q}_d(R)$ and $\mathcal{S}_d$ follow from the definition of $K_d$ and $\gamma(\mathcal{Q}_d)$ using \cref{lem:GeneralCircuitUpperBound}.
\end{proof}
The encoding in \eqref{eq:UpperBoundTransformation} also appears in work on the integrality number of an integer program \cite[Section~3]{Paat2022}. For any integer program in equality form $Ax=b,\; l\leq x\leq u$ with integer-valued $b,u,l$ and columns of $A$ in $\mathcal{Q}_d$ from the proof, it gives an integrality number at most $K_d$; cf. \cite[Lemma~5]{Paat2022}. The constraint matrix also has affine TU-dimension at most $K_d$ \cite{Bader2018}; a related knapsack construction appears in \cite[Example~8]{Bader2018}.

\begin{remark}[Conditional upper bound from the Euclidean Steinitz estimate]
    \label{rem:SteinitzCubeUpper}
    Assume the Euclidean Steinitz estimate stated in \cite[Theorem~1.2]{OpenAI2026a}: there is an absolute constant $C_{\mathrm{St}} > 0$ such that every finite zero-sum sequence in the Euclidean unit ball of $\R^d$ has an ordering whose partial sums all have Euclidean norm at most $C_{\mathrm{St}}\sqrt{d}$. Then there is an absolute constant $K > 0$ such that, for every $d \geq 1$ and $R \geq 1$,
    \begin{equation*}
        \Dvp(\calQ_d(R)) \leq (KR\sqrt{d})^d.
    \end{equation*}

    To see this, let $a_1,\ldots,a_L \in \calQ_d(R)$ form an irreducible zero-sum sequence. Since $\norm*{a_i}_2 \leq R\sqrt{d}$, we apply the assumed estimate to $a_i/(R\sqrt{d})$. Let $s_j$ be the partial sums of the original vectors in the resulting order. Rescaling gives $\norm*{s_j}_2 \leq C_{\mathrm{St}}Rd$. The lattice points $s_0,\ldots,s_{L-1}$ are pairwise distinct: an equality $s_i = s_j$ with $0 \leq i < j \leq L - 1$ would give a nonempty proper zero-sum subsequence.

    The unit cubes centered at these $L$ lattice points have disjoint interiors and lie in the Euclidean ball of radius $C_{\mathrm{St}}Rd + \sqrt{d}/2$. Writing $v_d$ for the volume of the Euclidean unit ball and using $v_d \leq (2\pi e/d)^{d/2}$, we obtain
    \begin{equation*}
        \begin{aligned}
            L \leq v_d\paren*{C_{\mathrm{St}}Rd + \frac{\sqrt{d}}{2}}^d 
            \leq \paren*{\sqrt{2\pi e}\paren*{C_{\mathrm{St}}R\sqrt{d} + \frac{1}{2}}}^d 
            \leq (KR\sqrt{d})^d,
        \end{aligned}
    \end{equation*}
    where $K \coloneqq \sqrt{2\pi e}(C_{\mathrm{St}} + 1/2)$ suffices. The same bound with $R = 1$ applies to $\calS_d \subseteq \calQ_d(1)$. 
    We use the external Steinitz estimate only in this remark.
\end{remark}

\section{The Average-Case Davenport Constant}
\label{sec:AverageCase}

We prove \cref{thm:AverageCase} using the first and second moment methods. Let $n \geq 2$ be even, and let $X_1,\ldots,X_n$ be independent and uniformly distributed on $\calS_d$. We write $\Prob_0$ and $\E_0$ for probability and expectation conditional on $\sum_{i=1}^n X_i = 0$, and define
\begin{equation*}
    Z_n \coloneqq \#\set*{\varnothing \neq I \subsetneq [n] : \sum_{i \in I} X_i = 0}.
\end{equation*}
Thus, $p_d(n) = \Prob_0(Z_n > 0)$. To anticipate the threshold, for each even $n \geq 4$ choose an even integer $s=s(n)$ with $\abs*{s-n/2} \leq 1$. There are $\binom{n}{s} = \exp(n\log{2} + O(\log{n}))$ such subsets, and each has conditional zero-sum probability $\exp(-\frac{d}{2}\log{n} + O(d))$. The leading terms balance when $n\log{2} \approx \frac{d}{2}\log{n}$, suggesting a transition at $n \sim d\log{d}/\log{4}$. The following moment estimates establish both sides of this transition.

We first show that below the threshold the expected number of proper zero-sum subsequences tends to zero.

\begin{lemma}[First-moment estimate]
    \label{lem:FirstMoment}
    For every fixed $\eta \in (0,1)$,
    \begin{equation*}
        \lim_{d \to \infty} \E_0 Z_n = 0,
    \end{equation*}
    uniformly over all even integers $n$ satisfying
    \begin{equation*}
        2 \leq n \leq (1 - \eta)\frac{d\log{d}}{\log{4}}.
    \end{equation*}
\end{lemma}

The proof is given in \cref{app:FirstMoment}. By Markov's inequality,
\begin{equation*}
    p_d(n) = \Prob_0(Z_n > 0) \leq \E_0 Z_n.
\end{equation*}
Thus, \cref{lem:FirstMoment} establishes the first probability limit in \cref{thm:AverageCase}.

Above the threshold, it suffices to count zero-sum subsequences whose lengths lie in a central interval. For $\gamma \in (0,1/2)$, define
\begin{equation*}
    \widehat{Z}_n(\gamma) \coloneqq \#\set*{I \subseteq [n] : \gamma n \leq \abs*{I} \leq (1 - \gamma)n, \sum_{i \in I} X_i = 0}.
\end{equation*}
For two subsets of fixed sizes $s$ and $t$, the joint conditional zero-sum probability depends only on their intersection size. We therefore average over the hypergeometric overlap distribution. A local approximation and concentration control the regular range, where all four regions determined by the subsets have sizes proportional to $n$; an entropy bound makes the boundary contribution negligible after normalization.

\begin{lemma}[Second-moment estimate]
    \label{lem:SecondMoment}
    For every fixed $\eta \in (0,1)$, there exists a constant $\gamma = \gamma(\eta) \in (0,1/2)$ such that
    \begin{equation*}
        \lim_{d \to \infty} \frac{\E_0\widehat{Z}_n(\gamma)^2}{(\E_0\widehat{Z}_n(\gamma))^2} = 1,
    \end{equation*}
    uniformly over all even integers $n$ satisfying
    \begin{equation*}
        n \geq (1 + \eta)\frac{d\log{d}}{\log{4}}.
    \end{equation*}
\end{lemma}

The proof is given in \cref{app:SecondMoment}. Set $W \coloneqq \widehat{Z}_n(\gamma)$. Cauchy--Schwarz under $\Prob_0$ gives
\begin{equation*}
    (\E_0 W)^2 = \bigl(\E_0\paren*{W\mathbf{1}_{\{W>0\}}}\bigr)^2 \leq \E_0 W^2\,\Prob_0\paren*{W>0}.
\end{equation*}
Since $\widehat{Z}_n(\gamma) \leq Z_n$, it follows that
\begin{equation*}
    p_d(n) \geq \Prob_0\paren*{\widehat{Z}_n(\gamma) > 0} \geq \frac{(\E_0\widehat{Z}_n(\gamma))^2}{\E_0\widehat{Z}_n(\gamma)^2}.
\end{equation*}
By \cref{lem:SecondMoment}, the last quotient tends to one. This establishes the second probability limit in \cref{thm:AverageCase}.

We finally derive the asserted limit for $D_\eps(\nu_d)$. Fix $\eps,\eta \in (0,1)$. For all sufficiently large $d$, the uniform first-moment estimate gives $p_d(n) < 1 - \eps$ for every even integer $n \leq (1 - \eta)d\log{d}/\log{4}$. The second-moment estimate gives $p_d(n) \geq 1 - \eps$ at the smallest even integer satisfying $n \geq (1 + \eta)d\log{d}/\log{4}$. By the definition of $D_\eps(\nu_d)$, it follows that
\begin{equation*}
    (1 - \eta)\frac{d\log{d}}{\log{4}} < D_\eps(\nu_d) \leq (1 + \eta)\frac{d\log{d}}{\log{4}} + 2.
\end{equation*}
Dividing by $d\log{d}$, letting $d \to \infty$, and then letting $\eta \downarrow 0$ proves
\begin{equation*}
    \lim_{d \to \infty} \frac{D_\eps(\nu_d)}{d\log{d}} = \frac{1}{\log{4}}.
\end{equation*}

\begin{remark}[Rademacher sums]
    \label{rem:RademacherSums}
    Fix an integer $m \geq 1$, let $\calR_m$ be the distribution of a sum of $m$ independent uniform signs, and put $\nu_{m,d} \coloneqq \calR_m^{\otimes d}$. For $n \geq 2$ with $mn \in 2\Z$, define $p_{m,d}(n)$ by the same conditional reducibility probability as $p_d(n)$, with independent $X_i \sim \nu_{m,d}$, and put $D_\eps(\nu_{m,d}) \coloneqq \min\{n \geq 2 : mn \in 2\Z,\ p_{m,d}(n) \geq 1-\eps\}$. With these quantities, the conclusions of \cref{thm:AverageCase} remain valid over these admissible lengths. The moment argument in \cref{app:RademacherSums} uses $mn$ underlying signs per coordinate, while subset and overlap counts still concern the $n$ vector blocks; subset sizes $1 \leq s < n$ contribute only when $ms$ is even. For $m > 1$, this distribution is nonuniform on $\{-m,-m+2,\ldots,m\}^d$.
\end{remark}

\section{Lattice Balls and Sparse Ternary Cubes}
\label{sec:LatticeBalls}

We prove \cref{thm:LatticeBalls} using circuit bounds and explicit irreducible zero-sum sequences. For $p \geq 2$, a Euclidean circuit estimate matches a sparse ternary construction; $p = 1$ requires a separate argument. Throughout, $p$ and $R$ are fixed as $d \to \infty$.

For fixed $p$ and $R$, put $k \coloneqq \floor*{R^p}$. Every vector in $\calB_{p,d}(R)$ has at most $k$ nonzero coordinates, each of magnitude at most $\lfloor R\rfloor$. Hence
\begin{equation*}
    \abs*{\calB_{p,d}(R)} \leq \sum_{j=0}^{\min\{k,d\}} \binom{d}{j}(2\floor*{R})^j = O_{p,R}(d^k).
\end{equation*}
An irreducible zero-sum sequence containing the zero vector has length one. For every other such sequence of length $L$, let $s$ be the number of distinct vectors and $r$ the rank of their column matrix. Since $s-r \leq s \leq \abs*{\calB_{p,d}(R)}$, a uniform circuit bound $F$ gives $L \leq (s-r)F \leq \abs*{\calB_{p,d}(R)}F$ by \cref{lem:CircuitTransfer}. The bound on $\abs*{\calB_{p,d}(R)}$ shows that this transfer incurs only a polynomial factor, whose $d$-th root tends to one. It therefore does not affect the exponential rate.

\paragraph{Euclidean upper bound.}

\begin{lemma}[Euclidean circuit bound]
    \label{lem:EuclideanCircuitBound}
    Let $k \geq 3$ be an integer, and let $a_1,\ldots,a_{h+1} \in \Z^d$ form a positive circuit of rank $h \geq 1$. If $\norm*{a_j}_2^2 \leq k$ for every $1 \leq j \leq h+1$, then its primitive positive integer relation $g \in \Z_{>0}^{h+1}$ satisfies
    \begin{equation*}
        \norm*{g}_1 \leq 3(h+1)\rho_k^h.
    \end{equation*}
\end{lemma}

The proof in \cref{app:EuclideanCircuitBound} combines the cofactor formula with an integer Gram determinant estimate that exploits the circuit structure.

Fix an integer $k \geq 3$. Every positive circuit supported on nonzero vectors in $\calB_{2,d}(\sqrt{k})$ has rank $h \geq 1$ because it contains no zero column, while $h \leq d$ follows from the ambient dimension. By \cref{lem:EuclideanCircuitBound} and $\rho_k > 1$, its primitive positive relation $g$ satisfies $\norm*{g}_1 \leq 3(h+1)\rho_k^h \leq 3(d+1)\rho_k^d$. By \cref{lem:CircuitTransfer} and the preceding bound on $\abs*{\calB_{2,d}(\sqrt{k})}$, there is a constant $C_k > 0$, depending only on $k$, such that $\Dvp(\calB_{2,d}(\sqrt{k})) \leq C_k d^{k+1}\rho_k^d$ for every $d \geq 1$. Taking $d$-th roots yields
\begin{equation}
    \label{eq:EuclideanUpperBound}
    \limsup_{d \to \infty} \Dvp(\calB_{2,d}(\sqrt{k}))^{1/d} \leq \rho_k.
\end{equation}

\paragraph{Ternary lower bound and the case $p \geq 2$.}

\begin{lemma}[Ternary lower bound]
    \label{lem:TernaryLowerBound}
    For every fixed integer $k \geq 3$,
    \begin{equation*}
        \liminf_{d \to \infty} \Dvp(\calT_d(k))^{1/d} \geq \rho_k.
    \end{equation*}
\end{lemma}

The block construction in \cref{app:TernaryConstruction} produces a large maximal minor while keeping the greatest common divisor of all maximal minors equal to one. After discarding zero-cofactor columns and changing column signs as needed, the resulting positive circuit relation is primitive and hence defines an irreducible zero-sum sequence. The block determinant is obtained by raising the terms of a second-order recurrence with dominant root $\rho_k^2$ to the fixed block size. Each recurrence step adds two blocks, giving rate $\rho_k$ per coordinate; padding with zero coordinates covers all sufficiently large dimensions.

Fix $2 \leq p < \infty$ and $R \geq 0$ with $k \coloneqq \floor*{R^p} \geq 3$. We have $\calT_d(k) \subseteq \calB_{p,d}(R) \subseteq \calB_{2,d}(\sqrt{k})$. The first inclusion follows from $\norm*{v}_p^p = \abs*{\supp(v)}$ for ternary vectors. For the second, use $\abs*{a}^2 \leq \abs*{a}^p$ for $a \in \Z$ and the integrality of $\norm*{v}_2^2$. By monotonicity, \cref{lem:TernaryLowerBound} gives
\begin{equation*}
    \liminf_{d \to \infty} \Dvp(\calB_{p,d}(R))^{1/d} \geq \rho_k,
\end{equation*}
while \eqref{eq:EuclideanUpperBound} gives the reverse inequality for the upper limit. This proves the first assertion of \cref{thm:LatticeBalls}.

\paragraph{The $p = 1$ case.}

For fixed $R \geq 3$, put $M \coloneqq \floor*{R} \geq 3$. Since integer vectors have integral $\ell_1$-norm, $\calB_{1,d}(R) = \calB_{1,d}(M)$. It therefore suffices to treat fixed integer radii $M \geq 3$.

\begin{lemma}[Bounds for $p = 1$]
    \label{lem:SpecialCase}
    For every fixed integer $M \geq 3$, there is a constant $C_M > 0$ such that, for every $d \geq 1$,
    \begin{equation*}
        \frac{(M-1)^{d+1}-1}{M-2} \leq \Dvp(\calB_{1,d}(M)) \leq C_M d^{M+1}(M-1)^d.
    \end{equation*}
\end{lemma}

The upper bound is proved in \cref{app:SpecialCase} by combining a circuit estimate with \cref{lem:CircuitTransfer}.

For the lower bound, we use the bidiagonal construction of \cite[Lemma~3]{Hunkenschroeder2026}. Put $P = M - 1$. The identity
\begin{equation*}
    P^d e_1 + \sum_{i=1}^{d-1} P^{d-i}(-Pe_i + e_{i+1}) - Pe_d = 0
\end{equation*}
specifies a zero-sum sequence whose vectors all have $\ell_1$-norm at most $M$.

The coordinate equations force the multiplicities of any zero-sum subsequence to be proportional to these counts. Since $-Pe_d$ occurs only once, the subsequence is either empty or the entire sequence. Its length is the geometric sum $(P^{d+1}-1)/(P-1)$.

Taking $d$-th roots in \cref{lem:SpecialCase} and letting $d \to \infty$ proves the second assertion of \cref{thm:LatticeBalls}.

\begin{proof}[Proof of \cref{cor:SparseTernaryProximity}]
    Put $\calX=\calT_d(k)$, $s=\abs*{\calX}=O_k(d^k)$, and $C=\gamma(\calX)$. By \cref{cor:CircuitDavenport}, $C\leq\Dvp(\calX)\leq sC$. We bound $\Pi(\calX)$ between $C-d-1$ and $sC$.

    For the upper bound, fix an optimal LP vertex $x^*$. For each column type retain a cheapest variable, choosing the one positive at $x^*$ if present. A vertex uses at most one copy, and optimality makes that copy cheapest. Moving mass to the representatives preserves feasibility and integrality and can only decrease cost. Both optimal values are therefore preserved. Write $x^*$ also for its restriction, since lifting by zeros preserves the distance. For a nearest integer optimum $z$ of the reduced program, decompose $z-x^*$ conformally into at most $s$ primitive circuits with nonnegative coefficients, by conic Carath\'eodory in its orthant (see \cite[Section~2]{Ekbatani2022}). A coefficient at least one would allow the feasible exchanges $x^*+g$ and $z-g$ for its circuit $g$. Optimality forces $c^\top g=0$, so $z-g$ would be a closer integer optimum. Thus each coefficient is less than one and the distance is at most $sC$.

    For the lower bound, restrict a circuit attaining $C$ to its support and orient its columns. The resulting matrix $B$ has columns in $\calX$ and $\ker B=\R g$ for a primitive $g>0$. Choose $j$ with $G=g_j=\norm*{g}_\infty$. The feasible points of $\min\{\bm{1}^\top x: Bx=-b_j,\ x\geq0\}$, where $b_j$ is column $j$ of $B$, are $\lambda g-e_j$ with $\lambda\geq1/G$. By primitivity, such a point is integral precisely when $\lambda$ is integral. Thus the unique LP and integer optima have $\lambda=1/G$ and $\lambda=1$, respectively. Their distance is $(1-1/G)C\geq C-d-1$, since $g$ has at most $d+1$ coordinates.

    The two bounds, $s=O_k(d^k)$, and \cref{cor:SparseTernaryCubes} with $\rho_k>1$ give the claimed limit.
\end{proof}

\appendix
\section{Affine Bernoulli Surjectivity}
\label[appendix]{app:AffineBernoulliSurjectivity}

We prove \cref{thm:AffineBernoulliSurjectivity}, which we restate below. For a prime $p$, write $\mathbb{F}_p \coloneqq \Z/p\Z$. An integer-valued random variable $Y$ is $\alpha$-balanced if, for every prime $p$,
\begin{equation*}
    \max_{r \in \mathbb{F}_p} \Prob\paren*{Y \equiv r \bmod p} \leq 1 - \alpha.
\end{equation*}
For a fixed prime $p$, a probability distribution $\nu$ on $\mathbb{F}_p$ is called $\alpha$-balanced if
\begin{equation*}
    \max_{r \in \mathbb{F}_p} \nu(\{r\}) \leq 1 - \alpha.
\end{equation*}
We use the same term for a random variable with this distribution, and coordinatewise for random vectors. In particular, reducing an $\alpha$-balanced integer-valued random variable modulo $p$ gives an $\alpha$-balanced $\mathbb{F}_p$-valued random variable.

\begin{restatedtheorem}[Affine Bernoulli Surjectivity]{thm:AffineBernoulliSurjectivity}
    Fix $0 < \eps < 1/6$ and $T > 0$, and set $q_n \coloneqq n^{-1/6+\eps}$. For every pair of deterministic matrices $A_n = (a_{ij}) \in \Z^{n \times (n+1)}$, $S_n = (s_{ij}) \in \{\pm1\}^{n\times(n+1)}$ satisfying $\abs*{a_{ij}} \leq n^T$, let $M_n = (m_{ij}) \in \mathbb Z^{n\times(n+1)}$ be given by $m_{ij} \coloneqq a_{ij} + s_{ij}\xi_{ij}$ with $\xi_{ij} \sim \Ber(q_n)$,
    where all the $\xi_{ij}$ are independent. Then, with uniform convergence in $A_n$ and $S_n$,
    \begin{equation*}
        \lim_{n \to \infty }\Prob\paren*{M_n \colon \Z^{n+1} \to \Z^n \text{ is surjective}} = \prod_{k=2}^\infty \zeta(k)^{-1}.
    \end{equation*}
\end{restatedtheorem}

For comparison with \cite[Theorem~1.1]{Nguyen2022}, take $u=1$ and $B = \lbrace 0\rbrace$. Set $\alpha_n \coloneqq q_n$. For all sufficiently large $n$, we have $\alpha_n\leq1/2$. Every entry of $M_n$ takes the two adjacent values $a_{ij}$ and $a_{ij} + s_{ij}$ with probabilities $1-\alpha_n$ and $\alpha_n$, respectively. The two values are distinct modulo every prime, so the largest residue probability is $1-\alpha_n$. Thus the entries are independent $\alpha_n$-balanced integer random variables.

We deduce \cref{thm:AffineBernoulliSurjectivity} from the following estimates for three ranges of primes. All constants below may depend only on $\eps$ and $T$, but not on $p$ or the deterministic matrices $A_n$ and $S_n$. We suppress the dependence on $\eps$ and $T$ in our notation.

\begin{proposition}[{cf. \cite[Propositions~2.1--2.3]{Nguyen2022}}]
	\label{prop:ThreePrimeRegimes}
	There exists a constant $c_0 > 0$ such that, with $U_n \coloneqq \exp(c_0\alpha_n n/2)$, for every fixed $L \geq 2$ the following statements hold uniformly in $A_n$ and $S_n$:
	\begin{enumerate}
		\item
		\label{item:SmallPrimes}
		\begin{equation*}
			\lim_{n \to \infty} \Prob\paren*{\rank_{\mathbb{F}_p}M_n = n \text{ for every prime } p \leq L} = \prod_{\substack{p \leq L \\ p \text{ prime}}} \prod_{k=2}^\infty (1-p^{-k}),
		\end{equation*}
		\item 
		\label{item:MediumPrimes}
		there exists a constant $C>0$, independent of $L$, such that, for all sufficiently large $n$, 
		\begin{equation*}
			\Prob\paren*{\rank_{\mathbb{F}_p} M_n < n \text{ for some prime } L < p \leq U_n} \leq \frac{C}{L} + C\exp\paren*{-c_0\alpha_n n/2},
		\end{equation*}
		\item
		\label{item:LargePrimes}
		\begin{equation*}
			\Prob\paren*{\rank_{\mathbb{F}_p} M_n < n \text{ for some prime } p \geq U_n} = o(1).
		\end{equation*}
	\end{enumerate}
\end{proposition}

\begin{proof}[Proof of \cref{thm:AffineBernoulliSurjectivity}]
	By Smith normal form (see, e.g., \cite[Section~4.4]{Schrijver1986}), the matrix $M_n$ defines a surjection $\Z^{n+1} \to \Z^n$ if and only if the greatest common divisor of its $n\times n$ minors is $1$. The latter is equivalent to the fact that $M_n \mod p$ has full row rank for all primes $p$. Hence, it suffices to consider the rank over all prime fields $\mathbb{F}_p$. Fix $L \geq 2$. We apply \cref{prop:ThreePrimeRegimes} when $n\to \infty$. Part~\labelcref{item:MediumPrimes} gives
	\begin{equation*}
		\limsup_{n\to \infty}\Prob\paren*{\rank_{\mathbb{F}_p} M_n < n \text{ for some prime } L < p \leq U_n} \leq \frac{C}{L},
	\end{equation*}
	as $\alpha_n n= q_n n= n^{5/6+\eps}\to \infty$. Part~\labelcref{item:LargePrimes} gives
	\begin{equation*}
		\lim_{n\to \infty}\Prob\paren*{\rank_{\mathbb{F}_p} M_n < n \text{ for some prime } p \geq U_n} = 0.
	\end{equation*}
	Let $\calE_n$ be the event that $\rank_{\mathbb{F}_p}M_n = n$ for all primes $p$. The two bounds above, part~\labelcref{item:SmallPrimes} of \cref{prop:ThreePrimeRegimes}, and the union bound imply the lower bound in
	\begin{equation*}
		\prod_{\substack{p \leq L \\ p \text{ prime}}} \prod_{k=2}^\infty (1-p^{-k}) - \frac{C}{L}\leq\liminf_{n \to \infty} \Prob\paren*{\calE_n}\leq\limsup_{n \to \infty} \Prob\paren*{\calE_n}\leq \prod_{\substack{p \leq L \\ p \text{ prime}}} \prod_{k=2}^\infty (1-p^{-k}).
	\end{equation*}
	The upper bound follows from part~\labelcref{item:SmallPrimes} because $\calE_n$ implies full row rank modulo every prime $p\leq L$.
	Letting $L$ go to infinity and interchanging the products proves $\lim_{n \to \infty}\Prob\paren*{\calE_n} = \prod_{k=2}^\infty \zeta(k)^{-1}$. The uniform convergence follows since all three statements of the proposition are uniform.
\end{proof}

The next three subsections prove \cref{prop:ThreePrimeRegimes} for small, medium, and large primes, respectively.

\subsection{Small Primes}

The small primes are treated jointly by applying \cite[Theorem~4.1]{Nguyen2022}, which generalizes \cite[Corollary~3.4]{Wood2019} using a refinement of its methods. We first introduce some group-theoretic notation; see \cite{Dummit2004} for background. For a finite abelian group $G$ and a finite set of primes $P$, we write $G_P := \bigoplus_{p\in P}G_p$, where $G_p$ is the Sylow $p$-subgroup of $G$. We consider the group $\Z^n / M_n \Z^{n+1}$, which is finitely generated but not necessarily finite. The following result gives the limiting probability that it is finite and its Sylow $p$-subgroups are trivial for every $p\in P$. We state the consequence of the proof of the cited theorem for $u = 1$ and $B = \lbrace 0\rbrace$.\footnote{Although the finiteness condition stated in \cref{thm:NguyenWoodSmallPrimes} is not explicit in \cite[Theorem~4.1]{Nguyen2022}, the final step of the proof in \cite[p.~25]{Nguyen2022} establishes the limit for $(\Z^n/M_n\Z^{n+1})\otimes\Z/a\Z=\{0\}$ with $a=\prod_{p\in P}p$. Since $P\neq\emptyset$, this implies that $\Z^n/M_n\Z^{n+1}$ is finite.}

\Needspace{9\baselineskip}
\begin{theorem}[{\cite[Theorem~4.1]{Nguyen2022}}]
	\label{thm:NguyenWoodSmallPrimes}
	Let $\alpha_n$ be a function of integers $n$ such that for any $\Delta > 0$ we have $\alpha_n n\geq \Delta\log n$ for $n$ sufficiently large. For every positive integer $n$, let $M_n \in \Z^{n \times (n+1)}$ be a random matrix with independent $\alpha_n$-balanced entries. Let $P$ be a nonempty finite set of primes. Then
	\begin{equation*}
		\lim_{n\to \infty} \Prob\paren*{\Z^n / M_n\Z^{n+1}\text{ is finite and }(\Z^n / M_n\Z^{n+1})_P \cong \lbrace 0\rbrace} = \prod_{p\in P}\prod_{k=2}^\infty (1 - p^{-k}).
	\end{equation*}
\end{theorem}

We deduce part~\labelcref{item:SmallPrimes} of \cref{prop:ThreePrimeRegimes} from \cref{thm:NguyenWoodSmallPrimes}.

\begin{proof}[Proof of \cref{prop:ThreePrimeRegimes}, part~\labelcref{item:SmallPrimes}]
	Set $G_n \coloneqq \Z^n/M_n\Z^{n+1}$, and let $P$ be the set of primes at most $L$. This set is nonempty because $L\geq 2$. The entries of $M_n$ are independent and $\alpha_n$-balanced, and $\alpha_n n/\log{n} = n^{5/6 + \eps}/\log{n} \to \infty$. Applying \cref{thm:NguyenWoodSmallPrimes} to $P$ gives
	\begin{equation*}
		\lim_{n\to \infty} \Prob\paren*{G_n\text{ is finite and }(G_n)_P \cong \lbrace 0\rbrace} = \prod_{p\in P}\prod_{k=2}^\infty (1 - p^{-k}).
	\end{equation*}
	For finite $G_n$, the property that $(G_n)_P \cong \lbrace 0\rbrace$ is equivalent to $(G_n)_p\cong\lbrace 0\rbrace$ for all $p\in P$. The latter is equivalent to $p\nmid |G_n|$. By Smith normal form and $P\neq \emptyset$, the event above is therefore equivalent to $M_n\mod p$ having full row rank for every $p\in P$.

	Finally, the convergence is uniform in $A_n$ and $S_n$. Otherwise, there exist $\delta > 0$, dimensions $n_1 < n_2 < \cdots$, and pairs $(A_{n_j},S_{n_j})$ whose probabilities differ from $\prod_{p\in P}\prod_{k=2}^\infty (1 - p^{-k})$ by at least $\delta$ for every $j\geq 1$. Complete this sequence by taking $A_n=0$ and setting $s_{ij}=1$ for all $i,j$ at every missing dimension. The resulting random matrices still have independent $\alpha_n$-balanced entries. So their full-rank probabilities converge to $\prod_{p\in P}\prod_{k=2}^\infty (1 - p^{-k})$ by \cref{thm:NguyenWoodSmallPrimes}. This contradicts the deviation of at least $\delta$ along the subsequence indexed by $n_j$.
\end{proof}

\subsection{Medium Primes}

We establish the fixed-prime rank estimate in \cref{lem:FixedPrimeRank} and sum it over the medium primes. Some intermediate estimates will also be used for large primes.

Fix an arbitrary prime $p$. We expose the columns successively using the following notation. Let $X_1,\ldots,X_{n+1} \in \mathbb{F}_p^n$ denote the columns of $M_n$ reduced modulo $p$. Thus, $(X_j)_i = a_{ij} + s_{ij}\xi_{ij} \in \mathbb{F}_p$, where the deterministic coefficients are understood modulo $p$. For $1 \leq m \leq n+1$, set
\begin{equation*}
    W_m \coloneqq \operatorname{span}_{\mathbb{F}_p}(X_1,\ldots,X_m).
\end{equation*}
For a subspace $V \subseteq \mathbb{F}_p^n$, let $V^\perp$ denote its orthogonal complement with respect to the standard dot product over $\mathbb{F}_p$, and write $\operatorname{codim} V = n - \dim V$.

For a subspace $V$ of codimension $r$, the containment probability for a uniform vector is $p^{-r}$. We control deviations from this benchmark by excluding sparse normal vectors, then using Fourier capture for moderate deviations and symmetrization and swapping for large ones. This gives \cref{cor:RandomSpanEstimate}, which allows us to compare the rank process with that of a uniform matrix.

We begin by recording two standard anti-concentration estimates that will be used repeatedly below. The first is a non-identically distributed variant of Odlyzko's lemma \cite{Kahn1995,Odlyzko1988}. Although Nguyen and Wood state the result for independent and identically distributed entries, their proof does not require the entries to be identically distributed.

\begin{lemma}[{Odlyzko's lemma, \cite[Lemma~3.1]{Nguyen2022}}]
    \label{lem:NonIdenticalOdlyzko}
    Let $\mathbb{F}$ be a field, and let $Y = (Y_1,\ldots,Y_N) \in \mathbb{F}^N$ have independent coordinates satisfying $\max_{x \in \mathbb{F}} \Prob\paren*{Y_i = x} \leq 1 - \alpha$ for every $1 \leq i \leq N$. If $H \subseteq \mathbb{F}^N$ is a linear subspace of codimension $r$, then
    \begin{equation*}
        \Prob\paren*{Y \in H} \leq (1 - \alpha)^r.
    \end{equation*}
\end{lemma}

The next anti-concentration estimate is an application of \cite[Theorem~6.3]{Nguyen2022} tailored to our specific affine Bernoulli setting. Nguyen and Wood refer to Maples \cite{Maples2010} for a proof of their result and mention that it is based on arguments by Halász.

\begin{lemma}[{Finite-field Erd\H{o}s--Littlewood--Offord estimate, cf. \cite[Theorem~6.3]{Nguyen2022}}]
    \label{lem:ErdosLittlewoodOfford}
    Let $I \subseteq [n]$, let $w_i \in \mathbb{F}_p \setminus \{0\}$ for $i \in I$, and suppose that $\alpha_n\abs*{I} \geq 4$. Then, for every $1 \leq j \leq n+1$ and $r \in \mathbb{F}_p$,
    \begin{equation*}
        \abs*{\Prob\paren*{\sum_{i \in I} w_i(X_j)_i = r} - \frac{1}{p}} \leq \frac{2}{\sqrt{\alpha_n\abs*{I}}}.
    \end{equation*}
\end{lemma}

\begin{proof}
    By the definition of the affine Bernoulli entries,
    \begin{equation}
        \label{eq:AffineLinearForm}
        \Prob\paren*{\sum_{i \in I} w_i(X_j)_i = r} = \Prob\paren*{\sum_{i \in I} s_{ij}w_i\xi_{ij} = r - \sum_{i \in I} a_{ij}w_i}.
    \end{equation}
    The coefficients $s_{ij}w_i$ are all nonzero. The result therefore follows from the finite-field Erd\H{o}s--Littlewood--Offord estimate of Nguyen and Wood \cite[Theorem~6.3]{Nguyen2022}.
\end{proof}

We next record the Fourier notation and the effect of the affine transformations on the Fourier coefficients; see \cite{Terras1999} for background on Fourier analysis over $\mathbb{F}_p$. For $t \in \mathbb{F}_p$, let
    \begin{equation*}
        e_p(t) \coloneqq \exp(2\pi i t/p), \qquad \phi_p(t) \coloneqq (1-\alpha_n) + \alpha_n e_p(t).
    \end{equation*}
Equivalently, $\phi_p(t) = \E e_p(t\xi)$ for $\xi \sim \Ber(\alpha_n)$. Thus, $\phi_p$ is the Fourier transform of the Bernoulli distribution on $\mathbb{F}_p$. The Fourier transform of the affine Bernoulli variable $a + s\xi$ is related to $\phi_p$ as follows. For every $a,t \in \mathbb{F}_p$ and $s \in \{\pm1\}$,
\begin{equation}
    \label{eq:AffineFourier}
    \E e_p\paren*{t(a + s\xi)} = e_p(at)\phi_p\paren*{st}, \qquad \abs*{\E e_p\paren*{t(a + s\xi)}} = \abs*{\phi_p(t)}.
\end{equation}
Indeed, the translation by $a$ contributes only the phase $e_p(at)$, while multiplication by $s \in \{\pm1\}$ either leaves $\phi_p(t)$ unchanged or replaces it by
\begin{equation*}
    \phi_p(-t) = \overline{\phi_p(t)}.
\end{equation*}
Thus, neither operation changes the absolute value of the Fourier coefficient. This allows us to deduce the following estimate for our affine Bernoulli vectors.

\begin{lemma}[Fourier estimate]
    \label{lem:FourierEstimate}
    Let $V \subseteq \mathbb{F}_p^n$ be a linear subspace of codimension $r$. Then
    \begin{equation}
        \label{eq:FourierEstimate}
        \max_{1 \leq j \leq n+1} \abs*{\Prob\paren*{X_j \in V} - p^{-r}} \leq p^{-r} \sum_{w \in V^\perp \setminus \{0\}} \prod_{i=1}^n \abs*{\phi_p(w_i)}.
    \end{equation}
    The right-hand side is independent of the deterministic matrices $A_n$ and $S_n$.
\end{lemma}

\begin{proof}
    By Fourier inversion on $\mathbb{F}_p^n$ applied to the indicator function of $V$ and then taking expected values, we get
    \begin{equation*}
        \Prob\paren*{X_j \in V} = p^{-r} \sum_{w \in V^\perp} \E e_p(w \cdot X_j)
    \end{equation*}
    for every $1 \leq j \leq n+1$. The term corresponding to $w = 0$ is $p^{-r}$. For $w \neq 0$, independence of the coordinates and \eqref{eq:AffineFourier} yield
    \begin{equation*}
        \abs*{\E e_p(w \cdot X_j)} = \prod_{i=1}^n \abs*{\E e_p(w_i(X_j)_i)} = \prod_{i=1}^n \abs*{\phi_p(w_i)}.
    \end{equation*}
    Consequently, the triangle inequality implies
    \begin{equation*}
        \abs*{\Prob\paren*{X_j \in V} - p^{-r}} \leq p^{-r} \sum_{w \in V^\perp \setminus \{0\}} \prod_{i=1}^n \abs*{\phi_p(w_i)}.
    \end{equation*}
    Since the right-hand side is independent of $j$, taking the maximum over $1 \leq j \leq n+1$ proves the claim.
\end{proof}

We first exclude spans with a sparse normal vector.

\begin{definition}[Sparse subspaces]
    Fix a sufficiently small constant $\delta_0 > 0$. A linear subspace $V \subseteq \mathbb{F}_p^n$ is called \emph{sparse} if $V^\perp$ contains a nonzero vector with at most $\delta_0n$ nonzero coordinates.
\end{definition}

The following lemma shows that the spans arising during the final column exposures are very unlikely to be sparse. It is an affine analogue of \cite[Lemma~8.6]{Nguyen2022}. The proof combines Odlyzko's lemma with the support argument used in \cite[Proof of Lemma~8.22]{Nguyen2022}.

\begin{lemma}[{Sparse normal vectors, cf. \cite[Lemma~8.6]{Nguyen2022}}]
    \label{lem:SparseNormals}
    There are constants $c,C,\eta > 0$ such that, whenever $(1-\eta)n \leq m \leq n$, we have
    \begin{equation*}
        \Prob\paren*{W_m \text{ is sparse}} \leq Ce^{-c\alpha_n n},
    \end{equation*}
    uniformly in $p$ and in the deterministic matrices $A_n$ and $S_n$.
\end{lemma}

\begin{proof}
    Suppose that $W_m$ is sparse, and choose a nonzero vector $w \in W_m^\perp$ with minimal support. Abbreviate $\sigma \coloneqq \supp(w)$ and $t \coloneqq \abs*{\sigma} \leq \delta_0 n$. Consider the restrictions of the exposed columns to the coordinates in $\sigma$, and put $H \coloneqq \operatorname{span}_{\mathbb{F}_p}(X_1\vert_\sigma,\ldots,X_m\vert_\sigma) \subseteq \mathbb{F}_p^\sigma$.

    Since $w\vert_\sigma$ is orthogonal to $H$, we have $\dim H \leq t-1$. In fact, $\dim H = t-1$. Otherwise $H^\perp$ would have dimension at least two and would contain a nonzero vector vanishing in at least one coordinate of $\sigma$. Extending this vector by zero outside $\sigma$ would produce a nonzero element of $W_m^\perp$ with strictly smaller support, contradicting the choice of $w$.

    Consequently, some $t-1$ of the restricted columns span $H$. Fixing their indices and conditioning on these columns, every remaining restricted column must belong to the resulting hyperplane $H$. We split the possible values of $t$ into two ranges.
	
    Suppose first that $t \leq 144/\alpha_n$. By Odlyzko's lemma, \cref{lem:NonIdenticalOdlyzko}, every remaining restricted column belongs to $H$ with probability at most $1-\alpha_n$. Taking a union bound over the support $\sigma$ and the $t-1$ columns spanning $H$, the probability that $W_m$ is sparse and that its smallest nonzero normal support has size at most $144/\alpha_n$ is therefore at most
    \begin{equation*}
        \sum_{1 \leq t \leq \min\{\delta_0n,144/\alpha_n\}} \binom{n}{t} \binom{m}{t-1} (1-\alpha_n)^{m-t+1} \leq e^{-c\alpha_n n}.
    \end{equation*}
    After taking logarithms, the two binomial factors contribute $O(\alpha_n^{-1}\log{n})$, while the last factor contributes $-\Omega(\alpha_n n)$. We have
    \begin{equation*}
        \alpha_n^2n/\log{n} = n^{2/3 + 2\eps}/\log{n} \to \infty,
    \end{equation*}
    and hence $\alpha_n^{-1}\log{n} = o(\alpha_n n)$. Thus, the exponential decay of the latter factor absorbs the contribution of the binomial coefficients. The additional factor $O(\alpha_n^{-1})$ arising from the sum over $t$ is absorbed in the same way, which proves the claimed bound when $t \leq 144/\alpha_n$.

    Now suppose that $144/\alpha_n < t \leq \delta_0n$. Since $H$ is a hyperplane in $\mathbb{F}_p^\sigma$, the vector $w\vert_\sigma$, whose $t$ coordinates are all nonzero, is a normal vector of $H$. By \cref{lem:ErdosLittlewoodOfford}, every remaining restricted column belongs to $H$ with probability at most $p^{-1}+ 2(\alpha_nt)^{-1/2} \leq 1/2 + 1/6  = 2/3$ as $\alpha_nt > 144$. Taking a union bound over the support $\sigma$ and the $t-1$ columns spanning $H$, the probability that $W_m$ is sparse and that its smallest nonzero normal support has size in this range is therefore at most
    \begin{equation*}
        \sum_{144/\alpha_n < t \leq \delta_0n} \binom{n}{t} \binom{m}{t-1} \paren*{\frac{2}{3}}^{m-t+1} \leq e^{-cn}.
    \end{equation*}
    Indeed, on the logarithmic scale, the product of the two binomial coefficients contributes at most $O(\delta_0n\log(1/\delta_0))$, whereas the last factor contributes at most $-(1-\eta-\delta_0)n\log(3/2)$. Thus, by choosing $\delta_0$ and $\eta$ sufficiently small, the negative contribution dominates. The number of possible values of $t$ is absorbed into the exponential bound.
    Combining the two ranges and adjusting the constants proves the upper bound.
    
    Uniformity follows from \cref{lem:NonIdenticalOdlyzko,lem:ErdosLittlewoodOfford} and $p^{-1}\leq 1/2$.
\end{proof}

Following \cite[Appendix~A]{Nguyen2020} and \cite[Sections~8.8 and 8.11]{Nguyen2022}, we classify the remaining subspaces by their deviations from the uniform containment probability. For Nguyen and Paquette's Appendix~A, we use the numbering of their \href{https://people.math.osu.edu/nguyen.1261/cikk/Sur-CPC.pdf}{author version} throughout the arguments below.

\begin{definition}[Saturation classes]
    Fix a constant $\lambda > 0$, to be chosen below, and let $V \subseteq \mathbb{F}_p^n$ be a non-sparse subspace of codimension $r$. Consider the maximal deviation
    \begin{equation*}
        \max_{1 \leq j \leq n+1} \abs*{\Prob\paren*{X_j \in V} - p^{-r}}.
    \end{equation*}
    We call $V$ \emph{saturated} if this quantity is at most $e^{-\lambda\alpha_n n}$, \emph{semi-saturated} if it is strictly greater than $e^{-\lambda\alpha_n n}$ and at most $16p^{-r}$, and \emph{unsaturated} if it exceeds both $e^{-\lambda\alpha_n n}$ and $16p^{-r}$.
\end{definition}

The maximum over $j$ accounts for the fact that the column distributions need not be identical. These three classes exhaust all non-sparse subspaces. By definition, saturated subspaces already have containment probabilities within $e^{-\lambda\alpha_n n}$ of $p^{-r}$. The next two lemmas show that the random span $W_m$ is unlikely to be semi-saturated or unsaturated.

For semi-saturated spans, compare \cite[Proposition~A.13]{Nguyen2020} and \cite[Lemma~8.10]{Nguyen2022}.

\begin{lemma}[Semi-saturated spans]
    \label{lem:Semisaturated}
    After decreasing $\eta$ if necessary, one can choose $\lambda > 0$ and constants $c,C > 0$ such that, whenever $(1-\eta)n \leq m \leq n$, we have
    \begin{equation*}
        \Prob\paren*{W_m \text{ is semi-saturated and } \operatorname{codim} W_m \leq \eta n} \leq Ce^{-cn},
    \end{equation*}
    uniformly in $p$ and in the deterministic matrices $A_n$ and $S_n$.
\end{lemma}

\begin{proof}
	Our first goal is to upper bound the number of semi-saturated subspaces with codimension $r$.
    Choose $\theta > 0$ so small that $17e \theta < 1$. We use the Fourier-capture argument in the proof of Nguyen and Paquette's Lemma~A.12 \cite{Nguyen2020}; see also the version of Nguyen and Wood \cite[Lemmas~8.9 and~8.10]{Nguyen2022} and \cite[Proposition~2.5]{Maples2010}. Applied with $\delta = \delta_0$, this argument provides a constant $\lambda \in (0,1)$ and a deterministic set $R \subseteq \mathbb{F}_p^n \setminus \{0\}$ such that $\abs*{R} \leq \theta^n p^n$ and $R$ contains every vector $w$ satisfying $\abs*{\supp(w)} > \delta_0 n$ and $\prod_{i=1}^n \abs*{\phi_p(w_i)} > e^{-\lambda\alpha_n n}$. To see this, the assumed product bound implies
    \begin{equation*}
        \sum_{i=1}^{n} \paren*{1 - \abs*{\phi_p(w_i)}^2} \leq -2\log\paren*{\prod_{i=1}^{n} \abs*{\phi_p(w_i)}} < 2\lambda\alpha_n n.
    \end{equation*}
    Consequently, all but at most $\delta_0 n/2$ coordinates satisfy $\abs*{\phi_p(w_i)} \geq 1 - 4\lambda\alpha_n/\delta_0$. Since $\abs*{\supp(w)} > \delta_0 n$, at least $\delta_0 n/2$ of these coordinates are nonzero. The Fourier level-set estimate and the subsequent counting argument in that proof of Lemma~A.12 bound the number of all vectors with these properties by $\theta^n p^n$, provided $\lambda$ is sufficiently small. Thus, $R$ can be chosen to contain every such non-sparse vector, independently of any particular subspace $V$.

    Let $V$ be a semi-saturated subspace of codimension $r$. Its defining lower bound and \cref{lem:FourierEstimate} give
    \begin{equation*}
        e^{-\lambda\alpha_n n} < \max_{1 \leq j \leq n+1} \abs*{\Prob\paren*{X_j \in V} - p^{-r}} \leq p^{-r} \sum_{w \in V^\perp \setminus \{0\}} \prod_{i=1}^n \abs*{\phi_p(w_i)}.
    \end{equation*}
    Since $V$ has codimension $r$, the sum contains fewer than $p^r$ terms. By the pigeonhole principle, there is a nonzero $w \in V^\perp$ satisfying $\prod_{i=1}^n \abs*{\phi_p(w_i)} > e^{-\lambda\alpha_n n}$. Moreover, $V$ is non-sparse by definition, so every nonzero vector in $V^\perp$ has more than $\delta_0 n$ nonzero coordinates. Consequently, $w \in R$, and hence $V^\perp \cap R \neq \varnothing$. Finally, \eqref{eq:AffineFourier} shows that the Fourier weight $\prod_{i=1}^n \abs*{\phi_p(w_i)}$ is independent of the deterministic matrices $A_n$ and $S_n$. Thus, the same set $R$ captures all column distributions simultaneously.

    We now follow the counting argument of Nguyen and Paquette \cite[Proposition~A.13]{Nguyen2020}. Fix $1 \leq r \leq \eta n$. Every semi-saturated subspace $V$ of codimension $r$ satisfies $V^\perp \cap R \neq \varnothing$. For a fixed $w \in R$, the $r$-dimensional subspaces $V^\perp$ containing $w$ correspond, after quotienting by $\operatorname{span}(w)$, to all $(r-1)$-dimensional subspaces of $\mathbb{F}_p^{n-1}$. Their number is at most
    \begin{equation*}
        \prod_{i=0}^{r-2} \frac{p^{n-1-i}-1}{p^{r-1-i}-1} \leq Cp^{(r-1)(n-r)}.
    \end{equation*}
    Since every semi-saturated subspace is captured by at least one $w \in R$, summing this bound over $w \in R$ shows that the number of semi-saturated subspaces of codimension $r$ is at most $C\abs*{R}p^{(r-1)(n-r)} \leq C\theta^np^{nr-r^2+r}$.

    We now convert this counting estimate into a probability bound. For every semi-saturated subspace $V$ of codimension $r$, semi-saturation gives $\Prob\paren*{X_j \in V} \leq 17p^{-r}$ for every $1 \leq j \leq n+1$. Moreover, if $W_m = V$, then $X_1,\ldots,X_m \in V$. Hence, by a union bound over the semi-saturated subspaces of codimension $r$ and by independence of the columns,
    \begin{align*}
        \Prob\paren*{W_m \text{ is semi-saturated and } \operatorname{codim} W_m = r} &\leq \sum_{\substack{V \text{ semi-saturated} \\ \operatorname{codim} V = r}} \Prob\paren*{X_1,\ldots,X_m \in V} \\
        &\leq C\theta^np^{nr-r^2+r}(17p^{-r})^m.
    \end{align*}

    Put $k = n-m$. Since $W_m$ is spanned by $m$ vectors, its codimension $r$ satisfies $r \geq k$. If $r \geq k+1$, then $r(k-r+1) \leq 0$, and hence
    \begin{equation*}
        \Prob\paren*{W_m \text{ is semi-saturated and } \operatorname{codim} W_m = r} \leq C\theta^n17^m \leq C(17\theta)^n.
    \end{equation*}
    In the borderline case $r = k \geq 1$, semi-saturation gives $e^{-\lambda\alpha_n n} < 16p^{-r}$ and hence $p^r < 16e^{\lambda\alpha_n n}$. Therefore, after adjusting $C$,
    \begin{equation*}
        \Prob\paren*{W_m \text{ is semi-saturated and } \operatorname{codim} W_m = r} \leq C\theta^n17^m e^{\lambda\alpha_n n} \leq C(17e\theta)^n.
    \end{equation*}
    The case $r = k = 0$ is impossible because $\mathbb{F}_p^n$ is saturated. Choosing $\theta > 0$ so that $17e\theta < 1$, and summing over the at most $n$ possible values of $r$, proves the claimed bound.
\end{proof}

We control unsaturated subspaces using the symmetrization and swapping techniques of \cite[Appendix~A]{Nguyen2020} and \cite[Section~8.11]{Nguyen2022}.

\begin{lemma}[{Unsaturated spans, cf. \cite[Proposition~A.17]{Nguyen2020} and \cite[Lemma~8.14]{Nguyen2022}}]
    \label{lem:Unsaturated}
    The constants $\lambda$ and $\eta$ in \cref{lem:Semisaturated} may be decreased, if necessary, so that there are constants $c,C > 0$ such that, whenever $(1-\eta)n \leq m \leq n$, we have
    \begin{equation*}
        \Prob\paren*{W_m \text{ is unsaturated and } \operatorname{codim} W_m \leq \eta n} \leq Ce^{-cn},
    \end{equation*}
    uniformly in $p$ and in the deterministic matrices $A_n$ and $S_n$.
\end{lemma}

\begin{proof}
    For each fixed unsaturated subspace, we construct an auxiliary family that gives a pointwise swapping estimate. We then sum these estimates over the subspaces.

    \emph{Reduction to a fixed subspace.}
    Since $W_m$ takes a unique value in the finite set of subspaces of $\mathbb{F}_p^n$, we have the disjoint decomposition
    \begin{equation}
        \label{eq:UnsaturatedDecomposition}
        \Prob\paren*{W_m \text{ is unsaturated and } \operatorname{codim} W_m \leq \eta n} = \sum_{\substack{V \text{ unsaturated} \\ \operatorname{codim} V \leq \eta n}} \Prob\paren*{W_m = V},
    \end{equation}
    where the sum ranges over linear subspaces $V \subseteq \mathbb{F}_p^n$. We first derive a bound for each summand and then sum these bounds using the disjointness of suitable span events.

    Fix an unsaturated linear subspace $V \subseteq \mathbb{F}_p^n$ with $\operatorname{codim} V \leq \eta n$, and put $r \coloneqq \operatorname{codim} V$. Let $j(V)$ be the least index maximizing $\Prob(X_j \in V)$, and write $\rho_V \coloneqq \Prob\paren*{X_{j(V)} \in V}$. Since probabilities are nonnegative, whenever $\Prob(X_j \in V) \leq p^{-r}$, we have
    \begin{equation*}
        \abs*{\Prob\paren*{X_j \in V} - p^{-r}} \leq p^{-r}.
    \end{equation*}
    On the other hand, since $V$ is unsaturated,
    \begin{equation*}
        \max_{1 \leq j \leq n+1} \abs*{\Prob\paren*{X_j \in V} - p^{-r}} > \max\set*{e^{-\lambda\alpha_n n},16p^{-r}}.
    \end{equation*}
    The maximum must therefore be attained by a positive deviation. Consequently,
    \begin{equation}
        \label{eq:UnsaturatedRho}
        \rho_V - p^{-r} > \max\set*{e^{-\lambda\alpha_n n},16p^{-r}}.
    \end{equation}

    \emph{Construction of the comparison vector.}
    We use the symmetrization argument in the proof of Nguyen and Paquette's Lemma~A.15 \cite{Nguyen2020}. They attribute their Lemma~A.15 to Maples \cite[Lemma~2.8]{Maples2010} and use the Halász-type construction from \cite[Proposition~3.6]{Maples2010}, recorded in their Proposition~A.24; see also \cite[Lemma~8.12]{Nguyen2022} for the version of Nguyen and Wood. The arguments provide an $\alpha_n/64$-balanced probability distribution $\nu$ on $\mathbb{F}_p$ whose Fourier transform is nonnegative. If $Y \in \mathbb{F}_p^n$ has independent $\nu$-distributed coordinates, then
    \begin{equation}
        \label{eq:UnsaturatedSymmetrization}
        p^{-r} \sum_{w \in V^\perp \backslash\lbrace 0\rbrace} \prod_{i=1}^n \abs*{\phi_p(w_i)} \leq \paren*{\frac{1}{2} + o(1)}\paren*{\Prob(Y \in V) - p^{-r}}.
    \end{equation}
    To justify this estimate, write
    \begin{equation*}
        f(w) \coloneqq \prod_{i=1}^n \abs*{\phi_p(w_i)}, \qquad g(w) \coloneqq \prod_{i=1}^n \widehat{\nu}(w_i).
    \end{equation*}
    The construction in Nguyen and Paquette's Proposition~A.24 gives $f(w) \leq g(w)^{16}$. For a suitable constant $c > 0$, put $\tau \coloneqq e^{-c\alpha_n n}$. On the range $f(w) \leq \tau$, we therefore have $f(w) \leq \tau^{15/16}g(w)$. For the remaining terms, the level-set comparison in their proof of Lemma~A.15, using that $V$ is non-sparse, gives
    \begin{equation*}
        \sum_{\substack{w \in V^\perp \setminus \{0\} \\ f(w) > \tau}} f(w) \leq \frac{1}{2} \sum_{w \in V^\perp \setminus \{0\}} g(w).
    \end{equation*}
    Combining the two ranges yields
    \begin{equation*}
        \sum_{w \in V^\perp \setminus \{0\}} f(w) \leq \paren*{\frac{1}{2} + \tau^{15/16}} \sum_{w \in V^\perp \setminus \{0\}} g(w).
    \end{equation*}
    Finally, Fourier inversion gives $p^{-r} \sum_{w \in V^\perp \setminus \{0\}} g(w) = \Prob(Y \in V) - p^{-r}$. Since $\alpha_n n \to \infty$, we have $\tau^{15/16} = o(1)$, which proves the claimed estimate. The comparison distribution constructed in Nguyen and Paquette's Proposition~A.24 depends only on the moduli of the Fourier coefficients of the entry distribution. Hence, \eqref{eq:AffineFourier} shows that, for fixed $n,p$, this single choice of $\nu$ works simultaneously for all unsaturated subspaces $V$ considered here and all deterministic matrices $A_n$ and $S_n$.

    Combining \cref{lem:FourierEstimate} with \eqref{eq:UnsaturatedSymmetrization}, and taking $n$ sufficiently large such that $\frac{1}{2} + o(1) \leq \frac{3}{5}$, gives
    \begin{equation}
        \label{eq:UnsaturatedGain}
        \Prob\paren*{Y \in V} \geq p^{-r} + \frac{5}{3}\paren*{\rho_V - p^{-r}} > \frac{3}{2}\rho_V.
    \end{equation}
    Thus, the vector $Y$ is at least $3/2$ times as likely to belong to $V$ as any of the original columns $X_1,\ldots,X_{n+1}$. The swapping argument below amplifies this fixed advantage over linearly many independent copies of $Y$, producing the required exponential decay.

    \emph{Construction of an auxiliary independent family.}
    We now amplify the gain in \eqref{eq:UnsaturatedGain} using linearly many independent copies of $Y$. We supplement these vectors with independent copies of $X_{j(V)}$ so that the resulting auxiliary family contains
    \begin{equation*}
        a = m - \Delta n + O(1)
    \end{equation*}
    vectors, where $\Delta > 0$ is a constant that we will specify below. Let us motivate the parameter choice: If the auxiliary vectors are linearly independent in $V$ and the original columns $X_1,\ldots,X_m$ span $V$, at most $m - a = \Delta n + O(1)$ of those original columns are needed to complete the auxiliary vectors to a generating system of $V$. The corresponding index set can be chosen in $\binom{m}{m-a}$ ways, and choosing $\Delta$ sufficiently small will ensure that this binomial factor is dominated by the exponential gain obtained from the independent copies of $Y$.

    To implement this construction, fix a constant $0< \delta_1 < 1/3$. Choose $\Delta > 0$ sufficiently small, to be specified at the end of the proof, and put $\delta_2 \coloneqq \delta_1 + \Delta < 1/3$. Decrease $\eta$ and $\lambda$, if necessary, so that $\eta < \Delta/2$ and $\lambda < \Delta/256$. Set
    \begin{equation*}
        r_0 \coloneqq \floor*{\delta_1n}, \qquad s \coloneqq m - \floor*{\delta_2n}, \qquad a \coloneqq r_0 + s.
    \end{equation*}
    Here $r_0$ is the number of independent copies of $Y$ and $s$ is the number of additional independent copies of $X_{j(V)}$. By construction, $m - a = \floor*{\delta_2n} - \floor*{\delta_1n} = \Delta n + O(1)$, and hence
    \begin{equation*}
        n - a = n - m + m - a \geq \Delta n - O(1).
    \end{equation*}
    Moreover, since $\operatorname{codim} V \leq \eta n$ and $\eta < \Delta/2$, we have $\dim V \geq a$ for all sufficiently large $n$.

    Let $Y_1,\ldots,Y_{r_0}$ be independent copies of $Y$, and let $Z_1,\ldots,Z_s$ be independent copies of $X_{j(V)}$. Take all these auxiliary vectors to be mutually independent and jointly independent of the original columns $X_1,\ldots,X_{n+1}$. Consider the event that all auxiliary vectors belong to $V$ and are linearly independent, and denote this event by $\calA_V$. We estimate its probability by exposing the auxiliary vectors $Y_1,\ldots,Y_{r_0},Z_1,\ldots,Z_s$ successively.

    Suppose that some of these vectors have already been exposed and are linearly independent in $V$, and let $H \subseteq V$ denote their span. If we condition on the exposed vectors, then the subspace $H$ is fixed and the next auxiliary vector remains independent and its underlying distribution does not change. The next vector extends the current linearly independent family inside $V$ precisely when it belongs to $V \setminus H$. Since $H \subseteq V$, this probability equals $\Prob\paren*{Y_i \in V \setminus H} = \Prob\paren*{Y_i \in V} - \Prob\paren*{Y_i \in H}$ if the next vector is $Y_i$, and analogously for $Z_i$.

    At every stage, we have $\dim H \leq a-1$. Since $Y$ is $\alpha_n/64$-balanced and each $Z_i$ is $\alpha_n$-balanced, \cref{lem:NonIdenticalOdlyzko} gives
    \begin{equation*}
        \Prob\paren*{Y_i \in H} \leq (1 - \alpha_n/64)^{n-a+1}, \qquad \Prob\paren*{Z_i \in H} \leq (1 - \alpha_n)^{n-a+1}.
    \end{equation*}
    By \eqref{eq:UnsaturatedRho}, we have $\rho_V > e^{-\lambda\alpha_n n}$. Since $n-a \geq \Delta n - O(1)$ and $\lambda < \Delta/256$, it follows that, for some constant $c_\star > 0$,
    \begin{equation}
        \label{eq:UnsaturatedAvoidance}
        \max\set*{\Prob\paren*{Y_i \in H},\Prob\paren*{Z_i \in H}} \leq \rho_V e^{-c_\star\alpha_n n}.
    \end{equation}
    On the other hand, \eqref{eq:UnsaturatedGain} gives
    \begin{equation*}
        \Prob\paren*{Y_i \in V} \geq \frac{3}{2}\rho_V, \qquad \Prob\paren*{Z_i \in V} = \rho_V.
    \end{equation*}
    Consequently, each $Y_i$ extends the current linearly independent family inside $V$ with probability at least $\rho_V(3/2 - e^{-c_\star\alpha_n n})$ and each $Z_i$ does so with probability at least $\rho_V(1 - e^{-c_\star\alpha_n n})$. Multiplying these successive conditional probabilities yields, for all sufficiently large $n$,
    \begin{equation}
        \label{eq:UnsaturatedIndependent}
        \Prob\paren*{\calA_V} \geq \rho_V^a \paren*{\frac{3}{2} - e^{-c_\star\alpha_n n}}^{r_0} \paren*{1 - e^{-c_\star\alpha_n n}}^s\geq (3/2)^{r_0-1} \rho_V^a.
    \end{equation}
    The factorization $(3/2 - e^{-c_\star\alpha_n n})^{r_0} = (3/2)^{r_0}(1-2/3\cdot e^{-c_\star\alpha_n n})^{r_0}$ leaves the correction product $(1-2/3\cdot e^{-c_\star\alpha_n n})^{r_0}(1-e^{-c_\star\alpha_n n})^s$. Since $r_0+s\leq n$ and $n e^{-c_\star\alpha_n n}=o(1)$, this product is $1-o(1)$ and hence at least $2/3$ for all sufficiently large $n$. This proves the second inequality.

    \emph{The pointwise swapping estimate.}
    We now use \eqref{eq:UnsaturatedIndependent} to derive a pointwise bound for $\Prob\paren*{W_m = V}$. 
    For every $I \subseteq \{1,\ldots,m\}$ with $\abs*{I} = m-a$, define the span
    \begin{equation*}
        S_I \coloneqq \operatorname{span}_{\mathbb{F}_p}(Y_1,\ldots,Y_{r_0},Z_1,\ldots,Z_s,(X_i)_{i \in I}).
    \end{equation*}
    On the event $\{W_m = V\} \cap \calA_V$, the $a$ auxiliary vectors are linearly independent in $V$, while the original columns $X_1,\ldots,X_m$ span $V$. Thus, some of the original columns complete the auxiliary vectors to a basis of $V$. This requires at most $m-a$ columns. By adding further original columns if necessary, we may therefore choose a set $I \subseteq \{1,\ldots,m\}$ with $\abs*{I} = m-a$ such that $S_I = V$. Consequently,
    \begin{equation*}
        \{W_m = V\} \cap \calA_V \subseteq \bigcup_{\substack{I \subseteq \{1,\ldots,m\} \\ \abs*{I} = m-a}} \paren*{\{S_I = V\} \cap \bigcap_{l \notin I} \{X_l \in V\}}.
    \end{equation*}
    The event $\calA_V$ depends only on the auxiliary vectors and is therefore independent of $W_m$. Taking probabilities and applying a union bound over $I$, we obtain
    \begin{equation*}
        \Prob\paren*{W_m = V}\Prob\paren*{\calA_V} \leq \sum_{\substack{I \subseteq \{1,\ldots,m\} \\ \abs*{I} = m-a}} \Prob\paren*{S_I = V} \prod_{l \notin I} \Prob\paren*{X_l \in V}.
    \end{equation*}
    Indeed, $S_I$ depends only on the auxiliary vectors and on $(X_i)_{i \in I}$, which are independent of the remaining original columns. By the definition of $\rho_V$, $\Prob\paren*{X_l \in V} \leq \rho_V$ for every $1 \leq l \leq m$. Moreover, exactly $a$ indices lie outside $I$. Thus,
    \begin{equation*}
        \Prob\paren*{W_m = V}\Prob\paren*{\calA_V} \leq \rho_V^a \sum_{\substack{I \subseteq \{1,\ldots,m\} \\ \abs*{I} = m-a}} \Prob\paren*{S_I = V}.
    \end{equation*}
    Combining this with \eqref{eq:UnsaturatedIndependent} gives the pointwise estimate
    \begin{equation}
        \label{eq:UnsaturatedSwapping}
        \Prob\paren*{W_m = V} \leq (3/2)^{-r_0+1} \sum_{\substack{I \subseteq \{1,\ldots,m\} \\ \abs*{I} = m-a}} \Prob\paren*{S_I = V}.
    \end{equation}

    \emph{Summation over the subspaces.}
    It remains to sum \eqref{eq:UnsaturatedSwapping} over $V$. Since the distribution of the auxiliary vectors $Z_1,\ldots,Z_s$ depends on the maximizing index $j(V)$, we group the subspaces according to this index. For $j \in \{1,\ldots,n+1\}$, write $S_I^{(j)}$ for the span $S_I$ constructed using independent copies of $X_j$ in place of $Z_1,\ldots,Z_s$. Thus, for the fixed subspace $V$, the random subspace $S_I$ in \eqref{eq:UnsaturatedSwapping} has the same distribution as $S_I^{(j(V))}$. Using \eqref{eq:UnsaturatedDecomposition} and \eqref{eq:UnsaturatedSwapping}, the probability that $W_m$ is unsaturated and $\operatorname{codim} W_m \leq \eta n$ is at most
    \begin{equation*}
        (3/2)^{-r_0+1} \sum_{j=1}^{n+1} \sum_{\substack{I \subseteq \{1,\ldots,m\} \\ \abs*{I} = m-a}} \sum_{\substack{V \text{ unsaturated} \\ \operatorname{codim} V \leq \eta n \\ j(V) = j}} \Prob\paren*{S_I^{(j)} = V} \leq (n+1) \binom{m}{m-a} (3/2)^{-r_0+1}.
    \end{equation*}
    Indeed, for fixed $j$ and $I$, the events $\{S_I^{(j)} = V\}$ are disjoint as $V$ varies, so the innermost sum is at most one. There are $n+1$ choices for $j$ and $\binom{m}{m-a}$ choices for $I$.

    Finally, $m-a = \Delta n + O(1)$ and $r_0 = \delta_1 n + O(1)$. Since $m \leq n$, the standard estimate $\binom{m}{k} \leq (em/k)^k$ gives
    \begin{equation*}
        \binom{m}{m-a} \leq \exp\paren*{\paren*{\Delta\log{e/\Delta} + o(1)}n}.
    \end{equation*}
    Choose $\Delta > 0$ sufficiently small that $\Delta\log(e/\Delta) < \frac{\delta_1}{2}\log(3/2)$. The exponential decay from $(3/2)^{-r_0}$ then dominates both the binomial coefficient and the factor $n+1$. Consequently, there are constants $c,C > 0$ such that
    \begin{equation*}
        \Prob\paren*{W_m \text{ is unsaturated and } \operatorname{codim} W_m \leq \eta n} \leq Ce^{-cn}.
    \end{equation*}
\end{proof}

\begin{corollary}[{cf. \cite[Theorems~5.2 and~8.2]{Nguyen2022}}]
    \label{cor:RandomSpanEstimate}
    Let $\calG_m$ be the event that $W_m$ is saturated or $\operatorname{codim} W_m > \eta n$. After decreasing $\eta$ if necessary, there are constants $c,C > 0$ such that, whenever $(1-\eta)n \leq m \leq n$,
    \begin{equation*}
        \Prob\paren*{\calG_m^c} \leq Ce^{-c\alpha_n n}.
    \end{equation*}
    Moreover, on $\calG_m$, if $r \coloneqq \operatorname{codim} W_m \leq \eta n$, then
    \begin{equation*}
        \abs*{\Prob\paren*{X_{m+1} \in W_m \mid X_1,\ldots,X_m} - p^{-r}} \leq e^{-\lambda\alpha_n n}.
    \end{equation*}
    The estimates are uniform in $p$ and in the deterministic matrices $A_n$ and $S_n$.
\end{corollary}

\begin{proof}
    If $\calG_m^c$ occurs, then $\operatorname{codim} W_m \leq \eta n$ and $W_m$ is not saturated. Thus, $W_m$ is either sparse or non-sparse and not saturated. In the latter case, it is semi-saturated or unsaturated. A union bound together with \cref{lem:SparseNormals,lem:Semisaturated,lem:Unsaturated} therefore gives, after adjusting the constants,
    \begin{equation*}
        \Prob\paren*{\calG_m^c} \leq Ce^{-c\alpha_n n}.
    \end{equation*}
    Now suppose that $\calG_m$ occurs and that $r \coloneqq \operatorname{codim} W_m \leq \eta n$. By the definition of $\calG_m$, the subspace $W_m$ is saturated. Since $X_{m+1}$ is independent of the exposed columns, conditioning on $X_1,\ldots,X_m$ fixes $W_m$ without changing the distribution of $X_{m+1}$. The claimed estimate therefore follows directly from the definition of saturation. The constants are uniform by \cref{lem:SparseNormals,lem:Semisaturated,lem:Unsaturated}.
\end{proof}

\begin{lemma}[{cf. \cite[Theorem~5.1]{Nguyen2022}}]
    \label{lem:FixedPrimeRank}
    There are constants $c_0,C > 0$ such that
    \begin{equation*}
        \Prob\paren*{\rank_{\mathbb{F}_p}[X_1 \ \cdots \ X_{n+1}] < n} \leq C\paren*{p^{-2} + e^{-c_0\alpha_n n}},
    \end{equation*}
    uniformly in $p$ and in the deterministic matrices $A_n$ and $S_n$.
\end{lemma}

\begin{proof}
    We expose a large independent block, then compare the remaining rank process with that of uniform vectors over $\mathbb{F}_p$.

    \emph{Construction of a large independent block.}
    Let $\eta$ be as in \cref{cor:RandomSpanEstimate}. Set $d \coloneqq \floor*{\eta n}$ and $N \coloneqq n-d$. We first show that $X_1,\ldots,X_N$ are linearly independent with exponentially high probability.

    Suppose that $X_1,\ldots,X_{j-1}$ are linearly independent. Their span then has dimension $j-1$ and hence codimension $n-j+1$. Odlyzko's lemma, \cref{lem:NonIdenticalOdlyzko}, gives
    \begin{equation*}
        \Prob\paren*{X_j \in \operatorname{span}_{\mathbb{F}_p}(X_1,\ldots,X_{j-1}) \mid X_1,\ldots,X_{j-1}} \leq (1-\alpha_n)^{n-j+1}
    \end{equation*}
    If $X_1,\ldots,X_N$ are linearly dependent, there must be a first column that belongs to the span of its predecessors. A union bound over this first dependence therefore yields
    \begin{equation}
        \label{eq:FixedPrimeIndependent}
        \Prob\paren*{X_1,\ldots,X_N \text{ are linearly dependent}} \leq \sum_{j=1}^N (1-\alpha_n)^{n-j+1} \leq \alpha_n^{-1}(1-\alpha_n)^{d+1} \leq e^{-c\alpha_n n}
    \end{equation}
    for some constant $c > 0$ and all sufficiently large $n$. In the last step, we used $d = \eta n + O(1)$ and $\alpha_n n/\log{n} = n^{5/6 + \eps}/\log{n} \to \infty$.

    Consequently, with exponentially high probability, the first $N$ columns form a linearly independent block. On the event that $X_1,\ldots,X_N$ are linearly independent, $\dim W_N = N = n-d$ and hence $\operatorname{codim} W_N = d$. Thus, the remaining exposure process starts in the small-codimension regime covered by \cref{cor:RandomSpanEstimate}.

    \emph{Comparison with the uniform rank process.}
    For $N \leq m \leq n+1$, let $R_m \coloneqq \operatorname{codim} W_m$. On the event that $X_1,\ldots,X_N$ are linearly independent, we have $R_N = d$ and $R_m \leq d \leq \eta n$ for every $m \geq N$. Moreover, for $N \leq m \leq n$, if $R_m = r$, then
    \begin{equation*}
        R_{m+1} =
        \begin{cases}
            r, & \text{if } X_{m+1} \in W_m, \\
            r-1, & \text{if } X_{m+1} \notin W_m. \\
        \end{cases}
    \end{equation*}

    By \cref{cor:RandomSpanEstimate} and a union bound,
    \begin{equation}
        \label{eq:FixedPrimeExceptional}
        \Prob\paren*{\calG_m^c \text{ for some } N \leq m \leq n} \leq Cn e^{-c\alpha_n n} \leq C e^{-c_1\alpha_n n}
    \end{equation}
    for some $c_1 > 0$. Whenever $X_1,\ldots,X_N$ are linearly independent and $\calG_m$ occurs and $R_m = r$, \cref{cor:RandomSpanEstimate} gives
    \begin{equation*}
        \abs*{\Prob\paren*{R_{m+1} = r \mid X_1,\ldots,X_m} - p^{-r}} \leq e^{-\lambda\alpha_n n}.
    \end{equation*}
    Equivalently, the conditional probability that the codimension decreases from $r$ to $r-1$ differs from $1-p^{-r}$ by at most $e^{-\lambda\alpha_n n}$.

    We now compare $(R_m)_{m=N}^{n+1}$ with the Markov chain $(D_t)_{t=0}^{d+1}$ defined by $D_0 = d$ and
    \begin{equation*}
        D_{t+1} =
        \begin{cases}
            r, & \text{with probability } p^{-r}, \\
            r-1, & \text{with probability } 1-p^{-r}, \\
        \end{cases}
    \end{equation*}
    whenever $D_t = r \geq 1$, and, once it hits $0$, it remains $0$. Although $(R_m)_{m=N}^{n+1}$ need not itself be a Markov chain, the preceding estimate holds conditional on the entire exposure history and therefore permits a step-by-step coupling with $(D_t)$. We run this coupling as long as the two codimensions agree and the relevant event $\calG_m$ occurs. At each of the $d+1$ transitions, the coupling fails with probability at most $e^{-\lambda\alpha_n n}$. Combining \eqref{eq:FixedPrimeIndependent} and \eqref{eq:FixedPrimeExceptional} with the step-by-step coupling, we obtain
    \begin{equation}
        \label{eq:FixedPrimeCoupling}
        \Prob\paren*{R_{n+1} > 0} \leq \Prob\paren*{D_{d+1} > 0} + Ce^{-c_1\alpha_n n} + (d+1)e^{-\lambda\alpha_n n} \leq \Prob(D_{d+1} > 0) + Ce^{-c_2\alpha_n n}
    \end{equation}
    for some $c_2 > 0$. Here the final inequality follows from $d \leq n$ and $\alpha_n n/\log{n} \to \infty$.

    \emph{Evaluation of the comparison chain.}
    The chain $(D_t)_{t=0}^{d+1}$ is precisely the process obtained by successively exposing $d+1$ independent uniform vectors in $\mathbb{F}_p^d$. Indeed, if their current span has codimension $r$, then the next uniform vector belongs to this span with probability $p^{-r}$. Consequently, $D_{d+1} = 0$ precisely when the corresponding uniform $d \times (d+1)$ matrix has full row rank. The standard finite-field rank formula therefore gives
    \begin{equation*}
        \Prob\paren*{D_{d+1} = 0} = \prod_{k=2}^{d+1} (1-p^{-k}).
    \end{equation*}
    It follows that
    \begin{equation}
        \label{eq:FixedPrimeUniformChain}
        \Prob\paren*{D_{d+1} > 0} = 1 - \prod_{k=2}^{d+1} (1-p^{-k}) \leq \sum_{k=2}^{d+1} p^{-k} \leq \sum_{k=2}^\infty p^{-k} \leq 2p^{-2}.
    \end{equation}
    Combining \eqref{eq:FixedPrimeCoupling} and \eqref{eq:FixedPrimeUniformChain}, we obtain
    \begin{equation*}
        \Prob\paren*{R_{n+1} > 0} \leq 2p^{-2} + Ce^{-c_2\alpha_n n} \leq C\paren*{p^{-2} + e^{-c_0\alpha_n n}},
    \end{equation*}
    where we may take $c_0 = c_2$ after enlarging $C$. The inequality $R_{n+1} > 0$ corresponds precisely to the event that $[X_1 \ \cdots \ X_{n+1}]$ fails to have full row rank over $\mathbb{F}_p$.
\end{proof}

\begin{proof}[Proof of \cref{prop:ThreePrimeRegimes}, part~\labelcref{item:MediumPrimes}]
    Let $c_0,C>0$ be the constants from \cref{lem:FixedPrimeRank}. The lemma and a union bound imply that the probability under consideration is at most
    \begin{equation*}
        C\sum_{\substack{p > L \\ p \text{ prime}}} p^{-2} + C U_n e^{-c_0\alpha_n n} \leq C\sum_{k > L} k^{-2} + C e^{-c_0\alpha_n n/2} \leq \frac{C}{L} + C e^{-c_0\alpha_n n/2}.
    \end{equation*}
\end{proof}

\subsection{Large Primes}

The error term in \cref{lem:FixedPrimeRank} prevents us from extending the medium-prime union bound to all primes $p \geq U_n$. Following \cite[Section 6, p.~31]{Nguyen2022}, we control them simultaneously by a \emph{watch list}: a prime enters the list when an exposed column first fails to increase the dimension of the span modulo that prime.

In this subsection, $X_i$ denotes an integer column of $M_n$ or its reduction modulo $p$, according to the ambient space.
\begin{definition}
    \label{def:WatchList}
    Let $X_1,\ldots,X_{n+1} \in \Z^n$ denote the columns of $M_n$. For every prime $p \geq U_n$, set $W_0(p) = \{0\}$ and, for $1 \leq k \leq n+1$, let
    \begin{equation*}
        W_k(p) \coloneqq \operatorname{span}_{\mathbb{F}_p}(X_1,\ldots,X_k).
    \end{equation*}
    We say that the $k$-th exposure has a \emph{failed rank increment} modulo $p$ if $X_k \in W_{k-1}(p)$, or equivalently if $\dim W_k(p) = \dim W_{k-1}(p)$. For $1 \leq k \leq n$, define the watch list after $k$ exposures by
    \begin{equation*}
        \mathfrak{W}_k \coloneqq \set*{p \text{ prime} : p \geq U_n \text{ and } \dim W_k(p) \leq k-1}.
    \end{equation*}
    Thus, $\mathfrak{W}_k$ consists precisely of the large primes for which at least one failed rank increment has occurred among the first $k$ exposures.
\end{definition}

If $M_n$ does not have full row rank over $\mathbb{F}_p$, that is, if $\rank_{\mathbb{F}_p} M_n < n$, then at least two of its $n+1$ column exposures fail to increase the dimension. Thus, any such prime enters the watch list during the first $n$ exposures and subsequently has a second failed rank increment.

The watch list itself is small with high probability. Indeed, every prime entering it during the first $n$ exposures divides the determinant of the square matrix formed by those columns. Provided that this determinant is nonzero, Hadamard's inequality bounds its absolute value, while every prime on the list is large; consequently, only a few distinct primes can belong to the watch list.  The next lemma makes this argument precise and establishes the required nonvanishing of the determinant with high probability.

\begin{lemma}[{cf. \cite[Section 6, pp.~31--32]{Nguyen2022}}]
    \label{lem:WatchList}
    The probability that $\det[X_1 \ \cdots \ X_n] = 0$ is $O(\exp\paren*{-c\alpha_n n})$ for some constant $c > 0$, uniformly in $A_n$ and $S_n$. Moreover, whenever this determinant is nonzero, for every $1 \leq k \leq n$,
    \begin{equation}
        \label{eq:WatchList}
        \abs*{\mathfrak{W}_k} \leq C\frac{\log{n}}{c_0\alpha_n},
    \end{equation}
    where $C > 0$ depends only on $T$ and the fixed constants.
\end{lemma}

\begin{proof}
	Fix a prime $\pi_n$ with $\pi_n \geq U_n$.
	To prove that
    \begin{equation*}
        \Prob\paren*{\det[X_1 \ \cdots \ X_n] = 0} = O(\exp\paren*{-c\alpha_n n})
    \end{equation*}
    for some constant $c > 0$, we use \cref{lem:FixedPrimeRank} applied to a prime $\pi_n\geq U_n$. There is the minor distinction that \cref{lem:FixedPrimeRank} is a statement about the probability of $n+1$ vectors failing to span $\mathbb{F}_{\pi_n}^n$, and our claim is just about $n$ vectors. This is not an issue because the proof of \cref{lem:FixedPrimeRank}, with the column exposure stopped after $X_n$ rather than $X_{n+1}$, still gives
	\begin{equation*}
		\Prob\paren*{\rank_{\mathbb{F}_{\pi_n}}[X_1 \ \cdots \ X_n] < n} \leq C\paren*{\pi_n^{-1} + e^{-c_0\alpha_n n}}.
	\end{equation*}
	Let us sketch why this estimate is a consequence of the proof of \cref{lem:FixedPrimeRank}. One can still do the coupling in \eqref{eq:FixedPrimeCoupling} for $n$ vectors with an upper bound depending on $\Prob\paren*{D_{d}  > 0}$. Using the same arguments from the paragraph ``Evaluation of the comparison chain'', one gets 
	\begin{align*}
		\Prob\paren*{D_d > 0} = 1 - \prod_{k=1}^d (1 - \pi^{-k}_n)\leq \sum_{k=1}^d \pi^{-k}_n\leq 2\pi^{-1}_n,
	\end{align*}
	which implies the estimate above.
	
    Now, if $\det[X_1 \ \cdots \ X_n] = 0$, this matrix remains singular after reduction modulo $\pi_n$. Since $\pi_n \geq U_n$, it follows that
    \begin{equation*}
        \Prob\paren*{\det[X_1 \ \cdots \ X_n] = 0} \leq Ce^{-c\alpha_n n}
    \end{equation*}
    for some $c > 0$.

    Suppose now that this determinant is nonzero, and let $p \in \mathfrak{W}_k$. The remaining $n-k$ columns can increase the dimension by at most $n-k$, and hence $\dim W_n(p) \leq n-1$. Consequently, the square matrix formed by the first $n$ columns is singular over $\mathbb{F}_p$, and therefore $p \mid \det[X_1 \ \cdots \ X_n]$. Since every prime $p \in \mathfrak{W}_k$ divides the nonzero determinant and these primes are distinct, their product divides the determinant. Since every entry of $M_n$ has absolute value at most $n^T + 1$, Hadamard's inequality gives
    \begin{equation*}
        \exp\paren*{\abs*{\mathfrak{W}_k}c_0\alpha_n n/2} = U_n^{\abs*{\mathfrak{W}_k}} \leq \prod_{p \in \mathfrak{W}_k} p \leq \abs*{\det[X_1 \ \cdots \ X_n]} \leq n^{n/2}(n^T + 1)^n.
    \end{equation*}
    Taking logarithms and using $\log(n^T + 1) = T\log{n} + O(1)$ proves \eqref{eq:WatchList}.
\end{proof}

Following the arguments in \cite[Section 6]{Nguyen2022}, we next bound the probability that some large prime has a second failed rank increment. Set
\begin{equation*}
    n_0 \coloneqq n - \floor*{\frac{3\log{n}}{\alpha_n}},
\end{equation*}
and distinguish whether some large prime accumulates two failed rank increments by time $n_0$ or afterwards. At the first time $k \leq n_0$ at which any large prime accumulates two failed rank increments, that prime has had exactly one previous failed rank increment, so $\dim W_{k-1}(p) = k-2$. Thus, $W_{k-1}(p)$ has sufficiently large codimension for Odlyzko's lemma, \cref{lem:NonIdenticalOdlyzko}. After time $n_0$, only $n+1-n_0 = O(\log{n}/\alpha_n)$ exposures remain, but their codimension may be too small for this estimate. We handle them using sparse-normal lifting and finite-field anti-concentration in \cref{lem:DenseNormals,lem:LateAnticoncentration}.

\begin{lemma}[{cf. \cite[Section 6, p.~32]{Nguyen2022}}]
    \label{lem:EarlyLosses}
    The probability that $\det[X_1 \ \cdots \ X_n] \neq 0$ and that some large prime accumulates two failed rank increments by time $n_0$ is $O(n^{-1/2})$, uniformly in $A_n$ and $S_n$.
\end{lemma}

\begin{proof}
    On the event $\det[X_1 \ \cdots \ X_n] \neq 0$, \cref{lem:WatchList} gives $\abs*{\mathfrak{W}_j} \leq C\frac{\log{n}}{c_0\alpha_n}$ for every $1 \leq j \leq n$. Consequently, if the event in the statement occurs, there is a first time $2 \leq k \leq n_0$ at which some large prime has accumulated two failed rank increments, and immediately before the $k$-th exposure the watch list has size $\abs*{\mathfrak{W}_{k-1}}\leq C\frac{\log{n}}{c_0\alpha_n}$. 
    
    Fix $2 \leq k \leq n_0$, and condition on an arbitrary realization of $X_1,\ldots,X_{k-1}$ for which $\abs*{\mathfrak{W}_{k-1}} \leq C\frac{\log{n}}{c_0\alpha_n}$ and no large prime has yet had two failed rank increments. Note that we assume $\abs*{\mathfrak{W}_{k-1}} \leq C\frac{\log{n}}{c_0\alpha_n}$ and do not condition on the event that $\det[X_1 \ \cdots \ X_n] \neq 0$ as this would alter the probabilities of the unexposed columns.

    Every $p \in \mathfrak{W}_{k-1}$ has then had exactly one failed rank increment, so $\dim W_{k-1}(p) = k-2$. Since $X_k$ is independent of the previously exposed columns, \cref{lem:NonIdenticalOdlyzko} gives
    \begin{equation*}
        \Prob\paren*{X_k \in W_{k-1}(p) \mid X_1,\ldots,X_{k-1}} \leq (1-\alpha_n)^{n-k+2}.
    \end{equation*}
    Note that in our notation $X_k$ denotes the reduction modulo $p$, whereas $X_1,\ldots,X_{k-1}$ denote the integer columns.
    A union bound over the watch list therefore shows that, uniformly over all such realizations, the conditional probability that the first time at which any large prime has a second failed rank increment equals $k$ is at most $(1-\alpha_n)^{n-k+2}C\frac{\log{n}}{c_0\alpha_n}$. Summing over its possible occurrence times yields
    \begin{equation*}
        C\frac{\log{n}}{c_0\alpha_n} \sum_{k=2}^{n_0} (1-\alpha_n)^{n-k+2} \leq C\frac{\log{n}}{c_0\alpha_n^2} (1-\alpha_n)^{n-n_0+2} = O(n^{-1/2}),
    \end{equation*}
    where the last estimate follows from the definition of $n_0$.
\end{proof}

\cref{lem:EarlyLosses} applies while the codimension of the exposed span is relatively large. To handle the final $O(\log{n}/\alpha_n)$ column exposures, we construct a single high-probability event that controls normal vectors simultaneously for all large primes $p \geq U_n$. To do so, fix $c_1 > 0$ sufficiently small in terms of $c_0$ and $T$, in particular so that $c_1(T + 1) < c_0/2$, and set
\begin{equation*}
    \ell \coloneqq \floor*{\frac{c_1\alpha_n n}{\log{n}}}.
\end{equation*}
The following lemma adapts \cite[Lemma~6.4 and Lemma~6.5(2)]{Nguyen2022} to the affine Bernoulli columns.

\begin{lemma}
    \label{lem:DenseNormals}
    Let $\calG$ be the event that $\operatorname{span}_{\R}(X_1,\ldots,X_{n_0})^\perp$ contains no nonzero vector supported on at most $\ell$ coordinates. Then, for sufficiently large $n$,
    \begin{equation*}
        \Prob\paren*{\calG^c} \leq \exp\paren*{-\alpha_n n/8},
    \end{equation*}
    uniformly in $A_n$ and $S_n$. 
    Moreover, on $\calG$ and for sufficiently large $n$, every nonzero vector in $W_k(p)^\perp$ has more than $\ell$ nonzero coordinates for every prime $p \geq U_n$ and every $n_0 \leq k \leq n$.
\end{lemma}

\begin{proof}
    If $\calG^c$ occurs, then, for some $1 \leq t \leq \ell$, a set of $t$ rows of $[X_1 \ \cdots \ X_{n_0}]$ is linearly dependent. Fix $S \subseteq [n]$ with $\abs*{S} = t$, and designate one of the rows indexed by $S$. After conditioning on the other $t-1$ rows, their span is a fixed subspace of $\R^{n_0}$ of dimension at most $t-1$, and hence of codimension at least $n_0-t+1$. The designated row is independent of the conditioned rows, and its coordinates are independent and have maximum atom at most $1-\alpha_n$. Therefore, \cref{lem:NonIdenticalOdlyzko} shows that the probability that the designated row belongs to this span is at most
    \begin{equation*}
        (1-\alpha_n)^{n_0-t+1}.
    \end{equation*}
    If the rows indexed by $S$ are linearly dependent, then at least one of them belongs to the span of the others.  A union bound over the choice of this row, over $S$, and over $1 \leq t \leq \ell$ therefore gives
    \begin{equation*}
        \Prob\paren*{\calG^c} \leq \sum_{t=1}^{\ell} \binom{n}{t} t(1-\alpha_n)^{n_0-t+1} \leq \exp\paren*{-\alpha_n n/8}.
    \end{equation*}
    The final estimate follows from the definitions of $n_0$ and $\ell$, the bound $\binom{n}{t} \leq (en/t)^t$, and a sufficiently small choice of $c_1$.

    It remains to prove the asserted consequence on $\calG$. Since $W_{n_0}(p) \subseteq W_k(p)$ for every $n_0 \leq k \leq n$, we have $W_k(p)^\perp \subseteq W_{n_0}(p)^\perp$. It therefore suffices to consider $k = n_0$. Suppose that $\calG$ occurs and that, for some prime $p \geq U_n$, $0 \neq w \in W_{n_0}(p)^\perp$ has support of size at most $\ell$. Let $S = \supp(w)$, with $\abs*{S} = t$, and let $A_S$ be the submatrix of $[X_1 \ \cdots \ X_{n_0}]$ formed by the rows indexed by $S$. For all sufficiently large $n$, we have $t \leq \ell < n_0$. Since $w\vert_S$ is a nonzero vector in the left kernel of $A_S$ over $\mathbb{F}_p$, we have $\rank_{\mathbb{F}_p}(A_S) < t$, and hence every $t \times t$ minor of $A_S$ is divisible by $p$. By Hadamard's inequality and the bound $\abs*{(X_j)_i} \leq n^T + 1$, the absolute value of each such minor is at most $t^{t/2}(n^T + 1)^t$. Since $t\leq\ell$, for all sufficiently large $n$ we have
    \begin{equation*}
        \begin{aligned}
            t^{t/2}(n^T + 1)^t &= \exp\paren*{t\paren*{\tfrac{1}{2}\log{t}+\log(n^T + 1)}} \\
            &\leq \exp\paren*{c_1(T + 1)\alpha_n n} \\
            &< \exp(c_0\alpha_n n/2) = U_n \leq p,
        \end{aligned}
    \end{equation*}
    where the strict inequality follows from $c_1(T + 1) < c_0/2$, and the last inequality follows from $p \geq U_n$. Divisibility by $p$ therefore forces all these minors to vanish over $\Z$. Thus, $\rank_{\R}(A_S) < t$, so these rows admit a nontrivial linear dependence over $\R$. Extending its coefficients by zero outside $S$ produces a nonzero vector in $\operatorname{span}_{\R}(X_1,\ldots,X_{n_0})^\perp$ supported on at most $\ell$ coordinates, contradicting $\calG$. 
\end{proof}

On the event $\calG$, every nonzero normal vector to a relevant late span has large support. The finite-field Littlewood--Offord estimate, \cref{lem:ErdosLittlewoodOfford}, therefore shows that an unexposed column is unlikely to lie in any fixed proper span. This gives the following anti-concentration estimate.

\Needspace{8\baselineskip}
\begin{lemma}[{cf. \cite[Lemma~6.2]{Nguyen2022}}]
    \label{lem:LateAnticoncentration}
    For all sufficiently large $n$, on the event $\calG$, the following bound holds for every prime $p \geq U_n$ and every integer $k$ satisfying $n_0 \leq k \leq n$ and $W_k(p) \neq \mathbb{F}_p^n$:
    \begin{equation}
        \label{eq:GoodSpanBound}
        \Prob\paren*{X_{k+1} \in W_k(p) \mid X_1,\ldots,X_k} = O\paren*{\frac{\sqrt{\log{n}}}{\alpha_n\sqrt{n}}},
    \end{equation}
    uniformly in $p,k,A_n,S_n$.
\end{lemma}

\begin{proof}
    Since $W_k(p)$ is a proper subspace, there exists a nonzero vector $w \in W_k(p)^\perp$. By \cref{lem:DenseNormals}, its support has size greater than $\ell$. After conditioning on $X_1,\ldots,X_k$, the space $W_k(p)$ and the vector $w$ are fixed, while $X_{k+1}$ remains independent of the exposed columns. Set $I \coloneqq \supp(w)$ and $t \coloneqq \abs*{I} > \ell$. Since $w \in W_k(p)^\perp$, membership in $W_k(p)$ implies orthogonality to $w$. Therefore,
    \begin{align*}
        \Prob\paren*{X_{k+1} \in W_k(p) \mid X_1,\ldots,X_k} &\leq \Prob\paren*{\sum_{i \in I} w_i(X_{k+1})_i = 0 \mid X_1,\ldots,X_k} \\
        &= \Prob\paren*{\sum_{i \in I} s_{i,k+1}w_i\xi_{i,k+1} = -\sum_{i \in I} a_{i,k+1}w_i}.
    \end{align*}
    The last equation holds because, conditional on the exposed columns, the vector $w$ is fixed and the variables $\xi_{i,k+1}$ remain independent. Moreover, all coefficients $s_{i,k+1}w_i$, $i \in I$, are nonzero. Since $\alpha_n\ell \to \infty$, we have $\alpha_nt \geq 4$ for all sufficiently large $n$. So we can apply \cref{lem:ErdosLittlewoodOfford}, which gives
    \begin{equation*}
        \Prob\paren*{X_{k+1} \in W_k(p) \mid X_1,\ldots,X_k} \leq \frac{1}{p} + \frac{2}{\sqrt{\alpha_n t}} \leq \frac{1}{p} + \frac{2}{\sqrt{\alpha_n\ell}} = O\paren*{\frac{\sqrt{\log{n}}}{\alpha_n\sqrt{n}}}.
    \end{equation*}
    The term $p^{-1}$ is absorbed into the final bound because $p \geq \exp(c_0\alpha_n n/2)$.
\end{proof}

Combining \cref{lem:LateAnticoncentration} with the watch list bound in \cref{lem:WatchList} controls second failed rank increments after $n_0$.

\begin{lemma}[{cf. \cite[Section~6.1, p.~34]{Nguyen2022}}]
    \label{lem:LateLosses}
    The probability that $\det[X_1 \ \cdots \ X_n] \neq 0$, that $\calG$ occurs, and that no large prime accumulates two failed rank increments by time $n_0$, but some large prime does so afterwards, is $O(\frac{\log^{5/2}n}{\alpha_n^3\sqrt{n}})$, uniformly in the deterministic affine data $A_n$ and $S_n$.
\end{lemma}

\begin{proof}
    Fix $n_0+1 \leq k \leq n+1$, and condition on an arbitrary realization of $X_1,\ldots,X_{k-1}$ for which $\calG$ occurs, $\abs*{\mathfrak{W}_{k-1}} \leq C\frac{\log{n}}{c_0\alpha_n}$, and no large prime has yet suffered two failed rank increments. As in the proof of \cref{lem:EarlyLosses}, we do not condition on $\det[X_1 \ \cdots \ X_n] \neq 0$.
    Every $p \in \mathfrak{W}_{k-1}$ has then suffered exactly one failed rank increment, and hence $\dim W_{k-1}(p) = k-2$. By \cref{lem:LateAnticoncentration}, the conditional probability that $X_k$ causes a second failed rank increment for a fixed $p \in \mathfrak{W}_{k-1}$ is $O\paren*{\frac{\sqrt{\log{n}}}{\alpha_n\sqrt{n}}}$.

    A prime outside $\mathfrak{W}_{k-1}$ can suffer only its first failed rank increment at time $k$. Therefore, uniformly over all relevant $X_1,\ldots,X_{k-1}$, the probability that the first time at which any large prime has a second failed rank increment equals $k$ is at most $\abs*{\mathfrak{W}_{k-1}}O(\frac{\sqrt{\log{n}}}{\alpha_n\sqrt{n}}) \leq O(\frac{\log^{3/2}n}{\alpha_n^2\sqrt{n}})$ by \cref{lem:WatchList}. Since $n+1-n_0 \leq 3\log{n}/\alpha_n + 1$, summing over all $n_0 + 1\leq k \leq n+1$ proves the claimed estimate.
\end{proof}

\begin{proof}[Proof of \cref{prop:ThreePrimeRegimes}, part~\labelcref{item:LargePrimes}]
    Let $\calE_{\mathrm{large}}$ be the event that $M_n$ fails to have full row rank over $\mathbb{F}_p$ for some prime $p \geq U_n$, and let $\calD$ be the event that $\det[X_1 \ \cdots \ X_n] \neq 0$. Let $\calE_{\mathrm{early}}$ be the event that some large prime accumulates two failed rank increments by time $n_0$ and $\calE_{\mathrm{late}}$ be the event that no large prime does so by time $n_0$, but some large prime does so afterwards. To not have full rank over $\mathbb{F}_p$, we require at least two failed rank increments for that prime. Consequently, 
    \begin{equation*}
        \calE_{\mathrm{large}} \subseteq \calD^c \cup (\calD \cap \calE_{\mathrm{early}}) \cup \calG^c \cup (\calD \cap \calG \cap \calE_{\mathrm{late}}).
    \end{equation*}
    Taking probabilities, applying the union bound, and using \cref{lem:WatchList,lem:EarlyLosses,lem:DenseNormals,lem:LateLosses}, respectively, gives
    \begin{equation*}
        \Prob\paren*{\calE_{\mathrm{large}}} \leq O(e^{-c\alpha_n n}) + O(n^{-1/2}) + e^{-\alpha_n n/8} + O\paren*{\frac{\log^{5/2}n}{\alpha_n^3\sqrt n}}.
    \end{equation*}
    Since $\alpha_n = n^{-1/6+\eps}$, we have $\frac{\log^{5/2}(n)}{\alpha_n^3\sqrt{n}} = \frac{\log^{5/2}(n)}{n^{3\eps}} = o(1)$ and all the remaining terms also tend to zero. Thus, $\Prob(\calE_{\mathrm{large}}) = o(1)$, uniformly in $A_n$ and $S_n$, as claimed.
\end{proof}

\section{Proofs of the Average-Case Estimates}
\label[appendix]{app:MomentEstimates}

We prove \cref{lem:FirstMoment,lem:SecondMoment} in the conditional probability model of \cref{sec:AverageCase}.

\subsection{First-Moment Estimate}
\label[appendix]{app:FirstMoment}

We first express the expectation as a sum over subset sizes and then estimate the contributions of small, intermediate, and central subsets.

For even $n \geq 4$ and even $s$ with $2 \leq s \leq n-2$, define
\begin{equation}
    \label{eq:SingleSubsetProbability}
    f_n(s) \coloneqq \binom{s}{s/2}\binom{n-s}{(n-s)/2}\binom{n}{n/2}^{-1}.
\end{equation}
In one coordinate, conditioning on a total sum of zero makes all sign sequences with exactly $n/2$ positive entries equally likely. Under this conditioning, the sum over a fixed $s$-element subset $I$ in this coordinate is zero precisely when $I$ and its complement each contain equally many positive and negative entries. Moreover, the $d$ coordinates remain independent under this conditioning. Thus,
\begin{equation*}
    \Prob_0\paren*{\sum_{i \in I} X_i = 0} = f_n(s)^d.
\end{equation*}
Odd subset sizes cannot contribute. Linearity of expectation therefore gives the exact identity
\begin{equation}
    \label{eq:FirstMomentRepresentation}
    \E_0 Z_n = \sum_{\substack{2 \leq s \leq n-2 \\ s \text{ even}}} \binom{n}{s}f_n(s)^d.
\end{equation}

We will use two elementary bounds on $f_n$. For even $s$ with $2 \leq s \leq n-4$, we have
\begin{equation*}
    \frac{f_n(s+2)}{f_n(s)} = \frac{(s+1)(n-s)}{(s+2)(n-s-1)}.
\end{equation*}
Together with the symmetry $f_n(s) = f_n(n-s)$, this shows that the maximum over even $2 \leq s \leq n-2$ is attained at $s=2$ and $s=n-2$. Thus,
\begin{equation}
    \label{eq:SingleSubsetConstantBound}
    f_n(s) \leq f_n(2) = \frac{n}{2(n-1)} \leq \frac{2}{3}.
\end{equation}
For every even integer $k \geq 2$, we have the elementary central binomial coefficient bounds
\begin{equation*}
    \frac{2^k}{\sqrt{2k}} \leq \binom{k}{k/2} \leq \frac{2^k}{\sqrt{k+1}}.
\end{equation*}
These follow by induction from $2^{-k-2}\binom{k+2}{(k+2)/2} = \frac{k+1}{k+2}2^{-k}\binom{k}{k/2}$. For every even integer $s$ with $2 \leq s \leq n/2$, substituting these bounds into \eqref{eq:SingleSubsetProbability} yields
\begin{equation}
    \label{eq:SingleSubsetDecayBound}
    f_n(s) \leq \sqrt{\frac{2n}{(s+1)(n-s+1)}} \leq \frac{2}{\sqrt{s}}.
\end{equation}

\begin{proof}[Proof of \cref{lem:FirstMoment}]
    Fix $\eta \in (0,1)$ and let $n \geq 2$ be even with $n \leq (1-\eta)d\log{d}/\log{4}$. For $n = 2$, we have $Z_n = 0$. For $4 \leq n \leq d/2$, \eqref{eq:FirstMomentRepresentation} and \eqref{eq:SingleSubsetConstantBound} give
    \begin{equation*}
        \E_0 Z_n \leq 2^n\paren*{\frac{2}{3}}^d \leq \paren*{\frac{8}{9}}^{d/2} = o(1).
    \end{equation*}

    It remains to consider $d/2 < n \leq (1-\eta)d\log{d}/\log{4}$. Uniformly in this range,
    \begin{equation}
        \label{eq:FirstMomentUniformScales}
        \log{n} \geq \log{d} - \log{2}, \qquad \sqrt{n}\log{n} = O\paren*{\sqrt{d}(\log{d})^{3/2}} = o(d).
    \end{equation}
    By symmetry,
    \begin{equation*}
        \E_0 Z_n \leq 2 \sum_{\substack{2 \leq s \leq n/2 \\ s \text{ even}}} \binom{n}{s} f_n(s)^d.
    \end{equation*}
    For all sufficiently large $d$, we split this sum at $s = \sqrt{n}$ and $s = n/10$.

    \emph{Small subsets.}
    For even $s$ with $2 \leq s \leq \sqrt{n}$, we use \eqref{eq:SingleSubsetConstantBound} and $\binom{n}{s}\leq n^s \leq n^{\sqrt{n}}$. Since there are at most $\sqrt{n}$ summands, we obtain
    \begin{equation*}
        \sum_{\substack{2 \leq s \leq \sqrt{n} \\ s \text{ even}}} \binom{n}{s}f_n(s)^d \leq n^{\sqrt{n} + 1/2} \paren*{\frac{2}{3}}^d = \exp\paren*{-d\log(3/2) + o(d)} = o(1).
    \end{equation*}

    \emph{Intermediate subsets.}
    For even $s$ with $\sqrt{n} < s \leq n/10$, \eqref{eq:SingleSubsetDecayBound} gives $f_n(s) \leq 2n^{-1/4}$. The binomial theorem yields
    \begin{equation*}
        \sum_{0 \leq s \leq n/10} \binom{n}{s} \leq 10^{n/10}(11/10)^n \leq (10e)^{n/10} \leq 2^{n/2}.
    \end{equation*}
    Consequently,
    \begin{equation*}
        \sum_{\substack{\sqrt{n} < s \leq n/10 \\ s \text{ even}}} \binom{n}{s}f_n(s)^d \leq2^{n/2}\paren*{2n^{-1/4}}^d \leq \exp\paren*{-\frac{\eta}{4}d\log{d} + O(d)} = o(1).
    \end{equation*}
    Here we used $n\log{2} \leq (1-\eta)d\log{d}/2$ and \eqref{eq:FirstMomentUniformScales}.

    \emph{Central subsets.}
    For even $s$ with $n/10 < s \leq n/2$, \eqref{eq:SingleSubsetDecayBound} gives $f_n(s) \leq 2\sqrt{10/n}$, and therefore
    \begin{equation*}
        \sum_{\substack{n/10 < s \leq n/2 \\ s \text{ even}}} \binom{n}{s}f_n(s)^d \leq 2^n\paren*{2\sqrt{10/n}}^d \leq \exp\paren*{-\frac{\eta}{2}d\log{d} + O(d)} = o(1).
    \end{equation*}
    All error terms above are uniform over the stated range of $n$. Combining the three estimates with the case $n \leq d/2$ proves the lemma.
\end{proof}

\subsection{Second-Moment Estimate}
\label[appendix]{app:SecondMoment}

Grouping pairs of subsets by their sizes and intersection size gives an exact representation with hypergeometric overlap weights.

Let $n \geq 4$ be even. For even $s,t$ with $2 \leq s,t \leq n-2$ and an integer $k$ satisfying $\max\{0,s+t-n\} \leq k \leq \min\{s,t\}$, define
\begin{equation}
    \label{eq:PairProbability}
    q_n(s,t,k) \coloneqq \binom{n}{n/2}^{-1} \sum_{j=0}^{k} \binom{k}{j}\binom{s-k}{s/2-j}\binom{t-k}{t/2-j}\binom{n-s-t+k}{(n-s-t)/2+j}.
\end{equation}
Here, for a nonnegative integer $a$, we set $\binom{a}{b} = 0$ whenever $b \notin \{0,\ldots,a\}$. We also define
\begin{equation}
    \label{eq:OverlapWeights}
    p_n(s,t,k) \coloneqq \binom{s}{k}\binom{n-s}{t-k}\binom{n}{t}^{-1}.
\end{equation}
For a fixed $s$-element subset $I$ and a uniformly chosen $t$-element subset $J$ of $[n]$, this is the probability that $\abs*{I \cap J} = k$.

\begin{lemma}[Second-moment representation]
    \label{lem:SecondMomentRepresentation}
    For every even $n \geq 4$ and every $\gamma \in (0,1/2)$,
    \begin{equation}
        \label{eq:SecondMomentRepresentation}
        \E_0\widehat{Z}_n(\gamma)^2 = \sum_{\substack{\gamma n \leq s,t \leq(1-\gamma)n \\ s,t\text{ even}}} \binom{n}{s}\binom{n}{t} \sum_{k=\max\{0,s+t-n\}}^{\min\{s,t\}} p_n(s,t,k)q_n(s,t,k)^d.
    \end{equation}
\end{lemma}

\begin{proof}
    Fix subsets $I,J \subseteq [n]$ with even sizes $s,t$ and $\abs*{I \cap J} = k$. In a single coordinate under $\Prob_0$, all sign sequences with exactly $n/2$ positive entries are equally likely. If exactly $j$ entries in $I \cap J$ are positive, then both subset sums vanish precisely when the numbers of positive entries in $I \setminus J$, $J \setminus I$, and $[n] \setminus (I \cup J)$ are
    \begin{equation*}
        \frac{s}{2} - j, \qquad \frac{t}{2} - j, \qquad \frac{n-s-t}{2} + j,
    \end{equation*}
    respectively. Summing the corresponding counts over $j$ and dividing by $\binom{n}{n/2}$ gives $q_n(s,t,k)$. Since the coordinates remain independent under the conditioning,
    \begin{equation*}
        \Prob_0\paren*{\sum_{i \in I} X_i = 0, \sum_{i \in J} X_i = 0} = q_n(s,t,k)^d.
    \end{equation*}
    The number of ordered pairs $(I,J)$ with these parameters is
    \begin{equation*}
        \binom{n}{s}\binom{s}{k}\binom{n-s}{t-k} = \binom{n}{s}\binom{n}{t}p_n(s,t,k).
    \end{equation*}
    Expanding $\widehat{Z}_n(\gamma)^2$ as a sum over ordered pairs and taking expectations gives \eqref{eq:SecondMomentRepresentation}. Only even subset sizes contribute, since a sum of an odd number of signs cannot vanish.
\end{proof}

We next compare the joint probability with the product of the individual probabilities. For the admissible parameters above, define
\begin{equation}
    \label{eq:OverlapCorrelation}
    \rho_n(s,t,k) \coloneqq \frac{nk-st}{\sqrt{s(n-s)t(n-t)}}.
\end{equation}
This parameter measures the deviation of the overlap from its mean $st/n$. The following lemma gives a global upper bound and a uniform approximation when all four regions determined by the two subsets have sizes proportional to $n$.

\begin{lemma}[Joint probability estimates]
    \label{lem:PairProbabilityEstimates}
    For all admissible $n,s,t,k$,
    \begin{equation}
        \label{eq:PairGlobalBound}
        \frac{q_n(s,t,k)}{f_n(s)f_n(t)} \leq \sqrt{n}.
    \end{equation}
    Moreover, fix $\delta \in (0,1/4)$. As $n \to \infty$, uniformly over all admissible $s,t,k$ satisfying $\min\{k,s-k,t-k,n-s-t+k\} \geq \delta n$,
    \begin{equation}
        \label{eq:PairLocalApproximation}
        \frac{q_n(s,t,k)}{f_n(s)f_n(t)} = \frac{1+O_\delta(n^{-1})}{\sqrt{1-\rho_n(s,t,k)^2}}.
    \end{equation}
\end{lemma}

\begin{proof}
    The joint probability satisfies $q_n(s,t,k) \leq \min\{f_n(s),f_n(t)\}$. The elementary central binomial coefficient bounds used in \cref{app:FirstMoment} give
    \begin{equation*}
        f_n(s) \geq \frac{\sqrt{n+1}}{2\sqrt{s(n-s)}} \geq n^{-1/2},
    \end{equation*}
    and the same holds for $f_n(t)$. This proves \eqref{eq:PairGlobalBound}.

    For the local approximation, write
    \begin{equation*}
        r_1 = k, \qquad r_2 = s-k, \qquad r_3 = t-k, \qquad r_4 = n-s-t+k,
    \end{equation*}
    and assume $r_i \geq \delta n$ for every $i$. Since $n,s,t$ are even, all four $r_i$ have the same parity. For a positive integer $r$ and $u \in \Z - r/2$, set
    \begin{equation*}
        b_r(u) \coloneqq 2^{-r}\binom{r}{r/2+u}.
    \end{equation*}
    With $\Lambda = \Z - k/2$, a change of index and binomial symmetry in \eqref{eq:PairProbability} yield
    \begin{equation*}
        q_n(s,t,k) = b_n(0)^{-1} \sum_{u \in \Lambda} \prod_{i=1}^{4} b_{r_i}(u).
    \end{equation*}
    Thus, $\Lambda$ is either the integer lattice or the half-integer lattice, both with spacing one.

    Stirling's formula and a Taylor expansion about $r/2$ give, uniformly for $\delta n \leq r \leq n$ and $u \in \Z - r/2$ with $\abs*{u} \leq n^{3/5}$,
    \begin{equation*}
        b_r(u) = \sqrt{\frac{2}{\pi r}}e^{-2u^2/r}\paren*{1 + O_\delta\paren*{\frac{1}{n} + \frac{u^4}{n^3}}}.
    \end{equation*}
    Set $H = \sum_{i=1}^{4} r_i^{-1}$. Since $\sum_i r_i = n$, we have $16/n \leq H \leq 4/(\delta n)$. Multiplying the four expansions gives
    \begin{equation*}
        \prod_{i=1}^{4} b_{r_i}(u) = \frac{4e^{-2Hu^2}}{\pi^2\sqrt{r_1r_2r_3r_4}}\paren*{1 + O_\delta\paren*{\frac{1}{n} + \frac{u^4}{n^3}}}.
    \end{equation*}
    To sum these errors, we use the Gaussian lattice estimates
    \begin{equation*}
        \sum_{u \in \Lambda} e^{-2Hu^2} = \sqrt{\frac{\pi}{2H}}\paren*{1 + O_\delta(n^{-1})}, \qquad \sum_{u \in \Lambda} u^4e^{-2Hu^2} = O_\delta(n^2) \sum_{u \in \Lambda} e^{-2Hu^2}.
    \end{equation*}
    These follow from the midpoint rule on the intervals $[u-1/2,u+1/2]$: for a twice continuously differentiable, rapidly decaying function $g$, its error is at most $\frac{1}{8} \int_{\R} \abs*{g''(x)} \,dx$. Apply this to $g(x) = e^{-2Hx^2}$ and $g(x) = x^4e^{-2Hx^2}$. The fourth-moment estimate shows that the accumulated relative error is $O_\delta(n^{-1} + n^2/n^3) = O_\delta(n^{-1})$.

    The omitted terms are negligible. Indeed, since $b_r(u) \leq 1$, the binomial Chernoff bound gives
    \begin{equation*}
        \sum_{\substack{u \in \Lambda \\ \abs*{u} > n^{3/5}}} \prod_{i=1}^{4} b_{r_i}(u) \leq \sum_{\substack{u \in \Lambda \\ \abs*{u} > n^{3/5}}} b_{r_1}(u) \leq 2e^{-2n^{1/5}}.
    \end{equation*}
    The corresponding Gaussian tails are exponentially small as well. Combining these estimates with $b_n(0) = \sqrt{2/(\pi n)}\paren*{1 + O(n^{-1})}$ yields
    \begin{equation*}
        q_n(s,t,k) = \frac{2}{\pi} \sqrt{\frac{n}{H r_1r_2r_3r_4}}\paren*{1 + O_\delta(n^{-1})}.
    \end{equation*}
    Finally, direct expansion gives
    \begin{equation*}
        H r_1r_2r_3r_4 = \frac{s(n-s)t(n-t)-(nk-st)^2}{n},
    \end{equation*}
    while Stirling's formula gives, uniformly in the stated range,
    \begin{equation*}
        f_n(s)f_n(t) = \frac{2n}{\pi\sqrt{s(n-s)t(n-t)}}\paren*{1 + O_\delta(n^{-1})}.
    \end{equation*}
    Dividing the two expressions proves \eqref{eq:PairLocalApproximation}.
\end{proof}

We now average the normalized joint probability over the overlap distribution. For even $n \geq 4$ and even $s,t$ with $2 \leq s,t \leq n-2$, define
\begin{equation}
    \label{eq:AveragedCorrelation}
    R_{d,n}(s,t) \coloneqq \sum_{k=\max\{0,s+t-n\}}^{\min\{s,t\}} p_n(s,t,k)\paren*{\frac{q_n(s,t,k)}{f_n(s)f_n(t)}}^d.
\end{equation}

For fixed $n,s,t$, let $K$ have probability mass function $\Prob(K = k) = p_n(s,t,k)$. Thus, $K$ is hypergeometric with mean $\mu \coloneqq st/n$. Probabilities and expectations involving $K$ below refer to this auxiliary sampling experiment. Hoeffding's inequality for sampling without replacement \cite[Theorem~1 and Section~6, Theorem~4]{Hoeffding1963}, applied to both tails, gives for every $u \geq 0$,
\begin{equation}
    \label{eq:OverlapConcentration}
    \Prob\paren*{\abs*{K - \mu} \geq u} \leq 2e^{-2u^2/t} \leq 2e^{-2u^2/n}.
\end{equation}
In particular, for every $0 \leq \lambda < 2$,
\begin{equation}
    \label{eq:OverlapExponentialMoment}
    \E\exp\paren*{\lambda\frac{(K - \mu)^2}{n}} \leq 1 + \frac{2\lambda}{2-\lambda}.
\end{equation}
Indeed, with $Y = (K - \mu)^2/n$, integration of the tail bound yields
\begin{equation*}
    \E e^{\lambda Y} = 1 + \lambda \int_0^\infty e^{\lambda y}\Prob(Y > y) \,dy \leq 1 + 2\lambda \int_0^\infty e^{-(2-\lambda)y} \,dy = 1 + \frac{2\lambda}{2-\lambda}.
\end{equation*}
Both bounds hold uniformly over all admissible $n,s,t$.

We will use these estimates to prove the following lemma, which is the remaining ingredient in the proof of \cref{lem:SecondMoment}.

\begin{lemma}[Averaged correlation estimate]
    \label{lem:AveragedCorrelation}
    For every fixed $\eta \in (0,1)$, there exists $\gamma = \gamma(\eta) \in (0,1/2)$ such that
    \begin{equation*}
        \lim_{d \to \infty} R_{d,n}(s,t) = 1,
    \end{equation*}
    uniformly over all even integers $n \geq 4$ and even integers $s,t$ satisfying
    \begin{equation*}
        n \geq (1+\eta)\frac{d\log{d}}{\log{4}}, \qquad \gamma n \leq s,t \leq (1-\gamma)n.
    \end{equation*}
\end{lemma}

\begin{proof}
    We split the overlap sum into a regular range, where the four regions $I \cap J$, $I \setminus J$, $J \setminus I$, and $[n] \setminus (I\cup J)$ each have size at least a small fixed fraction of $n$, and the remaining boundary range. The local approximation and exponential-moment bound control regular overlaps; an entropy bound makes the boundary contribution negligible.

    \emph{Parameter choice.}
    Fix $\eta \in (0,1)$. For $0 < x < 1$, let $h(x) = -x\log{x} - (1-x)\log(1-x)$ denote the binary entropy. Choose $\delta \in (0,1/16)$ sufficiently small and set $\gamma = 1/2 - \delta$ so that
    \begin{equation}
        \label{eq:EntropyGap}
        h(\gamma) - h(4\delta)> \frac{\log{2}}{1 + \eta/2}.
    \end{equation}
    Such a choice is possible because the left-hand side tends to $\log{2}$ as $\delta \downarrow 0$. This entropy gap will be used to control the boundary overlaps. Our choice also ensures $\delta < \gamma^2/2$.

    The lower bound on $n$ gives $d/n = o(1)$ uniformly in $n$. Moreover, since $x \mapsto (\log{x})/x$ is decreasing for $x > e$, evaluating it at $(1 + \eta)d\log{d}/\log{4}$ gives
    \begin{equation}
        \label{eq:UniformDimensionBound}
        \frac{d\log{n}}{2n} \leq \frac{\log{2}}{1 + \eta}\paren*{1 + o(1)} < \frac{\log{2}}{1 + \eta/2}
    \end{equation}
    for all sufficiently large $d$, uniformly over the stated range of $n$.

    \emph{Regular overlaps.}
    Fix admissible $n,s,t$, and let $K$ and $\mu = st/n$ be as above. Define
    \begin{equation*}
        \calK \coloneqq \set*{k \in \Z : \min\{k,s-k,t-k,n-s-t+k\} \geq \delta n}.
    \end{equation*}
    Thus, $K \in \calK$ is the event that the random overlap is regular. Let $R^{\mathrm{reg}}_{d,n}(s,t)$ denote the contribution of $k \in \calK$ to \eqref{eq:AveragedCorrelation}. For $k \in \calK$, the determinant identity used in the proof of \cref{lem:PairProbabilityEstimates} and the bounds on $s,t$ give
    \begin{equation*}
        1 - \rho_n(s,t,k)^2 \geq 64\delta^3, \qquad \rho_n(s,t,k)^2 \leq \frac{(k - \mu)^2}{\gamma^2(1 - \gamma)^2n^2}.
    \end{equation*}
    Taking logarithms in \eqref{eq:PairLocalApproximation} and using $-\log(1-x) \leq x/(1-x)$ for $0 \leq x < 1$, we obtain a constant $C = C(\eta) > 0$ such that
    \begin{equation}
        \label{eq:RegularCorrelationBounds}
        e^{-Cd/n} \leq \paren*{\frac{q_n(s,t,k)}{f_n(s)f_n(t)}}^d \leq \exp\paren*{\frac{Cd}{n} + \frac{Cd(k - \mu)^2}{n^2}}
    \end{equation}
    for all sufficiently large $d$, uniformly over the stated parameters and $k \in \calK$.

    Set $\lambda = Cd/n = o(1)$. Averaging the upper bound and applying \eqref{eq:OverlapExponentialMoment} gives
    \begin{equation*}
        R^{\mathrm{reg}}_{d,n}(s,t) \leq e^\lambda \E\exp\paren*{\lambda\frac{(K - \mu)^2}{n}} \leq e^\lambda\paren*{1 + \frac{2\lambda}{2-\lambda}} = 1 + o(1).
    \end{equation*}

    For the lower bound, each of the four expressions defining $\calK$ is at least $\gamma^2n$ at $k = \mu$. Since $\delta < \gamma^2/2$, the event $\abs*{K - \mu} \leq \gamma^2n/2$ implies $K \in \calK$. Thus, \eqref{eq:OverlapConcentration} yields
    \begin{equation*}
        \Prob\paren*{K \in \calK} \geq 1 - 2e^{-\gamma^4n/2}.
    \end{equation*}
    Together with the lower bound in \eqref{eq:RegularCorrelationBounds}, this gives
    \begin{equation*}
        R^{\mathrm{reg}}_{d,n}(s,t) \geq e^{-\lambda}\paren*{1 - 2e^{-\gamma^4n/2}} = 1 - o(1).
    \end{equation*}
    Consequently,
    \begin{equation}
        \label{eq:RegularOverlapContribution}
        R^{\mathrm{reg}}_{d,n}(s,t) = 1 + o(1),
    \end{equation}
    uniformly over the stated range of $n,s,t$. It remains to estimate the contribution of $k \notin \calK$.

    \emph{Boundary overlaps.}
    Fix $I \subseteq [n]$ with $\abs*{I} = s$, and choose $J$ uniformly among all $t$-element subsets of $[n]$. Then $K = \abs*{I \cap J}$ has probability mass function $p_n(s,t,k)$. Since $\gamma = 1/2 - \delta$, the bounds on $s,t$ imply
    \begin{equation*}
        \abs*{s - t} \leq 2\delta n, \qquad \abs*{n - s - t} \leq 2\delta n.
    \end{equation*}
    Writing $I^c = [n] \setminus I$ and using symmetric differences, we have
    \begin{equation*}
        \abs*{I \triangle J} = 2\min\{s-K,t-K\} + \abs*{s-t}, \qquad \abs*{I^c \triangle J} = 2\min\{K,n-s-t+K\} + \abs*{n-s-t}.
    \end{equation*}
    Thus, $K \notin \calK$ implies that either $\abs*{I \triangle J} < 4\delta n$ or $\abs*{I^c \triangle J} < 4\delta n$. There are exactly $\binom{n}{j}$ subsets whose symmetric difference with a given subset of $[n]$ has size $j$. Consequently,
    \begin{equation}
        \label{eq:BoundaryOverlapProbability}
        \Prob(K \notin \calK) \leq \frac{2\sum_{j=0}^{\lfloor4\delta n\rfloor} \binom{n}{j}}{\binom{n}{t}} \leq 2(n+1)\exp\paren*{-n\paren*{h(\gamma) - h(4\delta)}}.
    \end{equation}
    For the second inequality, we used $h(t/n) \geq h(\gamma)$ and the elementary binomial estimates
    \begin{equation*}
        \sum_{j=0}^{\lfloor4\delta n\rfloor} \binom{n}{j} \leq e^{nh(4\delta)}, \qquad \binom{n}{t} \geq \frac{e^{nh(t/n)}}{n+1}.
    \end{equation*}
    For the first estimate, each factor $(4\delta)^j(1 - 4\delta)^{n-j}$ in the indicated range is at least $e^{-nh(4\delta)}$, so the estimate follows from the binomial theorem. For the second, $t$ is a mode of the binomial distribution with parameters $n$ and $t/n$, and its probability is therefore at least $1/(n+1)$.

    By \eqref{eq:EntropyGap}, the constant
    \begin{equation*}
        c_\eta \coloneqq h(\gamma) - h(4\delta) - \frac{\log{2}}{1+\eta/2}
    \end{equation*}
    is positive. Combining \eqref{eq:BoundaryOverlapProbability} with the global bound \eqref{eq:PairGlobalBound} and \eqref{eq:UniformDimensionBound} gives, for all sufficiently large $d$,
    \begin{equation*}
        \begin{aligned}
            0 \leq R_{d,n}(s,t) - R^{\mathrm{reg}}_{d,n}(s,t) &\leq n^{d/2}\Prob(K \notin \calK) \\
            &\leq 2(n+1)\exp\paren*{\frac{d}{2}\log{n} - n\paren*{h(\gamma) - h(4\delta)}} \\
            &\leq 2(n+1)e^{-c_\eta n} = o(1),
        \end{aligned}
    \end{equation*}
    uniformly over the stated range of $n,s,t$. Together with \eqref{eq:RegularOverlapContribution}, this yields $R_{d,n}(s,t) = 1 + o(1)$ uniformly and completes the proof.
\end{proof}

\begin{proof}[Proof of \cref{lem:SecondMoment}]
    Fix $\eta \in (0,1)$ and choose $\gamma$ as in \cref{lem:AveragedCorrelation}. For sufficiently large $d$ and even $n \geq (1+\eta)d\log{d}/\log{4}$, \eqref{eq:SecondMomentRepresentation} and \eqref{eq:AveragedCorrelation} give
    \begin{equation*}
        \begin{aligned}
            \E_0\widehat{Z}_n(\gamma)^2 &= \sum_{\substack{\gamma n \leq s,t \leq (1-\gamma)n \\ s,t \text{ even}}} \binom{n}{s}\binom{n}{t} f_n(s)^d f_n(t)^d R_{d,n}(s,t) \\
            &= (1 + o(1))\Bigg(\sum_{\substack{\gamma n \leq s \leq (1-\gamma)n \\ s \text{ even}}} \binom{n}{s} f_n(s)^d\Bigg)^2 \\
            &= (1 + o(1))\paren*{\E_0\widehat{Z}_n(\gamma)}^2.
        \end{aligned}
    \end{equation*}
    The second equality follows from \cref{lem:AveragedCorrelation}, since its error is uniform in $n,s,t$ and all summation weights are nonnegative. The last equality follows by linearity of expectation. The resulting $o(1)$ is uniform over the stated range of $n$. For sufficiently large $d$, the summation range contains an even integer, so $\E_0\widehat{Z}_n(\gamma) > 0$. Dividing by its square proves the claim.
\end{proof}

\subsection{Rademacher Sums}
\label[appendix]{app:RademacherSums}

We justify \cref{rem:RademacherSums} for fixed $m\geq1$ and admissible $n\geq2$, meaning $mn\in2\Z$. Here the conditional expectations and counting variables refer to the law $\nu_{m,d}$. For proper subset sizes $s,t$ with $ms,mt$ even and overlap $k$, the one-coordinate probabilities are
\begin{equation*}
    F_{m,n}(s)=f_{mn}(ms), \qquad Q_{m,n}(s,t,k)=q_{mn}(ms,mt,mk).
\end{equation*}
Other subset sizes contribute zero. The counting weights remain $\binom ns$ and $p_n(s,t,k)$, since subsets select whole vector blocks.

For contributing sizes, $ms$ and $m(n-s)$ are even and at least two. Thus \eqref{eq:SingleSubsetConstantBound} and \eqref{eq:SingleSubsetDecayBound} give
\begin{equation*}
    F_{m,n}(s)\leq\frac23, \qquad F_{m,n}(s)\leq\frac{2}{\sqrt{ms}}\leq\frac{2}{\sqrt{s}}\quad(s\leq n/2).
\end{equation*}
Consequently, $\E_0 Z_n\leq2^n(2/3)^d=o(1)$ uniformly for admissible $2\leq n\leq d/2$, including $n=2$. The same three-range argument handles the remaining subcritical lengths, now summing over $1\leq s<n$ with $ms$ even. This includes $s=1$ and odd $n$ when $m$ is even.

The correlation \eqref{eq:OverlapCorrelation} is unchanged under $(n,s,t,k)\mapsto(mn,ms,mt,mk)$. On regular overlaps, the local error in \eqref{eq:PairLocalApproximation} becomes $O_\delta((mn)^{-1})$. The four scaled region sizes have the same parity, so the summation lattice $\Z-mk/2$ still has spacing one. Hence the regular-overlap argument carries over.

The global bound \eqref{eq:PairGlobalBound} becomes $Q_{m,n}(s,t,k)/(F_{m,n}(s)F_{m,n}(t))\leq\sqrt{mn}$, adding only the factor $m^{d/2}$ to the boundary estimate. For fixed $m$, its exponent $(d/2)\log m=o(n)$ uniformly for $n\geq(1+\eta)d\log d/\log4$, and the positive entropy gap absorbs this factor.

Both moment estimates therefore persist. The passage to $D_\eps(\nu_{m,d})$ is unchanged, since admissible lengths have spacing at most two.

\Needspace{8\baselineskip}
\section{Proofs for Lattice Balls and Sparse Ternary Cubes}
\label[appendix]{app:LatticeBalls}

We prove \cref{lem:EuclideanCircuitBound,lem:TernaryLowerBound} and the upper bound in \cref{lem:SpecialCase}, using the notation of \cref{sec:LatticeBalls}.

\subsection{Euclidean Circuit Bound}
\label[appendix]{app:EuclideanCircuitBound}

We bound positive definite integer Gram determinants and apply the estimate to circuit cofactors.

For an integer $k \geq 3$, put
\begin{equation*}
    \lambda_k \coloneqq \rho_k^2, \qquad \kappa_k \coloneqq 
    \frac{k}{\lambda_k(1 - \lambda_k^{-4})^2}.
\end{equation*}
Thus $\lambda_k + \lambda_k^{-1} = k$. For a symmetric matrix $G \in \R^{n \times n}$, its Gram graph has vertex set $\{1,\ldots,n\}$ and an edge between distinct vertices $i$ and $j$ whenever $G_{ij} \neq 0$.

\begin{lemma}[Gram determinant bound]
    \label{lem:GramDeterminantBound}
    Let $k \geq 3$ be an integer and let $G \in \Z^{n \times n}$ be symmetric positive definite with $G_{ii} \leq k$ for every $i$. If its Gram graph has $c$ connected components, then
    \begin{equation*}
        \det G \leq \kappa_k^c \lambda_k^n.
    \end{equation*}
\end{lemma}

\begin{proof}
    After simultaneously permuting rows and columns, the connected components form diagonal blocks. It therefore suffices to prove the assertion when the Gram graph is connected. Write $\lambda \coloneqq \lambda_k$ and $\kappa \coloneqq \kappa_k$.

    Choose a rooted spanning tree and order its vertices so that every parent precedes its children. Represent $G$ as the Gram matrix of linearly independent vectors $v_1,\ldots,v_n$, in this order. Let $\ell_i$ be the squared distance of $v_i$ from $\operatorname{span}\{v_1,\ldots,v_{i-1}\}$. Gram--Schmidt orthogonalization gives
    \begin{equation*}
        \det G = \prod_{i=1}^n \ell_i.
    \end{equation*}
    Call a pivot large if $\ell_i > \lambda$.

    \emph{Earlier neighbors.}
    Suppose that $v_i$ has two earlier neighbors. Let $H$ be their $2 \times 2$ Gram matrix and let $b$ contain their inner products with $v_i$. Independence and integrality imply that the off-diagonal entry of $H$ has absolute value at most $k - 1$. Hence the largest eigenvalue of $H$ is at most $2k - 1$, while $\norm*{b}_2^2 \geq 2$. Projection onto these two vectors gives
    \begin{equation*}
        \ell_i \leq k - b^\top H^{-1} b \leq k-\frac{2}{2k-1} \leq \lambda.
    \end{equation*}
    The last inequality follows from $\lambda > k - \tfrac{1}{2}$.

    Consequently, a large nonroot pivot has exactly one earlier neighbor, its parent $u$. Let $E$ be the span of all preceding vectors other than $v_u$, and put
    \begin{equation*}
        s_u \coloneqq \norm*{v_u - \operatorname{proj}_E v_u}_2^2.
    \end{equation*}
    Since $v_i$ is orthogonal to $E$ and $s_u \leq \ell_u$,
    \begin{equation*}
        \ell_i = \norm*{v_i}_2^2 - \frac{\inp{v_i}{v_u}^2}{s_u} \leq k - \frac{1}{\ell_u}.
    \end{equation*}
    For $f(t) \coloneqq k - t^{-1}$, we have $f(\lambda) = \lambda$. Thus, the parent of a large pivot is large as well.

    \emph{Branching of large pivots.}
    The large children of a vertex $u$ are pairwise orthogonal, and each is orthogonal to every vector preceding $v_u$. Suppose that $m$ such children precede another large child $v_i$. Projecting the residual of $v_u$ onto these $m$ orthogonal vectors reduces its squared norm by at least $m/k$. In the preceding notation, this gives
    \begin{equation*}
        0 < s_u \leq \ell_u - \frac{m}{k}, \qquad \ell_i \leq k - \frac{1}{\ell_u - m/k}.
    \end{equation*}
    If $u$ is not the root, it has an earlier neighbor, so $\ell_u \leq k - 1/k$. A second large child would therefore imply
    \begin{equation*}
        \ell_i \leq f(k - 2/k) < \lambda,
    \end{equation*}
    a contradiction. At the root, where $\ell_u\le k$, the same argument excludes a third large child. Here we used $k - 2/k < \lambda$ and the monotonicity of $f$.

    Thus, if any large pivots occur, they form the root and at most two chains descending from it. A pivot at distance $j$ along either chain is at most $f^j(k)$.

    \emph{Product estimate.}
    Define
    \begin{equation*}
        P_0 = 1, \qquad P_1 = k, \qquad P_{j+1} = kP_j - P_{j-1}.
    \end{equation*}
    Since $\lambda + \lambda^{-1} = k$,
    \begin{equation*}
        P_j = \frac{\lambda^{j+1} - \lambda^{-j-1}}{\lambda - \lambda^{-1}}, \qquad f^j(k) = \frac{P_{j+1}}{P_j}.
    \end{equation*}
    Consequently, for every integer $t \geq 0$,
    \begin{equation*}
        \prod_{j=1}^{t}\frac{f^j(k)}{\lambda} = \frac{P_{t+1}}{k\lambda^t} = \frac{1 - \lambda^{-2t - 4}}{1 - \lambda^{-4}} \leq \frac{1}{1 - \lambda^{-4}}.
    \end{equation*}
    All remaining pivots are at most $\lambda$. Accounting for the root and the two possible chains therefore yields
    \begin{equation*}
        \det G \leq \lambda^n\frac{k}{\lambda}\paren*{\frac{1}{1 - \lambda^{-4}}}^2 = \kappa\lambda^n.
    \end{equation*}
    If there are no large pivots, the same bound follows directly from $\det G \leq \lambda^n$ and $\kappa > 1$. Multiplying over the connected components proves the lemma.
\end{proof}

\begin{proof}[Proof of \cref{lem:EuclideanCircuitBound}]
    Choose $h$ linearly independent coordinate rows of $[a_1 \ \cdots \ a_{h+1}]$, and denote the resulting matrix by
    \begin{equation*}
        B = [b_1 \ \cdots \ b_{h+1}] \in \Z^{h \times (h+1)}.
    \end{equation*}
    This restriction preserves the kernel and hence the primitive positive relation $g$. Moreover, $\norm*{b_j}_2^2 \leq k$ for every $j$.

    For $1 \leq j \leq h + 1$, let $B^{(j)}$ be obtained by deleting column $j$, and put
    \begin{equation*}
        \Delta_j \coloneqq (-1)^{j-1} \det B^{(j)}, \qquad \gamma \coloneqq \gcd\{\abs*{\Delta_1},\ldots,\abs*{\Delta_{h+1}}\}.
    \end{equation*}
    Every $\Delta_j$ is nonzero because the columns form a circuit. The cofactor formula gives
    \begin{equation*}
        g_j = \frac{\abs*{\Delta_j}}{\gamma}, \qquad \norm*{g}_1 = \frac{1}{\gamma} \sum_{j=1}^{h+1} \abs*{\Delta_j}.
    \end{equation*}

    \emph{Components after deleting a column.}
    The Gram graph of the entire circuit is connected. Indeed, if it were disconnected, the weighted sums over its components would lie in mutually orthogonal subspaces. Since their total sum is zero, each component sum would be zero, contradicting minimal dependence.

    Fix $j$, let $T \coloneqq B^{(j)}$, and denote the number of connected components of the Gram graph of $T^\top T$ by $c$. Each component contains a column $v_\ell$ with $\inp{b_j}{v_\ell} \neq 0$; otherwise it would also be a component of the entire circuit graph. The selected vectors $v_1,\ldots,v_c$ are pairwise orthogonal. Their inner products with $b_j$ are nonzero integers, so Bessel's inequality gives
    \begin{equation*}
        \frac{c}{k} \leq \sum_{\ell=1}^{c} \frac{\inp{b_j}{v_\ell}^2}{\norm*{v_\ell}_2^2} \leq \norm*{b_j}_2^2 \leq k.
    \end{equation*}
    Thus $c \leq k^2$. \cref{lem:GramDeterminantBound} now yields
    \begin{equation*}
        \abs*{\det T} = \sqrt{\det(T^\top T)} \leq \kappa_k^{c/2}\lambda_k^{h/2} \leq \kappa_k^{k^2/2}\rho_k^h.
    \end{equation*}

    \emph{The uniform constant.}
    Write $\lambda \coloneqq \lambda_k$ and $x \coloneqq \lambda^{-2}$. Since $k \geq 3$, we have $0 < x < 1/6$, and
    \begin{equation*}
        k^2 = \frac{(1+x)^2}{x}, \qquad \kappa_k = \frac{1}{(1-x)^2(1+x)}.
    \end{equation*}
    Using $-\log(1-x) \leq x/(1-x)$ and $\log(1+x) \geq x/(1+x)$, we obtain
    \begin{equation*}
        \frac{k^2}{2}\log{\kappa_k} = \frac{(1+x)^2}{2x}\paren*{-2\log(1-x) - \log(1+x)} \leq \frac{(1+x)(1+3x)}{2(1-x)} \leq \frac{21}{20} < \log{3}.
    \end{equation*}
    For the final inequality, note that
    \begin{equation*}
        \log{3} = 2 \int_0^{1/2} \frac{dt}{1-t^2} \geq 2 \int_0^{1/2}(1+t^2) \,dt = \frac{13}{12} > \frac{21}{20}.
    \end{equation*}
    Hence $\abs*{\Delta_j} \leq 3 \rho_k^h$ for every $j$. Summing the cofactors and using $\gamma \geq 1$ proves
    \begin{equation*}
        \norm*{g}_1 \leq \frac{3(h+1)\rho_k^h}{\gamma} \leq 3(h+1)\rho_k^h,
    \end{equation*}
    as claimed.
\end{proof}

\subsection{Sparse Ternary Construction}
\label[appendix]{app:TernaryConstruction}

We prove \cref{lem:TernaryLowerBound} by adjoining a unit column to an integer matrix with a large determinant so that the maximal minors have greatest common divisor one, preserving the ternary alphabet and support bound.

\begin{proof}[Proof of \cref{lem:TernaryLowerBound}]
    Fix an integer $k \geq 3$ and put $t \coloneqq k - 2$ and $b \coloneqq 2^t$. For every integer $m \geq 1$, we construct an irreducible zero-sum sequence in dimension $n \coloneqq 2bm$.

    \emph{A sparse skew-symmetric block.}
    Index the standard basis of $\R^b$ by the subsets $I \subseteq \{0,\ldots,t-1\}$, with the empty set labeled $0$ and $\{0\}$ labeled $1$. For $0 \leq j < t$, define
    \begin{equation*}
        E_j e_I \coloneqq (-1)^{\abs*{I \cap \{0,\ldots,j-1\}} + \one_{\{j \in I\}}} e_{I \triangle \{j\}}, \qquad S \coloneqq \sum_{j=0}^{t-1} E_j.
    \end{equation*}
    The indicator $\one_{\{j \in I\}}$ equals $1$ if $j \in I$ and $0$ otherwise; $I_b$ denotes the $b \times b$ identity matrix. Here $\triangle$ denotes symmetric difference, and $\{0,\ldots,j-1\}$ is empty when $j = 0$. Each $E_j$ is a signed permutation matrix. Toggling $j$ twice changes the sign once, whereas exchanging two distinct toggles changes their combined sign once. Consequently,
    \begin{equation*}
        E_j^2 = -I_b, \qquad E_j^\top = -E_j, \qquad E_jE_\ell = -E_\ell E_j \quad (j \neq \ell).
    \end{equation*}
    It follows that $S^\top = -S$ and $S^2 = -tI_b$. For each $I$, the vectors $E_je_I$ occupy distinct coordinates, so every column of $S$ is ternary with exactly $t$ nonzero entries. Moreover, $e_1^\top Se_0 = 1$.

    To preserve the determinant of the square matrix while making a maximal minor of the augmented matrix coprime to it, we will modify one off-diagonal block. Define a nilpotent matrix $N$ by
    \begin{equation*}
        Ne_0 = 0, \qquad Ne_i = e_{i+1} \quad (1 \leq i < b - 1), \qquad Ne_{b-1} = e_0.
    \end{equation*}
    Thus $N^b = 0$, its diagonal is zero, and each column has at most one nonzero entry. The only basis vector mapped to $e_0$ by $N^{b-1}$ is $e_1$, whence
    \begin{equation}
        \label{eq:TernaryTerminalEntry}
        e_0^\top N^{b-1}Se_0 = 1.
    \end{equation}

    \emph{The block matrix and its column supports.}
    Divide $\R^n$ into $2m$ blocks of size $b$, indexed by $0,\ldots,2m-1$, and let $U_i$ embed $\R^b$ into block $i$. Let $B_0$ be the block tridiagonal matrix with blocks
    \begin{equation*}
        (B_0)_{ij} = \begin{cases}
            S, & i = j, \\
            I_b, & \abs*{i - j} = 1, \\
            0, & \text{otherwise}.
        \end{cases}
    \end{equation*}
    Set
    \begin{equation*}
        B \coloneqq B_0 + U_1NU_0^\top, \qquad f \coloneqq U_0e_0.
    \end{equation*}
    Thus the first lower off-diagonal block is changed from $I_b$ to $I_b + N$. Endpoint columns of $B_0$ have $t + 1$ nonzero entries, while interior columns have $t + 2 = k$. The modification adds at most one entry to each column of the first block. Since $N$ has zero diagonal, these entries do not overlap with those of $I_b$. Hence all columns of $[B \mid -f]$ belong to $\calT_n(k)$. Every column of $B$ has at least $t + 1 \geq 2$ nonzero entries, whereas $f$ is a unit vector.

    \emph{The determinant and the required inverse blocks.}
    Define positive integers $u_j,v_j$ by
    \begin{equation*}
        u_0 = v_0 = 1, \qquad u_j = u_{j-1} + tv_{j-1}, \qquad v_j = u_j + v_{j-1} \quad (j \geq 1).
    \end{equation*}
    Their transition matrix is
    \begin{equation*}
        \begin{pmatrix}
            1 & t \\
            1 & t + 1
        \end{pmatrix},
    \end{equation*}
    which has determinant one and trace $k$. In particular,
    \begin{equation*}
        u_1 = k - 1, \qquad u_j = ku_{j-1} - u_{j-2} \quad (j \geq 2).
    \end{equation*}
    The unimodular transition preserves $\gcd(u_j,v_j) = 1$, and the recurrence for $u_j$ preserves coprimality of consecutive terms. Therefore, writing $u \coloneqq u_m$, $w \coloneqq u_{m-1}$, and $v \coloneqq v_{m-1}$, we have
    \begin{equation}
        \label{eq:TernaryCoprimalities}
        \gcd(u,w) = \gcd(u,v) = 1.
    \end{equation}

    To compute the block determinant, put
    \begin{equation*}
        F_0 = I_b, \qquad F_1 = S, \qquad F_j = SF_{j-1} - F_{j-2} \quad (j \geq 2).
    \end{equation*}
    The identity $S^2 = -tI_b$ and the scalar recurrences give
    \begin{equation*}
        F_{2j} = (-1)^j u_jI_b, \qquad F_{2j-1} = (-1)^{j-1}v_{j-1}S \quad (j \geq 1).
    \end{equation*}
    All these matrices are invertible and commute. Successive block elimination in $B_0$ has diagonal blocks $F_jF_{j-1}^{-1}$, for $1 \leq j \leq 2m$. Their determinants telescope, and $b$ is even, so
    \begin{equation*}
        \det B_0 = \det F_{2m} = u^b.
    \end{equation*}
    The first block row of $B_0^{-1}$ is given by $(B_0^{-1})_{0j} = (-1)^jF_{2m-1-j}F_{2m}^{-1}$; this follows directly by multiplication with $B_0$ and the recurrence for $F_j$. In particular,
    \begin{equation}
        \label{eq:TernaryInverseBlocks}
        (B_0^{-1})_{00} = -\frac{v}{u} S, \qquad (B_0^{-1})_{01} = \frac{w}{u} I_b.
    \end{equation}
    The determinant identity for a finite-rank update yields
    \begin{equation*}
        \begin{aligned}
            \Delta \coloneqq \det B &= \det B_0\,\det\!\paren*{I_b + U_0^\top B_0^{-1}U_1N} \\
            &= u^b\det\!\paren*{I_b + \frac{w}{u} N} = u^b,
        \end{aligned}
    \end{equation*}
    because $N$ is nilpotent.

    \emph{Primitivity of the augmented cofactor relation.}
    Let $z \coloneqq \operatorname{adj}(B)f \in \Z^n$, so $Bz = \Delta f$. To compute its first block, put $y \coloneqq U_0^\top B^{-1}f$. From $B_0B^{-1}f = f - U_1Ny$ and \eqref{eq:TernaryInverseBlocks}, we obtain
    \begin{equation*}
        y = -\frac{v}{u} Se_0 - \frac{w}{u} Ny.
    \end{equation*}
    Expanding the inverse of $I_b + (w/u)N$ as a finite geometric series therefore gives
    \begin{equation*}
        U_0^\top z = -v\sum_{j=0}^{b-1} (-w)^j u^{b-1-j}N^jSe_0.
    \end{equation*}
    For the coordinate $z_* \coloneqq e_0^\top U_0^\top z$, \eqref{eq:TernaryTerminalEntry} implies
    \begin{equation*}
        z_* \equiv (-1)^bvw^{b-1} \pmod{u}.
    \end{equation*}
    By \eqref{eq:TernaryCoprimalities}, this coordinate is coprime to $u$ and hence to $\Delta = u^b$. The maximal minors of $[B \mid -f]$ are, up to signs, the coordinates of $z$ and $\Delta$. Their greatest common divisor is therefore one.

    Write $B_j$ for column $j$ of $B$ and let $J \coloneqq \{j : z_j \neq 0\}$. For $j \in J$, use the vector $\sgn(z_j)B_j$ with multiplicity $\abs*{z_j}$, and use $-f$ with multiplicity $\Delta$. These vectors are still in $\calT_n(k)$ and have total sum zero. They are distinct: the columns of $B$ are linearly independent, and each has more than one nonzero coordinate. The kernel of $[B \mid -f]$ is one-dimensional, so after the sign changes and deletion of zero-cofactor columns its cone of nonnegative kernel vectors is the ray generated by
    \begin{equation*}
        g \coloneqq \paren*{(\abs*{z_j})_{j \in J},\Delta}.
    \end{equation*}
    This vector is primitive. Any integer multiplicity vector $0 \leq x \leq g$ describing a zero-sum subsequence must satisfy $x = \theta g$ for some real $\theta$. B\'ezout's identity and $\gcd(g) = 1$ imply $\theta \in \Z$, while the multiplicity bounds give $0 \leq \theta \leq 1$. Thus we obtain an irreducible zero-sum sequence of length at least $\Delta = u_m^b$.

    \emph{Passage to arbitrary dimensions.}
    For $d \geq 2b$, take $m \coloneqq \floor*{d/(2b)}$ and append $d - 2bm$ zero coordinates to every constructed vector. This preserves all zero-sum relations and the support bound, and hence
    \begin{equation*}
        \Dvp(\calT_d(k)) \geq u_{\floor*{d/(2b)}}^{\,b}.
    \end{equation*}
    Write $\lambda \coloneqq \lambda_k = \rho_k^2$. The recurrence and initial values of $u_j$ give
    \begin{equation*}
        u_j = \frac{\lambda^{j+1} + \lambda^{-j}}{\lambda + 1}, \qquad \log u_j = j\log\lambda + O_k(1).
    \end{equation*}
    Since $b$ is fixed and $bm/d \to 1/2$, the preceding lower bound implies
    \begin{equation*}
        \liminf_{d \to \infty} \Dvp(\calT_d(k))^{1/d} \geq \sqrt{\lambda} = \rho_k.
    \end{equation*}
\end{proof}

\subsection{The \texorpdfstring{$p = 1$}{p = 1} Case}
\label[appendix]{app:SpecialCase}

For the upper bound in \cref{lem:SpecialCase}, we first treat primitive columns and then account for arbitrary column gcds. The lower bound was proved in \cref{sec:LatticeBalls}.

\begin{lemma}[Circuit bound for primitive columns]
    \label{lem:PrimitiveL1CircuitBound}
    Let $M \geq 3$ be an integer, put $P \coloneqq M - 1$ and $c \coloneqq 1 + P^{-2}$, and let $a_1,\ldots,a_{h+1} \in \Z^h$ form a positive circuit of rank $h \geq 1$. Suppose that every $a_i$ is primitive and satisfies $\norm*{a_i}_1 \leq M$. Then its primitive positive integer relation $g$ satisfies
    \begin{equation*}
        \norm*{g}_1 \leq c^M(h+1)P^h.
    \end{equation*}
\end{lemma}

\begin{proof}
    Fix a column $a_*$ of the circuit and let $T$ be the square matrix formed by the remaining $h$ columns. We show that $\abs*{\det T} \leq c^MP^h$.

    \emph{Large columns and their graph.}
    A primitive column of $T$ cannot contain an entry of magnitude $M$: the norm bound would force it to be a nonprimitive multiple of a unit vector. If it contains an entry of magnitude $P = M - 1$, primitivity and the norm bound force exactly one further nonzero entry, of magnitude one. Call such a column large. All other columns have entries of magnitude at most $P - 1$.

    Represent each large column by an edge between its two nonzero coordinate rows, retaining its coefficients in the matrix. This defines an undirected multigraph on the $h$ rows; parallel edges are allowed, but loops are not. In a connected component with $v$ vertices and $e$ edges, connectivity gives $e \geq v - 1$, while independence of the columns gives $e \leq v$. Thus every component is a tree or has exactly one cycle. Isolated rows are included as trees with no edges.

    Let $a$ be the number of large columns and let $q \coloneqq h - a$. Write $v_1,\ldots,v_q$ for the remaining columns and let $\nu$ be the number of unicyclic components. Summing $v - e$ over the components shows that there are exactly $q$ tree components.

    \emph{Determinant reduction.}
    For a tree component $C$ with $a_C$ edges, let $H_C$ be its $(a_C + 1) \times a_C$ large-column matrix and let $w_C$ be its signed row-deletion cofactor vector. For any chosen root, leaf elimination expresses the corresponding cofactor as a product of $a_C$ entries, each of magnitude $1$ or $P$. Hence all entries of $w_C$ have magnitudes $P^j$ with $0 \leq j \leq a_C$. At most one entry has magnitude $P^{a_C}$: rooting the tree at two different vertices reverses the assignment on every edge of their connecting path, and the two endpoint coefficients of each such edge have magnitudes $P$ and $1$. Consequently, the normalized row $P^{-a_C}w_C^\top$ has all coefficients of magnitude at most one and at most one of magnitude greater than $1/P$. For an isolated row we use $a_C = 0$ and $w_C = (1)$.

    For a unicyclic component with $a_C$ edges, successively remove its attached leaves. The remaining cycle has length $\ell \geq 2$ and exactly two determinant terms, of magnitudes $P^j$ and $P^{\ell-j}$ for some $0 \leq j \leq \ell$. Since $P^j + P^{\ell-j} \leq P^\ell + 1$, its full determinant has magnitude at most
    \begin{equation*}
        P^{a_C-\ell}(P^\ell + 1) \leq cP^{a_C}.
    \end{equation*}
    The block is nonsingular because its columns belong to the basis $T$.

    Expand $\det T$ first along the square unicyclic blocks and then along the large columns of each tree. Exactly one row per tree remains for the other columns. By multilinearity, the resulting tree-cofactor sums form a $q \times q$ matrix $F$ with entries
    \begin{equation*}
        F_{Cj} \coloneqq P^{-a_C}w_C^\top(v_j|_C),
    \end{equation*}
    where $C$ ranges over the tree components. Therefore
    \begin{equation}
        \label{eq:L1DeterminantReduction}
        \abs*{\det T} \leq P^a c^\nu\abs*{\det F}.
    \end{equation}
    For $q = 0$, the empty determinant is one. For $q > 0$, the exact expansion and $\det T \neq 0$ imply that $F$ is nonsingular. This reduction is used only to evaluate the determinant; it does not replace the original integer columns or their multiplicity constraints.

    \emph{Norms of the reduced columns.}
    Let $\calU$ be the union of the row sets of the unicyclic components and put
    \begin{equation*}
        b_j \coloneqq \sum_{i \in \calU} \abs*{(v_j)_i} \qquad (1 \leq j \leq q).
    \end{equation*}
    Since the $j$th column $F_j$ is nonzero, $v_j$ has a nonzero integer entry on a tree row. Thus $0 \leq b_j \leq P$, and its total absolute coordinate sum on tree rows is at most $P + 1 - b_j$. We will show that $\nu \leq M + \sum_j b_j$. To absorb the factors $c^{b_j}$ in the determinant bound, we aim for $c^{b_j}\norm*{F_j}_2 \leq P$ for every $j$. We claim that
    \begin{equation}
        \label{eq:L1ReducedColumnNorms}
        \norm*{F_j}_2 \leq \begin{cases}
            P, & b_j = 0, \\
            P - 1 + P^{-1}, & b_j = 1, \\
            P + 1 - b_j, & 2 \leq b_j \leq P.
        \end{cases}
    \end{equation}
    The last bound follows from the triangle inequality, since each normalized tree coefficient has magnitude at most one.

    For the first two bounds, replace entries and coefficients by their absolute values, which can only increase the norm of the resulting nonnegative vector. For $P \geq 3$, maximize this convex norm over
    \begin{equation*}
        0 \leq x_i \leq P - 1, \qquad \sum_i x_i \leq H,
    \end{equation*}
    with $H = P + 1$ or $H = P$. At an extreme point, all but at most one nonzero coordinate are equal to $P - 1$. With two nonzero coordinates, their magnitudes are $(P - 1,2)$ for $H = P + 1$ and $(P - 1,1)$ for $H = P$. The remaining vertices have at most one nonzero coordinate and satisfy both claimed bounds directly. If the two entries lie in distinct trees, their reduced norms are at most
    \begin{equation*}
        \sqrt{(P - 1)^2 + 4} \leq P, \qquad \sqrt{(P - 1)^2 + 1} \leq P - 1 + P^{-1},
    \end{equation*}
    respectively. If they lie in the same tree, at most one coefficient exceeds $1/P$, so the corresponding bounds are $P - 1 + 2/P \leq P$ and $P - 1 + 1/P$. The larger entry can be assigned the larger coefficient when maximizing either expression.

    For $P = 2$, the entries have magnitude at most one. Three unit entries distributed among trees give the bounds $2$, $\sqrt{(3/2)^2 + 1}$, or $\sqrt{3}$, according to whether they occupy one, two, or three trees; all are at most $2$. Two unit entries give at most $3/2$ in one tree and $\sqrt{2} \leq 3/2$ in two trees. This proves \eqref{eq:L1ReducedColumnNorms} in every case.

    \emph{Controlling the number of unicyclic components.}
    Every unicyclic component is met by a nonzero entry of $a_*$ or of at least one of $v_1,\ldots,v_q$. Otherwise the positive circuit relation restricted to its rows would be a nontrivial relation among the columns of its nonsingular square block. Each such occurrence contributes at least one to the corresponding absolute coordinate sum, and the components have disjoint row sets. Hence
    \begin{equation*}
        \nu \leq \norm*{a_*}_1 + \sum_{j=1}^q b_j \leq M + \sum_{j=1}^q b_j.
    \end{equation*}
    For $P \geq 2$, direct calculation gives
    \begin{equation*}
        c(P - 1 + P^{-1}) \leq P, \qquad c^2(P - 1) \leq P.
    \end{equation*}
    Moreover, $c^b(P + 1 - b)$ is nonincreasing for integers $2 \leq b \leq P$, since the ratio of consecutive terms is $c(P - b)/(P + 1 - b) \leq 1$. Together with \eqref{eq:L1ReducedColumnNorms}, this yields $c^{b_j}\norm*{F_j}_2 \leq P$ for every $j$. Hadamard's inequality and \eqref{eq:L1DeterminantReduction} now give
    \begin{equation*}
        \abs*{\det T} \leq P^a c^{M+\sum_j b_j} \prod_{j=1}^q\norm*{F_j}_2 \leq c^MP^{a+q} = c^MP^h.
    \end{equation*}
    This argument also covers $q = 0$, with empty products equal to one. Since the deleted column was arbitrary, every maximal minor of the circuit matrix has this bound. The primitive positive relation consists of their absolute values divided by their common gcd, which is at least one. Summing the $h + 1$ cofactors proves the lemma.
\end{proof}

\begin{proof}[Proof of the upper bound in \cref{lem:SpecialCase}]
    Fix $M \geq 3$, put $P \coloneqq M - 1$, $c \coloneqq 1 + P^{-2}$, and let
    \begin{equation*}
        \Lambda_M \coloneqq \operatorname{lcm}(1,\ldots,M).
    \end{equation*}
    Consider a positive circuit of rank $h \geq 1$ whose original columns belong to $\calB_{1,d}(M)$, and let $g$ be its primitive positive integer relation. Retain $h$ linearly independent original coordinate rows. This preserves the kernel and decreases column norms. Write the resulting nonzero columns as
    \begin{equation*}
        a_i = c_i u_i \quad (1 \leq i \leq h + 1),
    \end{equation*}
    where $c_i$ is the positive gcd of the coordinates of $a_i$ and $u_i \in \Z^h$ is primitive. Then $1 \leq c_i \leq M$ and $\norm*{u_i}_1 \leq M$. The primitive columns $u_i$ still form a positive circuit; let $z$ be their primitive positive integer relation.

    \emph{Recovering the original multiplicities.}
    An integer relation $n$ among the $a_i$ satisfies $c_i n_i = \tau z_i$ for a common real scalar $\tau$. Since $\gcd(z) = 1$, B\'ezout's identity forces $\tau \in \Z$. All coordinates $n_i = \tau z_i/c_i$ are integral precisely when $c_i/\gcd(c_i,z_i)$ divides $\tau$ for every $i$. Consequently, the primitive positive relation $g$ is given by
    \begin{equation*}
        g_i = \frac{\tau_0z_i}{c_i}, \qquad \tau_0 \coloneqq \operatorname{lcm}_{1 \leq i \leq h + 1} \frac{c_i}{\gcd(c_i,z_i)}.
    \end{equation*}
    Each factor in this least common multiple divides some $c_i \leq M$, so $\tau_0 \mid \Lambda_M$. By \cref{lem:PrimitiveL1CircuitBound},
    \begin{equation*}
        \norm*{g}_1 \leq \Lambda_M\norm*{z}_1 \leq \Lambda_M c^M(h+1)P^h \leq \Lambda_M c^M(d+1)P^d.
    \end{equation*}

    \Needspace{12\baselineskip}
    \emph{Transfer to arbitrary irreducible zero-sum sequences.}
    The alphabet-size estimate in \cref{sec:LatticeBalls} gives
    \begin{equation*}
        \abs*{\calB_{1,d}(M)} = O_M(d^M).
    \end{equation*}
    An irreducible zero-sum sequence containing zero has length one. For any other irreducible zero-sum sequence, apply \cref{lem:CircuitTransfer} with the preceding uniform circuit bound and $s - r \leq s \leq \abs*{\calB_{1,d}(M)}$. It follows that
    \begin{equation*}
        \Dvp(\calB_{1,d}(M)) \leq C_M d^{M+1}(M-1)^d
    \end{equation*}
    for every $d \geq 1$, after choosing $C_M > 0$ depending only on $M$. This is the required upper bound.
\end{proof}

\section*{Acknowledgements}
The authors used OpenAI's ChatGPT and Codex, with models from the GPT-5 and GPT-6 families, to assist with literature research, mathematical exploration, proof development and checking, computational experiments, and manuscript preparation. The authors retain full responsibility for the content of the manuscript, including all mathematical claims and references.

\bibliographystyle{alpha}
\bibliography{References}

@article{Hunkenschroeder2026,
  author  = {Christoph Hunkenschr{\"o}der and
             Kim-Manuel Klein and
             Martin Kouteck{\'y} and
             Alexandra Lassota and
             Asaf Levin},
  title   = {Tight lower bounds for block-structured integer programs},
  journal = {Mathematical Programming},
  volume  = {216},
  pages   = {521--538},
  year    = {2026},
  doi     = {10.1007/s10107-025-02296-z}
}

@article{Grinberg1980,
  author  = {Grinberg, V. S. and Sevastyanov, S. V.},
  title   = {The value of the {Steinitz} constant},
  journal = {Functional Analysis and Its Applications},
  volume  = {14},
  year    = {1980},
  pages   = {125--126},
  doi     = {10.1007/BF01086559}
}

@article{Steinitz1913,
  author  = {Steinitz, Ernst},
  title   = {Bedingt konvergente Reihen und konvexe Systeme},
  journal = {Journal f{\"u}r die reine und angewandte Mathematik},
  volume  = {143},
  year    = {1913},
  pages   = {128--176},
  doi     = {10.1515/crll.1913.143.128}
}

@article{DeLoera2008,
  author  = {Jes{\'u}s A. De Loera and
             Raymond Hemmecke and
             Shmuel Onn and
             Robert Weismantel},
  title   = {$N$-fold Integer Programming},
  journal = {Discrete Optimization},
  volume  = {5},
  number  = {2},
  pages   = {231--241},
  year    = {2008},
  doi     = {10.1016/j.disopt.2006.06.006}
}

@article{Eisenbrand2019,
  author  = {Friedrich Eisenbrand and
             Christoph Hunkenschr{\"o}der and
             Kim-Manuel Klein and
             Martin Kouteck{\'y} and
             Asaf Levin and
             Shmuel Onn},
  title   = {An Algorithmic Theory of Integer Programming},
  journal = {CoRR},
  volume  = {abs/1904.01361},
  year    = {2019},
  eprint  = {1904.01361v3},
  archivePrefix = {arXiv},
  primaryClass   = {math.OC},
  note    = {\href{https://arxiv.org/abs/1904.01361v3}{Version~3}}
}

@article{Kozlov1997,
	author  = {Kozlov, Dmitry N. and Vu, Van H.},
	title   = {Coins and Cones},
	journal = {Journal of Combinatorial Theory, Series A},
	volume  = {78},
	number  = {1},
	pages   = {1--14},
	year    = {1997},
	doi     = {10.1006/jcta.1996.2724}
}

@article{Olson1969,
	author  = {Olson, John E.},
	title   = {A Combinatorial Problem on Finite Abelian Groups, I},
	journal = {Journal of Number Theory},
	volume  = {1},
	number  = {1},
	pages   = {8--10},
	year    = {1969},
	doi     = {10.1016/0022-314X(69)90021-3}
}

@misc{Brent2013,
  author        = {Brent, Richard P. and Osborn, Judy-anne H. and
                   Smith, Warren D.},
  title         = {Lower Bounds on Maximal Determinants of
                   {$\{\pm1\}$}-Matrices via the Probabilistic Method},
  year          = {2013},
  eprint        = {1211.3248},
  archivePrefix = {arXiv},
  primaryClass  = {math.CO}
}

@article{Nguyen2020,
  author  = {Nguyen, Hoi H. and Paquette, Elliot},
  title   = {Surjectivity of Near-Square Random Matrices},
  journal = {Combinatorics, Probability and Computing},
  volume  = {29},
  number  = {2},
  pages   = {267--292},
  year    = {2020},
  doi     = {10.1017/S0963548319000348}
}

@article{Nguyen2022,
  author  = {Nguyen, Hoi H. and Wood, Melanie Matchett},
  title   = {Random integral matrices: universality of surjectivity and the cokernel},
  journal = {Inventiones Mathematicae},
  volume  = {228},
  number  = {1},
  pages   = {1--76},
  year    = {2022},
  doi     = {10.1007/s00222-021-01082-w}
}

@article{Tropp2012,
  author  = {Tropp, Joel A.},
  title   = {User-Friendly Tail Bounds for Sums of Random Matrices},
  journal = {Foundations of Computational Mathematics},
  volume  = {12},
  number  = {4},
  pages   = {389--434},
  year    = {2012},
  doi     = {10.1007/s10208-011-9099-z}
}

@article{Hoeffding1963,
  author  = {Hoeffding, Wassily},
  title   = {Probability Inequalities for Sums of Bounded Random Variables},
  journal = {Journal of the American Statistical Association},
  volume  = {58},
  number  = {301},
  pages   = {13--30},
  year    = {1963},
  doi     = {10.1080/01621459.1963.10500830}
}

@book{Schrijver1986,
	author    = {Alexander Schrijver},
	title     = {Theory of Linear and Integer Programming},
	publisher = {John Wiley \& Sons},
	year      = {1986}
}

@article{Wood2019,
	author        = {Wood, Melanie Matchett},
	title         = {Random integral matrices and the {Cohen--Lenstra} heuristics},
	journal       = {American Journal of Mathematics},
	year          = {2019},
	volume        = {141},
	number        = {2},
	pages         = {383--398}
}

@book{Dummit2004,
	author    = {Dummit, David S. and Foote, Richard M.},
	title     = {Abstract Algebra},
	edition   = {3},
	publisher = {John Wiley \& Sons},
	address   = {Hoboken, NJ},
	year      = {2004}
}

@article{Kahn1995,
	author  = {Kahn, Jeff and Koml{\'o}s, J{\'a}nos and Szemer{\'e}di, Endre},
	title   = {On the probability that a random {$\pm 1$}-matrix is singular},
	journal = {Journal of the American Mathematical Society},
	year    = {1995},
	volume  = {8},
	number  = {1},
	pages   = {223--240}
}

@article{Odlyzko1988,
	author  = {Odlyzko, Andrew M.},
	title   = {On subspaces spanned by random selections of {$\pm 1$} vectors},
	journal = {Journal of Combinatorial Theory, Series A},
	year    = {1988},
	volume  = {47},
	number  = {1},
	pages   = {124--133}
}

@misc{Maples2010,
	author        = {Maples, Kenneth},
	title         = {Singularity of Random Matrices over Finite Fields},
	year          = {2010},
	url           = {https://arxiv.org/abs/1012.2372}
}

@book{Terras1999,
	author    = {Terras, Audrey},
	title     = {Fourier Analysis on Finite Groups and Applications},
	series    = {London Mathematical Society Student Texts},
	volume    = {43},
	publisher = {Cambridge University Press},
	address   = {Cambridge},
	year      = {1999}
}

@article{Paat2020,
  author  = {Paat, Joseph and Weismantel, Robert and Weltge, Stefan},
  title   = {Distances between optimal solutions of mixed-integer programs},
  journal = {Mathematical Programming},
  volume  = {179},
  number  = {1--2},
  pages   = {455--468},
  year    = {2020},
  doi     = {10.1007/s10107-018-1323-z}
}

@article{Paat2022,
  author  = {Paat, Joseph and Schl{\"o}ter, Miriam and Weismantel, Robert},
  title   = {The integrality number of an integer program},
  journal = {Mathematical Programming},
  volume  = {192},
  pages   = {271--291},
  year    = {2022},
  doi     = {10.1007/s10107-021-01651-0}
}

@article{Bader2018,
  author  = {Bader, J{\"o}rg and Hildebrand, Robert and Weismantel, Robert and
             Zenklusen, Rico},
  title   = {Mixed integer reformulations of integer programs and the affine {TU}-dimension of a matrix},
  journal = {Mathematical Programming},
  volume  = {169},
  pages   = {565--584},
  year    = {2018},
  doi     = {10.1007/s10107-017-1147-2}
}

@article{Gao2006,
  author  = {Weidong Gao and Alfred Geroldinger},
  title   = {Zero-sum problems in finite abelian groups: A survey},
  journal = {Expositiones Mathematicae},
  volume  = {24},
  number  = {4},
  pages   = {337--369},
  year    = {2006},
  doi     = {10.1016/j.exmath.2006.07.002}
}

@article{Baeth2014,
  author  = {Nicholas R. Baeth and Alfred Geroldinger},
  title   = {Monoids of modules and arithmetic of direct-sum decompositions},
  journal = {Pacific Journal of Mathematics},
  volume  = {271},
  number  = {2},
  pages   = {257--319},
  year    = {2014},
  doi     = {10.2140/pjm.2014.271.257}
}

@book{Onn2010,
  author    = {Shmuel Onn},
  title     = {Nonlinear Discrete Optimization: An Algorithmic Theory},
  series    = {Zurich Lectures in Advanced Mathematics},
  publisher = {European Mathematical Society},
  year      = {2010},
  doi       = {10.4171/093}
}

@incollection{Ekbatani2022,
	author    = {Ekbatani, Farbod and Natura, Bento and V{\'e}gh, L{\'a}szl{\'o} A.},
	title     = {Circuit Imbalance Measures and Linear Programming},
	booktitle = {Surveys in Combinatorics 2022},
	editor    = {Nixon, Anthony and Prendiville, Sean},
	series    = {London Mathematical Society Lecture Note Series},
	pages     = {64--114},
	publisher = {Cambridge University Press},
	year      = {2022},
	doi       = {10.1017/9781009093927.004}
}

@article{Brianski2024,
  author  = {Marcin Bria{\'n}ski and Martin Kouteck{\'y} and Daniel Kr{\'a}l' and
             Krist{\'y}na Pek{\'a}rkov{\'a} and Felix Schr{\"o}der},
  title   = {Characterization of matrices with bounded {Graver} bases and depth
             parameters and applications to integer programming},
  journal = {Mathematical Programming},
  volume  = {208},
  pages   = {497--531},
  year    = {2024},
  doi     = {10.1007/s10107-023-02048-x}
}

@inproceedings{Alon1996,
  author    = {Noga Alon and Dmitry N. Kozlov and Van H. Vu},
  title     = {The geometry of coin-weighing problems},
  booktitle = {Proceedings of the 37th Annual Symposium on Foundations of Computer Science},
  pages     = {524--532},
  publisher = {IEEE},
  year      = {1996}
}

@article{Plagne2015,
  author        = {Alain Plagne and Salvatore Tringali},
  title         = {The {Davenport} constant of a box},
  journal       = {Acta Arithmetica},
  volume        = {171},
  number        = {3},
  pages         = {197--219},
  year          = {2015},
  doi           = {10.4064/aa171-3-1},
  eprint        = {1405.4363},
  archivePrefix = {arXiv}
}

@article{Alon1997,
  author  = {Noga Alon and Van H. V{\~u}},
  title   = {Anti-{Hadamard} matrices, coin weighing, threshold gates, and
             indecomposable hypergraphs},
  journal = {Journal of Combinatorial Theory, Series A},
  volume  = {79},
  number  = {1},
  pages   = {133--160},
  year    = {1997},
  doi     = {10.1006/jcta.1997.2780}
}

@article{Buchin2012,
  author  = {Kevin Buchin and Ji{\v r}{\'i} Matou{\v s}ek and Robin A. Moser and
             D{\"o}m{\"o}t{\"o}r P{\'a}lv{\"o}lgyi},
  title   = {Vectors in a box},
  journal = {Mathematical Programming},
  volume  = {135},
  pages   = {323--335},
  year    = {2012},
  doi     = {10.1007/s10107-011-0474-y}
}

@misc{OpenAI2026a,
  author       = {{OpenAI}},
  title        = {The {Euclidean Steinitz--Bergstr\"om} theorem},
  year         = {2026},
  howpublished = {\href{https://github.com/openai/math/tree/main/preprints/The-Euclidean-Steinitz-Bergstrom-theorem-September-24-2026}{OpenAI mathematics repository}},
  note         = {Preprint, September 24, 2026, Theorem 1.2; used only as a conditional input here}
}

@article{Balachandran2021,
  author  = {Niranjan Balachandran and Eshita Mazumdar},
  title   = {Zero sums in restricted sequences},
  journal = {Discrete Mathematics},
  volume  = {344},
  number  = {7},
  pages   = {112394},
  year    = {2021},
  doi     = {10.1016/j.disc.2021.112394}
}

@article{Bushaw2024,
  author  = {Neal Bushaw and Glenn Hurlbert},
  title   = {Thresholds for zero-sums with small cross numbers in abelian groups},
  journal = {Integers},
  volume  = {24},
  pages   = {A93},
  year    = {2024},
  url     = {https://math.colgate.edu/~integers/y93/y93.pdf}
}

@misc{Girard2025,
  author        = {Benjamin Girard and Alain Plagne},
  title         = {The {Davenport} constant of balls and boxes},
  year          = {2025},
  eprint        = {2510.20412},
  archivePrefix = {arXiv},
  primaryClass  = {math.NT},
  note          = {arXiv:\allowbreak2510.20412v2. Revised 8 October 2026}
}

@misc{Girard2026,
  author        = {Benjamin Girard and Alain Plagne},
  title         = {The {Davenport} constant of low-dimensional boxes},
  year          = {2026},
  eprint        = {2610.16562},
  archivePrefix = {arXiv},
  primaryClass  = {math.NT},
  note          = {arXiv:\allowbreak2610.16562v1. 12 October 2026}
}

@article{Eisenbrand2020,
  author  = {Friedrich Eisenbrand and Robert Weismantel},
  title   = {Proximity Results and Faster Algorithms for Integer Programming
             Using the {Steinitz} Lemma},
  journal = {ACM Transactions on Algorithms},
  volume  = {16},
  number  = {1},
  pages   = {5:1--5:14},
  year    = {2020},
  doi     = {10.1145/3340322}
}

\end{document}